 \documentclass[11pt]{article}
 \usepackage{fullpage}
\date{}

\usepackage{tikz}
\usetikzlibrary{calc}
\usetikzlibrary{matrix}
\usepackage{amsmath,amssymb,amsfonts,amsthm,subfigure}
\usepackage{color}                    
\usepackage{hyperref}                 

\usepackage{graphicx}
\usepackage{multirow}
\allowdisplaybreaks

\DeclareMathOperator{\tr}{tr}

\newcommand{\eps}{{\epsilon}}
\newcommand{\bey}{\begin{eqnarray}}
\newcommand{\eey}{\end{eqnarray}}

\newcommand{\beq}{\begin{equation}}
\newcommand{\eeq}{\end{equation}}
\theoremstyle{plain}

\theoremstyle{definition}

\theoremstyle{remark}

\title{Moving Mesh Finite Element Simulation for Phase-Field Modeling of Brittle Fracture and Convergence of Newton's Iteration}

\author{Fei Zhang%
\thanks{College of Petroleum Engineering,  China University of Petroleum -- Beijing, 18 Fuxue Road, Changping,
Beijing 102249, China ({\it fzhang\_cup@outlook.com})},
\and Weizhang Huang%
\thanks{Department of Mathematics, the University of Kansas, Lawrence, Kansas 66049, U.S.A.
({\it whuang@ku.edu})},
\and Xianping Li%
\thanks{Department of Mathematics and Statistics, University of Missouri -- Kansas City,
5120 Rockhill Road, Kansas City, Missouri 64110, U.S.A. ({\it lixianp@umkc.edu})},
\and Shicheng Zhang%
\thanks{College of Petroleum Engineering,  China University of Petroleum -- Beijing, 18 Fuxue Road, Changping,
Beijing 102249, China ({\it zhangsc@cup.edu.cn})}
}
\begin{document}
\vskip 1cm
\maketitle
\begin{abstract}
A moving mesh finite element method is studied for the numerical solution of a phase-field model for brittle fracture.
The moving mesh partial differential equation approach is employed to dynamically track crack propagation. 
Meanwhile, the decomposition of the strain tensor into tensile and compressive components is essential for
the success of the phase-field modeling of brittle fracture but results in a non-smooth elastic energy and stronger
nonlinearity in the governing equation. This makes the governing equation much more difficult to solve and, in particular,
Newton's iteration often to fail to converge.
Three regularization methods are proposed to smooth out the decomposition of the strain tensor.
Numerical examples of fracture propagation under quasi-static load demonstrate that all of the methods can
effectively improve the convergence of Newton's iteration for relatively small values of the regularization parameter
but without comprising the accuracy of the numerical solution. They also show that the moving
mesh finite element method is able to adaptively concentrate the mesh elements around propagating cracks
and handle multiple and complex crack systems.
\end{abstract}

\noindent{\bf AMS 2010 Mathematics Subject Classification.} 65M50, 65M60, 74B99

\noindent{\bf Key Words.}
brittle fracture, phase-field model, Newton's iteration, moving mesh, mesh adaptation,
finite element method

\noindent{\bf Abbreviated title.}
Moving Mesh Finite Element Simulation of Brittle Fracture

\section{Introduction}

Brittle fracture is the fracture of a metallic object or other elastic material where plastic deformation is strongly limited.
It usually occurs very rapidly and can be catastrophic in engineering practice;
e.g., see Pokluda and \v{S}andera \cite{PS10}. Understanding the initiation and propagation of brittle fracture
and preventing fracture failure are vital to the engineering design, where numerical simulation of fracture processes
has become a powerful tool. Computational approaches for studying brittle fracture can be roughly categorized
into two groups, discrete crack models and smeared crack models. In the former group, discontinuous fields
are introduced into the numerical model and cracks are described as moving boundaries. One major challenge
for those models is to track moving boundaries. A commonly used strategy is to change the mesh geometry
by introducing new boundaries at each time step together with adaptive remeshing;
e.g., see \cite{AK10, BBCT00, KAM08, PD04}. The mesh regenerating and boundary updating
not only increase computational cost but also further complicate the implementation of boundary conditions.
In order to avoid complex remeshing, Mo\"{e}s et al. \cite{MDB99, MB02} propose
the extended finite element method, which enriches the finite element spaces with discontinuous fields
based on the partition-of-unity concept and allows the propagation of cracks along element interfaces.
A drawback for the method is that it requires explicit description of crack patterns and thus has difficulty
in dealing with complex cracks and unforeseen patterns of crack propagation.  

In the second group of computational approaches, smeared crack models approximate cracks with continuous fields
and do not rely on explicit description of cracks. The phase-field model based on the variational approach
proposed by Francfort and Marigo \cite{FM98} is a commonly used type of smeared crack model.
In the phase-field modeling, a phase-field variable $d$, which depends on a parameter $l$ describing the actual
width of the smeared cracks, is introduced to indicate where the material is damaged.
One of the major advantages of this model is that the initiation and propagation of cracks are completely
determined by a coupled system of partial differential equations based on the energy functional.
Another advantage is that the generation and propagation of fracture networks do not require
explicitly tracking fracture interfaces. The phase-field modeling is used in this work.

The phase-field modeling has been successfully applied in many other fields including image segmentation \cite{AT90},
dendritic crystal growth \cite{Kob93,WMS93}, and multiple-fluid hydrodynamics \cite{LiuShen2003,Shen2014,Shen2015c,LiuShen2006}.
Since its first application in brittle fracture simulation by Bourdin et al. \cite{BFM00},
significant progress has been made in this area; e.g., see
\cite{AMM09,Bor12,BHLV14,KM10,LLZ16,MVB15,MHW10,MWH10,NYZBC15,VMBV14}.
However, there still exist challenges. In particular, the strain tensor has to be decomposed along eigen-directions
into tensile and compressive components in the presence of cracks, with only the former component
contributing to generation and propagation of cracks. This decomposition of the strain tensor
is introduced by Miehe et al. \cite{MHW10} to account for the reduction of the stiffness of the elastic solid
by cracks and to rule out unrealistic branching. Unfortunately, it also makes the elastic energy non-smooth
and increases the nonlinearity of the governing equation.
As a consequence, the Jacobian matrix of the governing equation may not exist at places and
Newton's iteration can often fail to converge \cite{LLZ16, NYZBC15}.

The phase-field modeling is governed by a coupled system for the phase-field variable $d$
and the displacement $u$ which can be solved in the monolithic or staggered approach.
In the monolithic approach (e.g., see \cite{GL16,Wic17a,Wic17b}), the system is solved simultaneously
for both $d$ and $u$, and it is often challenging to obtain a convergent solution due to
the non-convexity nature of the energy functional.
A damped Newton method with linear search \cite{GL16}  and an error-oriented
Newton method \cite{Wic17a} have been developed to overcome convergence issues
while a modified Newton scheme with Jacobian modification has been proposed by Wick \cite{Wic17b}.
The staggered approach solves the coupled system sequentially for $d$ and $u$
and has been used by a number of researchers; e.g., see \cite{AGL15, AMM09, FM17, MHW10, MWH10}.
By fixing one variable such as the phase-field $d$, the underlying problem becomes convex in the other
unknown variable $u$. Since $d$ and $u$ are not coupled, the procedure becomes simpler and more robust.
The main disadvantage of this approach is that at a loading step many staggered iterations are required
to reach convergence, which can be costly.
The effects of the number of staggered iterations in relation to the size of the loading increment
has been studied by Ambati et al. \cite{AGL15}. Their results show that an insufficient number
of $d$-$u$ iterations can lead to inaccurate results when large loading increments are used.
However, the number of staggered iterations usually has relatively little performance impact
on the shape of the load-displacement curves when loading increments are sufficiently small.
This implies that the staggered approach without iteration can be used as long as
small loading increments are taken.
Since our focus in this work is on mesh adaptation and convergence of Newton's iteration,
we use the staggered approach without iteration for quasi-static brittle fracture problems
with small loading increments.

It is worth mentioning that for the staggered approach with/without iteration, the equation for $u$
remains highly nonlinear due to the decomposition of the strain tensor, which can make
Newton's iteration fail or be slow to converge; see \cite{AGL15} or the numerical results
in Section~\ref{SEC:numerics}. On the other hand, like implicit schemes versus explicit
schemes for ordinary differential equations, the staggered approach with iteration or the monolithic
approach can be more robust than the staggered approach without iteration which is used
in the current work. Comparison of their performances with mesh adaptation and regularization
of the decomposition of the strain tensor (see discussion below) can be an interesting
topic for future investigations.

The model parameter that describes the width of smeared cracks is needed to be small
for the phase-field model to be a reasonably accurate approximation of the original problem.
This in turn requires that the mesh elements be small at least in the crack regions, meaning 
that mesh adaptation is necessary to improve the computational accuracy and efficiency.
The mesh adaptation should be dynamical too since cracks can propagate under continuous load.

Moving mesh methods are well suited for the numerical simulation of the phase-field modeling
of brittle fracture. Although they have been successfully applied to phase-field models for other applications,
e.g., see \cite{DZ08,MR02,SY09,WLT08,YFLS06, YCD08}, they have not been employed for brittle fracture simulation,
which has distinguished challenges associated with the above mentioned decomposition of the strain tensor.
As a matter of fact, mesh adaptation has rarely been employed in brittle fracture simulation so far and
there are only a few published studies on the topic.
Noticeably, Heister et al. \cite{HWW15} develop a predictor-corrector local mesh
adaptivity scheme that allows the mesh to refine around cracks. Artina et al. \cite{AFMP15} present
an a posteriori error estimator for anisotropic mesh adaptation that generates thin, anisotropic elements around
cracks and isotropic elements away from the cracks, but they use an early model that does not decompose
the strain tensor and therefore does not distinguish between fracture caused by tension and compression.

The objective of this paper is twofold. The first is to study the MMPDE
(moving mesh partial differential equation) moving mesh method \cite{BHR09,HRR94a, HRR94b, HR11}
for the phase-field modeling of brittle fracture. The MMPDE method is a type of dynamic mesh adaptation method
specially designed for time dependent problems. It employs a mesh PDE to move the mesh continuously
in time to follow and adapt to evolving structures in the solution. A new formulation of the MMPDE method
was developed recently in \cite{HK15a}, which provides a simple, compact analytical formula for the nodal mesh
velocities (cf. (\ref{MMPDE2}) below), and this makes its implementation relatively easy.
It is also shown in \cite{HK15b} that the mesh governed by the underlying mesh equation
stays nonsingular if it is nonsingular initially.
The MMPDE moving mesh method is combined with a piecewise linear finite element discretization
to solve the governing equation in a quasi-static condition, with the quasi-time being introduced to represent
the load increments.

It is noted that a number of other moving mesh methods have also been developed
in the past and there is a large literature in the area. The interested reader is referred to the books
or review articles \cite{Bai94a,Baines-2011,BHR09,HR11,Tan05} and references therein.

The second goal of the paper is to investigate the convergence of Newton's iteration used for solving
the nonlinear algebraic system resulting from the finite element discretization of the displacement equation.

As mentioned above, Newton's iteration often fails to converge due to the non-existence of the Jacobian matrix
and nonlinearity of the governing equation caused by the decomposition of the strain tensor.
To avoid this difficulty, we propose to smooth out the decomposition with regularization
and discuss three regularization methods. Numerical examples show that all of them can effectively
improve the convergence of Newton's iteration for relatively small values of the regularization parameter
but without comprising the accuracy of the numerical solution. 
Note that no special treatment is needed for Newton's iteration to guarantee its convergence.
This is in contrast to previous works such as \cite{GL16,Wic17a,Wic17b} where modified Newton's schemes
are used along with some modification in the computation of the Jacobian matrix.
Our numerical results also show that
the proposed moving mesh finite element method improves the computational accuracy and efficiency
and is able to handle multiple and complex crack systems.

The rest of this paper is organized as follows. Section~\ref{SEC:phase-field}
provides a brief introduction of the phase-field model for brittle fracture.
The moving mesh finite element method and
the three regularization methods for the decomposition of the strain tensor are described in Section~\ref{SEC:mmfem}.
Numerical results obtained for three two-dimensional examples are presented
in Section~\ref{SEC:numerics}, and conclusions are drawn in Section~\ref{SEC:conclusion}.

\section{Phase-field models for brittle fracture}
\label{SEC:phase-field}

\subsection{Variational approach to elasticity models without cracks}

We first describe small strain isotropic elasticity models without cracks. We consider an elastic body
occupying a bounded domain $\Omega \subset \mathbb{R}^2$ with the boundary
$\partial \Omega = \partial \Omega_t \cup \partial \Omega_u$,
where the surface traction $\overline t$ is specified on $\partial \Omega_t$ and
the displacement $\overline u$ is given on $\partial \Omega_u$. The strain tensor is given by
\[
\epsilon = \frac{1}{2} \left (\nabla u + (\nabla u)^T\right ),
\]
where $\nabla {u}$ is the displacement gradient tensor. For isotropic material without damage, the elastic energy per unit volume, or strain energy density, is given by Hooke's law as
\begin{equation}
W_e(\epsilon) = \frac{\lambda}{2} \left (\tr({\epsilon})\right )^2 + \mu \tr({\epsilon}^2) ,
\label{energy-e-1}
\end{equation}
where $\lambda$ and $\mu$ are the Lam\'{e} constants. The stress tensor is given
by the derivative of the strain energy with respect to the strain tensor, i.e.,
\begin{equation}
\sigma : =  \frac{\partial W_e}{\partial {\epsilon}} = \lambda \tr({\epsilon}) I + 2 \mu {\epsilon} ,
\label{sigma-1}
\end{equation}
where the symbol $:=$ stands for ``definition''.
Then the total strain energy stored in the elastic solid is given by
\begin{equation}
\mathcal{W}_{e}({\epsilon}) = \int_\Omega W_e({\epsilon}) \, d \Omega .
\label{energy-0}
\end{equation}
The variation of $\mathcal{W}_e$ is 
\begin{equation}
\delta \mathcal{W}_e = \int_\Omega \frac{\partial W_e}{\partial \epsilon}:\delta \epsilon \, d \Omega
= \int_\Omega \sigma :  \epsilon(\delta u) \, d \Omega ,
\label{variation}
\end{equation}
where $A:B$ is the inner product of tensors $A$ and $B$, i.e., $A:B = \sum\limits_{i,j} A_{i,j} B_{i,j}$.
Define the function spaces
\begin{gather*}
V_{u} = \left \{ \; \varphi \; | \; \varphi \in H^1(\Omega),\;  \varphi = \overline u \; \text{on} \; \partial \Omega_u \right \},\\
V_{u} ^0= \left \{ \; \varphi \; | \; \varphi \in H^1(\Omega),\; \varphi = 0 \; \text{on} \; \partial \Omega_u \right \},
\end{gather*}
where $H^1(\Omega)$ is a Sobolev space, viz.,
\[
H^1(\Omega) = \left \{ \; \varphi \; | \;  \int_{\Omega} \varphi^2 \, d \Omega < +\infty,\;  
\int_{\Omega} (\nabla \varphi)^2 \, d \Omega < + \infty \right \}.
\]
If we add the boundary traction $\overline{t}$ and body force $f$, then the variational formulation of the elasticity model
is to find $u \in V_u$ such that 
\begin{equation}
\int_\Omega \sigma :  {\epsilon}(\delta {u}) \, d \Omega = 
\int_{\partial \Omega_t} \overline{t} \cdot \delta {u} \, dS +
\int_\Omega f \cdot \delta {u} \, d \Omega ,
\quad \forall \; \delta u \in V_u^0.
\label{pde-1}
\end{equation}

\subsection{Phase-field approach to elasticity models with cracks}

We now consider the situation with cracks in the elastic body.
We use the approach proposed by Francfort and Marigo \cite{FM98} where the total energy of the body
with a given crack $\Gamma$ is given by
\[
\mathcal{W}(\eps,\Gamma) = \mathcal{W}_{e}(\eps,\Gamma) + \mathcal{W}_{c}(\Gamma)
 : = \int_{\Omega \setminus \Gamma} W_{e}(\eps,\Gamma) \, d \Omega + \int_{\Gamma} g_c \, dS,
\]
where $\mathcal{W}_{e}(\eps,\Gamma)$ represents the energy stored in the bulk of the elastic body, 
$\mathcal{W}_{c}(\Gamma)$ is the energy required to create the crack according to the Griffith criterion,
and $g_c$ is the fracture energy density (also referred to as the fracture toughness) which is the amount of
energy needed to create a unit area of fracture surface.

\begin{figure} 
\centering 
\subfigure[sharp crack $\Gamma$ in the solid]{\label{fig:subfig:sharp crack}
\includegraphics[width=0.3\linewidth]{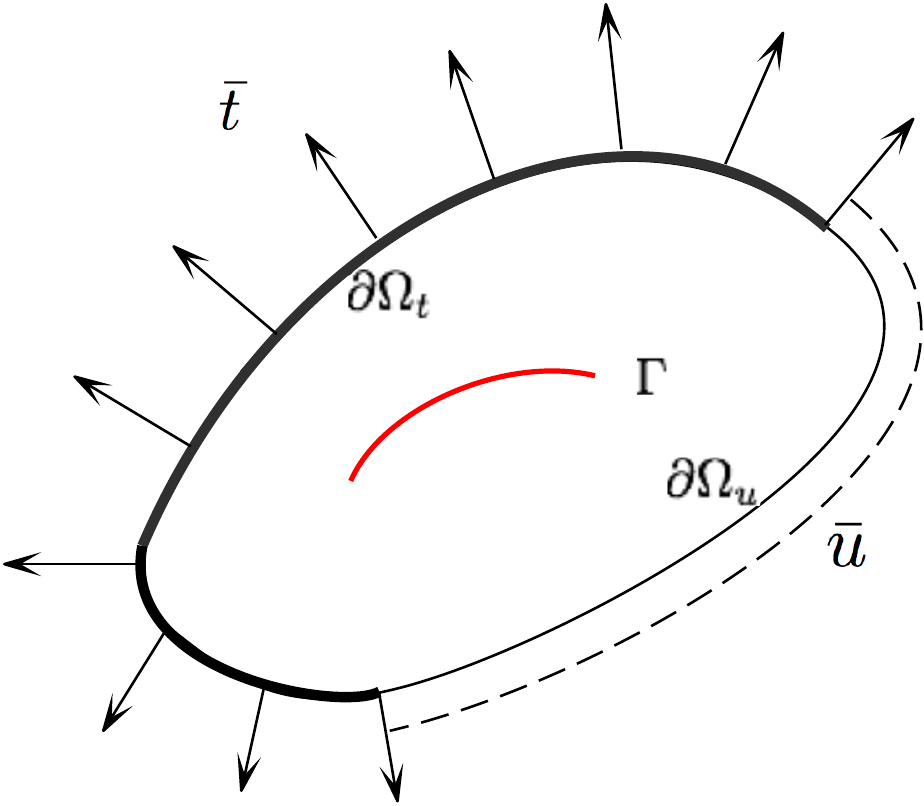}}
\hspace{10mm}
\subfigure[regularized crack $\Gamma_l(d)$ by the phase-field approximation]{\label{fig:subfig:regularized crack}
\includegraphics[width=0.3\linewidth]{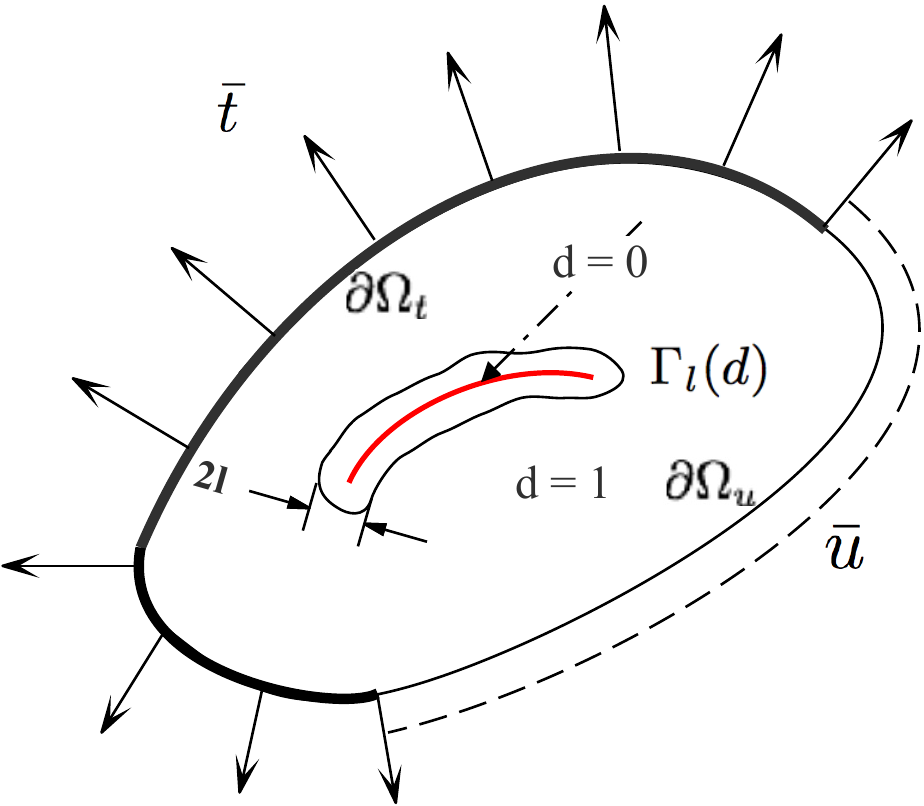}}
\caption{A sketch of the problem setting for brittle crack.}
\label{sketch-of-crack}
\end{figure}

In the phase-field modeling, the fracture surface is approximated by a phase-field variable $d(x,t)$,
which depends on a parameter $l$ describing the width of the smooth approximation of the crack.
This function is smooth with the value 0 or close to 0 near the crack and 1 away from the crack
(see Fig. \ref{sketch-of-crack}). The fracture energy $\mathcal{W}_{c}(\Gamma)$ is approximated
by the smeared total fracture energy \cite{BFM00} as
\begin{equation}
\label{energy-c}
\mathcal{W}_{c}^{l}  = \int_\Omega \frac{g_c}{4 l} \left ((d-1)^2 + 4 l^2 |\nabla d|^2 \right ) d \Omega.
\end{equation}

The elastic energy $\mathcal{W}_{e}(\eps,\Gamma)$ needs to be modified to reflect the loss
of material stiffness in the damage zone. We follow the approach by Miehe et al. \cite{MHW10}. 
Define the decomposition of a scalar function $f$ as
\[
f = f^+ + f^-,\quad f^+ = \frac{f+|f|}{2}, \quad f^- = \frac{f-|f|}{2}.
\]
Assuming that the strain tensor has the eigen-decomposition
$\epsilon = Q \text{diag}(\lambda_1, ..., \lambda_n) Q^T$, we define
\begin{equation}
\label{strain-decomp}
\epsilon = \epsilon^+ + \epsilon^-,
\quad 
\epsilon^+ = Q \text{diag}(\lambda_1^+, ..., \lambda_n^+) Q^T,
\quad 
\epsilon^- = Q \text{diag}(\lambda_1^-, ..., \lambda_n^-) Q^T .
\end{equation}
The \textit{tensile} strain component, $\epsilon^+$, contributes to the damage process
resulting in crack initiation and propagation whereas the \textit{compression} strain component, $\epsilon^-$,
does not contribute to the damage process. A commonly used damage model is given as 
\begin{equation}
W_e =  g(d) \Psi^+(\epsilon) +  \Psi^-(\epsilon),
\label{energy-e}
\end{equation}
where $g(d)$ is a degradation function that describes the reduction of the stiffness of the bulk of the solid and
\[
\Psi^+(\epsilon) = \frac{\lambda}{2} \left ( (\tr(\epsilon))^+\right )^2 + \mu \tr((\epsilon^+)^2) ,
\quad \Psi^-(\epsilon) = \frac{\lambda}{2} \left ( (\tr(\epsilon))^-\right )^2 + \mu \tr((\epsilon^-)^2) .
\]
The degradation function $g(d)$ is required to satisfy the property (e.g., see \cite{KSM2015})
\[
\begin{cases}
g(0) = 0: &\quad \text{Damage occurred for $d = 0$ and this part should vanish}; \\
g(1) = 1: &\quad \text{No damage occurs for $d = 1$}; \\
g'(0) = 0: &\quad \text{No more changes after the fully broken state}; \\
g'(1) \ne 0: &\quad \text{The damage has to be initiated at the onset.}
\end{cases}
\]
We take a commonly used quadratic degradation function $g(d) = d^2$ in our computation.
Combining \eqref{energy-c} and \eqref{energy-e}, we obtain the total energy as
\begin{equation}
\mathcal{W}^l =  \mathcal{W}_e + \mathcal{W}_c^l =
 \int_{\Omega} \left(
(d^2+k_l)\Psi^+(\epsilon) + \Psi^-(\epsilon) +
 \frac{g_c}{4l}\left ((d-1)^2+4l^2 |\nabla d|^2 \right ) 
 \right) d \Omega,
\label{energy-d-2}
\end{equation}
where $k_l \ll l$ is the (small) regularization parameter for avoiding degeneracy.
The variation of the energy is 
\[
\delta \mathcal{W}^l =
\int_{\Omega} \frac{\partial W}{\partial d} \delta d \; d \Omega +
\int_{\Omega} \frac{\partial W}{\partial \nabla d} \cdot \nabla \delta d \; d \Omega +
\int_{\Omega} \frac{\partial W}{\partial \epsilon} : \epsilon(\delta u) \; d \Omega
\]
where
\[
W = (d^2+k_l)\Psi^+(\epsilon) + \Psi^-(\epsilon) +
 \frac{g_c}{4l}\left ((d-1)^2+4l^2 |\nabla d|^2 \right ) .
 \]
The weak formulation is to find $d \in V_d = H^1(\Omega)$ (with $d$ satisfying a homogeneous
Neumann boundary condition) and $u \in V_{u}$ such that 
\begin{align}
& \int_\Omega \left (
 \left (2 d \mathcal{H} + \frac{g_c (d-1) }{2 l} \right ) \delta d
 +  2 g_c l \nabla d \cdot \nabla \delta d \right) d \Omega = 0 ,
 \quad \forall\; \delta d \in V_d 
\label{weak-2-d}
\\
& \int_\Omega \sigma : \epsilon(\delta u) d \Omega = 
\int_{\Omega_t} \overline{t} \cdot \delta {u} \; dS +
\int_\Omega f \cdot \delta {u} \; d \Omega ,
\quad \forall\; \delta u \in V_u^0
\label{weak-2-u}
\end{align}
where $\mathcal{H} = \Psi^+(\epsilon)$ and
\begin{equation}
\sigma := \frac{\partial W}{\partial \epsilon} = (d^2+k_l) \left (\lambda (\tr(\epsilon))^+ I + 2 \mu \epsilon^+ \right )
+ \left (\lambda (\tr(\epsilon))^- I + 2 \mu \epsilon^- \right ) .
\label{sigma-2}
\end{equation}
To ensure crack irreversibility in the sense that the cracks can only grow,
we replace $\mathcal{H} = \Psi^+(\epsilon)$ in (\ref{weak-2-d}) by
\begin{equation}
\label{eqn-H}
\mathcal{H} = \max_{s \le t} \Psi^+(\epsilon)(s),
\end{equation}
where $t$ is the quasi-time corresponding to the load increments.

In addition to the above described decomposition by Miehe et al. \cite{MHW10},
there is another often-used decomposition by Amor et al. \cite{AMM09} that is based on
the split of the strain tensor into volumetric and deviatoric parts.
Discussion of the advantages and disadvantages of both decompositions can be
found in \cite{AGL15,Bor12}.

\section{The moving mesh finite element approximation}
\label{SEC:mmfem}

\subsection{Finite element discretization and solution procedure}
\label{SEC:procedure}

We consider a simplicial mesh $\mathcal{T}_h$ for the domain $\Omega$ and denote the number
of its elements and vertices by $N$ and $N_v$, respectively. The piecewise linear approximations
of the function spaces $V_d$, $V_u$ and $V_u^0$ are given by
\begin{align*}
& V_d^h =  \left \{ 
\varphi_h\; | \; \varphi_h \in C^0(\overline{\Omega}) \cap V_d; \;
\varphi_h |_K \in P_1(K),\;  \forall K \in \mathcal{T}_h(\Omega) 
\right \},\\
& V_u^h =  \left \{ 
\varphi_h \; | \; \varphi_h \in C^0(\overline{\Omega}) \cap V_u; \;
\varphi_h |_K \in P_1(K),\;  \forall K \in \mathcal{T}_h(\Omega) 
\right \},\\
& V_u^{0,h} =  \left \{ 
\varphi_h \; | \; \varphi_h \in C^0(\overline{\Omega}) \cap V_u^0; \;
\varphi_h |_K \in P_1(K),\;  \forall K \in \mathcal{T}_h(\Omega) 
\right \},
\end{align*}
where $P_1(K)$ is the set of polynomials of degree less than or equal to 1 defined on $K$.
The linear finite element approximation of the phase-field problem (\ref{weak-2-d}) and (\ref{weak-2-u}) is
to find $d_h \in V_d^h$ and $u_h \in V_u^h$ such that
\begin{align}
& \int_\Omega \left (
 \left (2 d_h \mathcal{H} + \frac{g_c (d_h-1) }{2 l} \right ) \varphi_h
 +  2 g_c l \nabla d_h \cdot \nabla \varphi_h \right) d \Omega = 0 ,
 \quad \forall\; \varphi_h \in V_d^h 
\label{discrete-d}
\\
& \int_\Omega \sigma(u_h) : \epsilon(\varphi_h) d \Omega = 
\int_{\Omega_t} \overline{t} \cdot \varphi_h \; dS +
\int_\Omega f \cdot \varphi_h \; d \Omega ,
\quad \forall\; \varphi_h \in V_u^{0,h}.
\label{discrete-u}
\end{align}

Note that $u_h$ and $d_h$ are strongly coupled through (\ref{discrete-d}) and (\ref{discrete-u}).
Common procedures for solving these equations can be roughly categorized into two groups,
simultaneous (or called monolithic) solution and alternating (or called staggered) solution.
In the former group (such as see \cite{GL16,HWW15,Wic17a,Wic17b})
the phase-field variable and displacement
are solved simultaneously often by Newton's method. This procedure has to deal with a large system and
the highly nonlinear coupling between $d_h$ and $u_h$. Moreover, the decomposition of the strain tensor
further increases the nonlinearity (see Section \ref{SEC:convergence} for detailed
discussion). It is still challenging to obtain a convergent solution using the simultaneous solution procedure.
The second solution procedure (such as see \cite{AGL15, AMM09, MHW10, MWH10})
is to solve (\ref{discrete-d}) for $d_h$ and (\ref{discrete-u}) for $u_h$ alternatingly.
One of the advantages of this approach is that $d_h$ and $u_h$ are not coupled.
Moreover, (\ref{discrete-d}) is linear about $d_h$ and its solution does not need Newton's iteration.

In our computation, we use the alternating solution procedure and consider the problem in a quasi-static condition.
The quasi-time $t$ is introduced to represent the load increments. The solution procedure from  $t^n$ to $t^{n+1}$
is described as follows.

\begin{enumerate}
\item [(i)] We assume that the mesh $\mathcal{T}_h^n$ at time $t^n$ and
the history field $\mathcal{H}_h^n$ in \eqref{eqn-H} (defined on $\mathcal{T}_h^n$) are given.
\item [(ii)] Compute the phase-field variable $d_h^{n+1}$ and new mesh $\mathcal{T}_h^{n+1}$ as follows.
	\begin{itemize}
	\item Set $\mathcal{T}_h^{n+1,1} = \mathcal{T}_h^n$;
	\item For $k = 1:kk$
		\begin{itemize}
		\item[-] Compute $\mathcal{H}$ on $\mathcal{T}_h^{n+1,k}$ using linear interpolation of $\mathcal{H}_h^n$;
		\item[-] Compute $d_h^{n+1,k}$ using (\ref{discrete-d}) and $\mathcal{H}$ on $\mathcal{T}_h^{n+1,k}$;
		\item[-] If $k < kk$, compute the new mesh $\mathcal{T}_h^{n+1,k+1}$ by the MMPDE moving mesh method
		based on $\mathcal{T}_h^{n+1,k}$ and $d_h^{n+1,k}$.  (See Section~\ref{SEC:MMPDE})
		\end{itemize}
	\item Set $\mathcal{T}_h^{n+1} = \mathcal{T}_h^{n+1,kk}$ and $d_h^{n+1} = d_h^{n+1,kk}$.
	\end{itemize}
\item [(iii)] Compute the displacement field $u_h^{n+1}$ by solving the nonlinear system (\ref{discrete-u}) based on $d_h^{n+1}$ and $\mathcal{T}_h^{n+1}$. (See Section~\ref{SEC:convergence})

\item [(iv)] Compute $\Psi^{+,n+1}_h (\epsilon(u_h^{n+1}))$ and
interpolate $\mathcal{H}_h^n$ from the old mesh $\mathcal{T}_h^{n}$
to the new mesh $\mathcal{T}_h^{n+1}$ using linear interpolation and denote it by $\tilde{\mathcal{H}}_h^{n+1}$.
Let $\mathcal{H}_h^{n+1} = \max\{\Psi^{+,n+1}_h (\epsilon(u_h^{n+1})), \tilde{\mathcal{H}}_h^{n+1}\}$.
\end{enumerate}

In the above procedure, the parameter $kk$ can affect the adaptivity of the mesh $\mathcal{T}_h^{n+1}$
to the phase-field  variable $d_h^{n+1}$. From our computational experience, we choose $kk=5$, for which
we have found that $\mathcal{T}_h^{n+1}$ is sufficiently adaptive to $d_h^{n+1}$ but without compromising
the computational efficiency too much.

\subsection{Convergence of Newton's iteration} 
\label{SEC:convergence}

\begin{figure} [!htb]
\centering 
\subfigure[$\lambda^+(\lambda)$]{\label{fig:subfig:lambda+}
\includegraphics[width=0.32\linewidth]{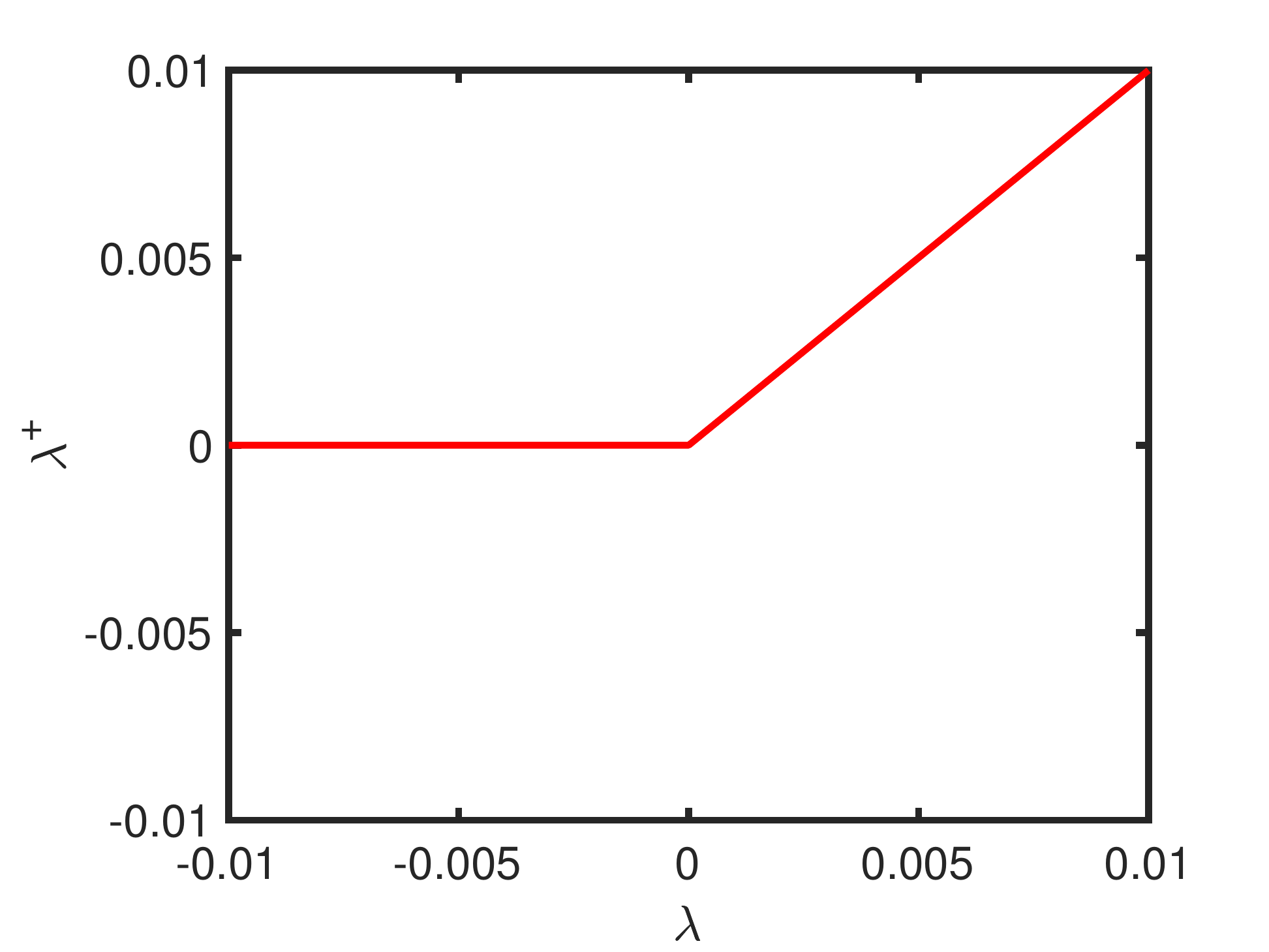}}
\subfigure[$\lambda^-(\lambda)$]{\label{fig:subfig:lambda-}
\includegraphics[width=0.32\linewidth]{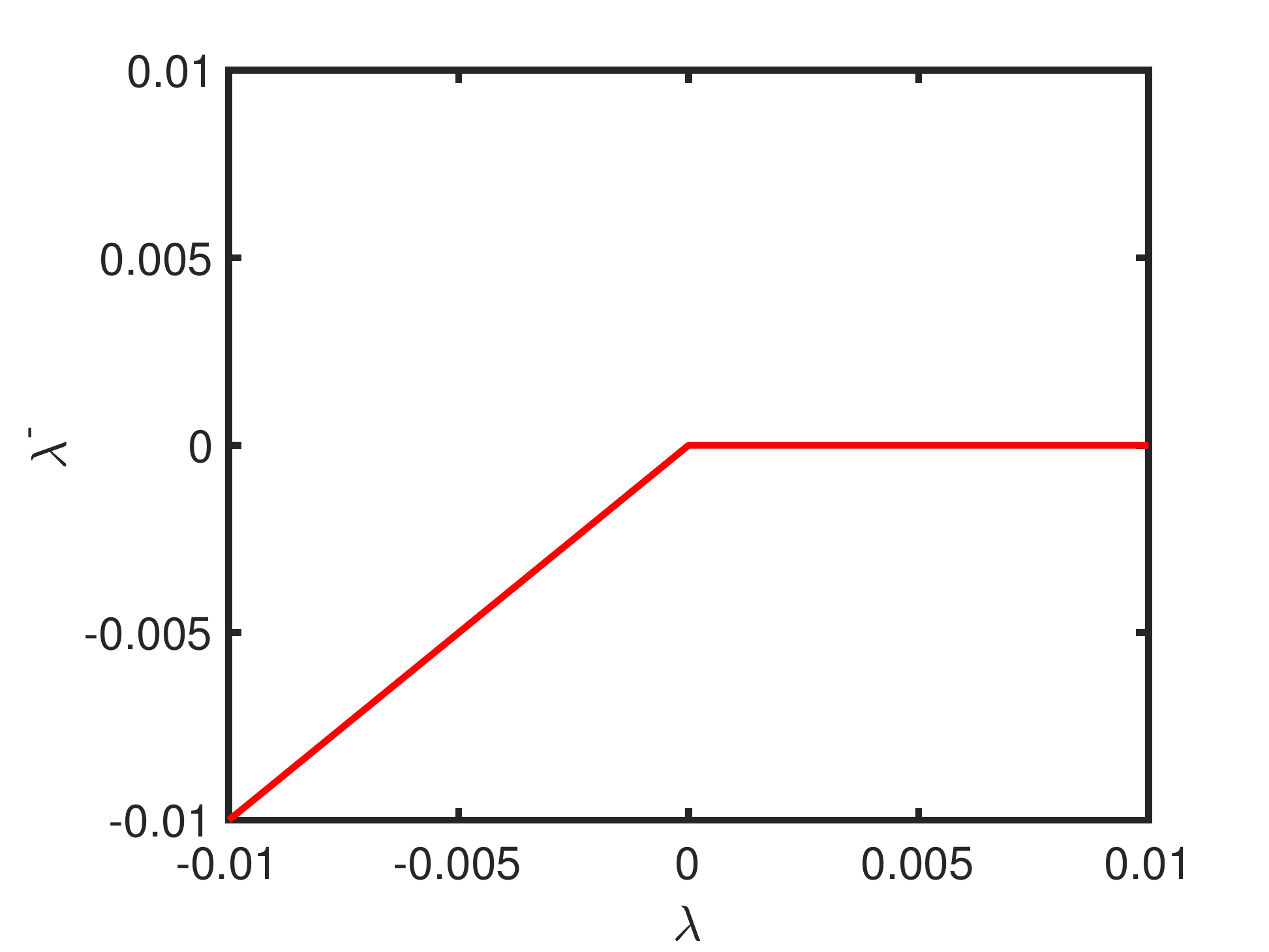}}
\subfigure[Regularization]{\label{fig:subfig:global}
\includegraphics[width=0.32\linewidth]{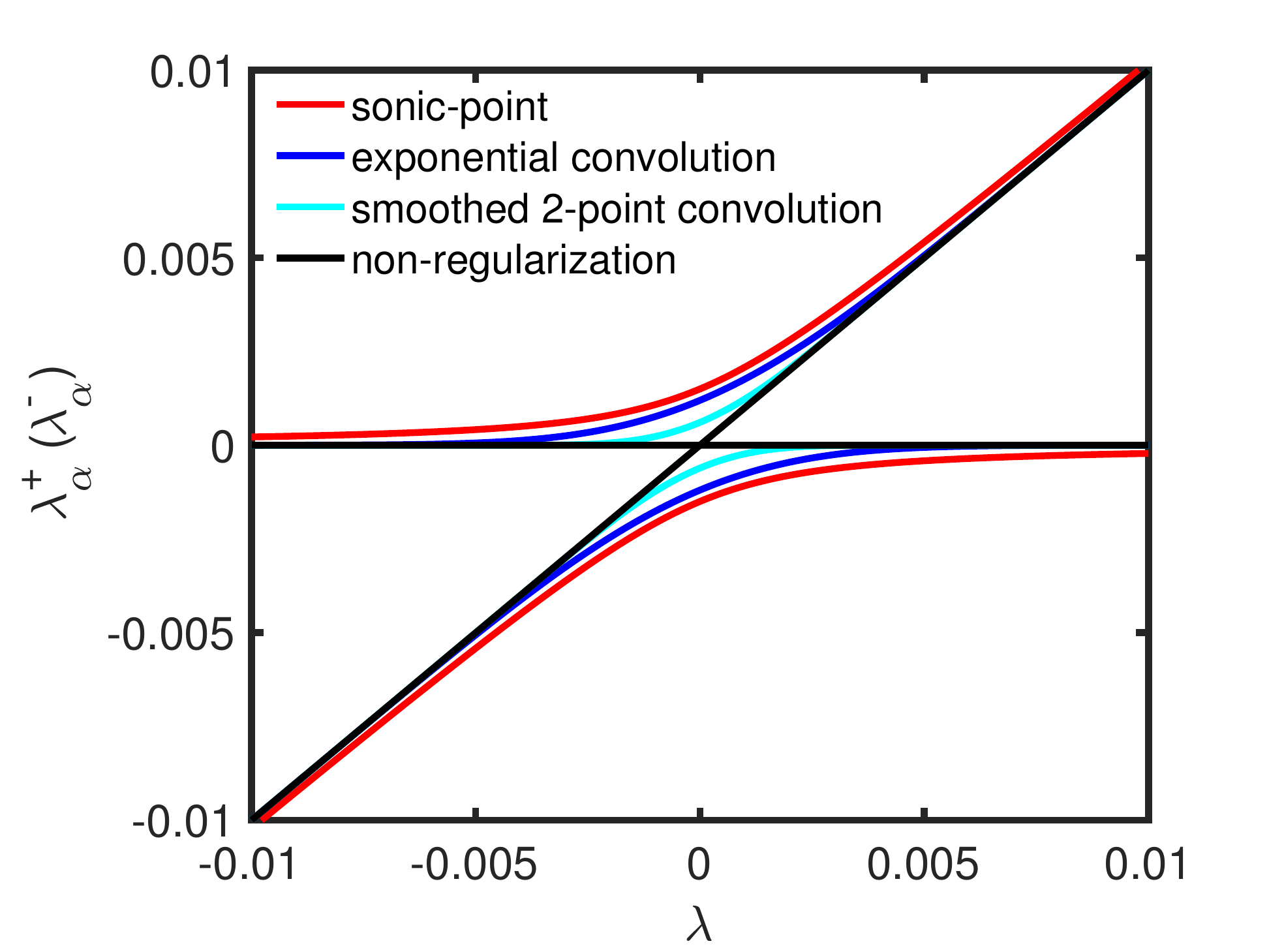}}
\caption{Decomposition and regularization of eigenvalue function with the regularization parameter $\alpha = 3 \times 10^{-3}$.}
\label{fig:lambda}
\end{figure}

We recall that the strain tensor is decomposed into tensile and compressive components
(cf. \eqref{strain-decomp}), with only the former component contributing to generation of cracks. 
The tensile and compressive components are nonsmooth functions of the displacement, which
has two major effects on Newton's iteration applied to the solution of the equation (\ref{discrete-u}),
the existence and computation of the Jacobian matrix and the convergence of Newton's iteration.
 As can be seen in Fig. \ref{fig:subfig:lambda+} and Fig. \ref{fig:subfig:lambda-}, 
the Jacobian matrix of equation (\ref{discrete-u}) does not exist in the places where any
of the eigenvalues vanishes. In other places where the eigenvalues do not vanish, 
the computation of the Jacobian matrix can also be tricky.
In principle, the analytical expressions of the derivatives of $\epsilon^{+}$ and $\epsilon^{-}$ with respect
to $\epsilon$ can be obtained through the derivatives of the eigenvalues and eigenvectors with respect to
$\epsilon$ (e.g, see \cite{Magnus-1985}).
However, these expressions are too complicated to be used in practical computation.
Special algorithms are needed to compute these derivatives using analytical formulas; see Miehe \cite{Mie98}
and Miehe and Lambrecht \cite{ML01}.

Another effect of the decomposition of the strain tensor is that Newton's iteration often fails to converge;
see the numerical examples in Section~\ref{SEC:effect}. The degeneracy of the equation (\ref{discrete-u})
can further complicate the situation. In principle, the regularization parameter $k_l$ can be chosen sufficiently large to
make Newton's iteration to converge. However, a large value of $k_l$ will smear the crack and lead to
an unphysically overestimated bulk energy. 

To overcome the above mentioned difficulties, we propose to smooth out the decomposition \eqref{strain-decomp}
and compute the Jacobian matrix of the equation (\ref{discrete-u}) using finite differences.
The latter is straightforward so we will not elaborate it here.
For the former, there are a variety of regularization methods. We consider here three of these methods
where the positive and negative eigenvalue functions are smoothed using a switching technique or
convolution with a smoothed delta function (or a mollifier).
\begin{itemize}
\item {\it The sonic-point regularization method.}
We first consider the so-called eigenvalue-switching technique which is used
to obtain a smooth transition of the solution through the sonic point \cite{ATL86} in computational fluid dynamics.
The regularized positive and negative eigenvalue functions are defined as
\[
\lambda_{\alpha}^+ = \frac{\lambda + \sqrt{\lambda^2 + \alpha^2}}{2} , \quad
\lambda_{\alpha}^- = \frac{\lambda - \sqrt{\lambda^2 + \alpha^2}}{2},
\]
where $\alpha > 0$ is the regularization parameter. 

\item {\it The exponential convolution method.}
The second regularization method is to take the convolution of $\lambda^+$ and $\lambda^-$
with the exponential delta function
\[
\delta_{\alpha}(\lambda) = \frac{1}{\sqrt{2\pi} \alpha} e^{-\frac{\lambda^2}{2 \alpha^2}}.
\]
Then the regularized positive and negative eigenvalue functions are given by
\begin{align*}
& \lambda_{\alpha}^{+}
= \int_{-\infty}^{\infty} \lambda^+(\eta) \delta_{\alpha}(\lambda-\eta) d \eta
= \frac{\lambda}{2} \left( 1+ \text{erf}\left (\frac{\lambda}{\sqrt{2}\alpha}\right ) \right)
+ \frac{\alpha}{\sqrt{2\pi}} e^{-\frac{\lambda^2}{2\alpha^2}},
\\
& \lambda_{\alpha}^{-}
= \int_{-\infty}^{\infty} \lambda^-(\eta) \delta_{\alpha}(\lambda-\eta) d \eta 
= \frac{\lambda}{2} \left( 1+ \text{erf}\left (-\frac{\lambda}{\sqrt{2}\alpha}\right ) \right) 
- \frac{\alpha}{\sqrt{2\pi}} e^{-\frac{\lambda^2}{2\alpha^2}}.
\end{align*}
\item {\it The smoothed 2-point convolution method.}
The third method is similar to the exponential convolution method but with a localized, smoothed 2-point delta function 
\[
\delta_{\alpha}(\lambda) = 
\begin{cases}
\frac{1}{\alpha} \left( \frac{3}{4} - \frac{\lambda^2}{\alpha^2} \right), \quad & \text{ for } |\lambda| \leqslant 0.5 \alpha
\\
\frac{1}{\alpha} \left( \frac{9}{8} - \frac{3}{2}\frac{|\lambda|}{\alpha} + \frac{1}{2}\frac{\lambda^2}{\alpha^2} \right), \quad
& \text{ for } 0.5\alpha \leqslant |\lambda| \leqslant 1.5 \alpha
\\
0, \quad & \text{ for } |\lambda| \geqslant 1.5 \alpha .
\end{cases}
\]
We then have
\begin{align*}
& \lambda_{\alpha}^+ = 
\begin{cases}
0,\quad & \text{ for } \lambda \leqslant -1.5 \alpha
\\
\frac{1}{\alpha} \left( \frac{1}{24\alpha^2} \lambda^4 + \frac{1}{4\alpha}\lambda^3 + \frac{9}{16}\lambda^2+ \frac{9\alpha}{16} \lambda + \frac{27\alpha^2}{128} \right),\quad & \text{ for } -1.5\alpha \leqslant \lambda \leqslant -0.5\alpha 
\\
\frac{1}{\alpha} \left( -\frac{1}{12\alpha^2}\lambda^4 + \frac{3}{8}\lambda^2 + \frac{\alpha}{2}\lambda + \frac{13\alpha^2}{64} \right), \quad & \text{ for } -0.5\alpha \leqslant \lambda \leqslant 0.5\alpha
\\
\frac{1}{\alpha} \left( \frac{1}{24\alpha^2} \lambda^4 - \frac{1}{4\alpha}\lambda^3 + \frac{9}{16}\lambda^2+ \frac{7\alpha}{16} \lambda + \frac{27\alpha^2}{128} \right),\quad & \text{ for } 0.5\alpha \leqslant \lambda \leqslant 1.5\alpha
\\
\lambda, \quad & \text{ for } \lambda \geqslant 1.5\alpha
\end{cases}
\\
& \lambda_{\alpha}^- = 
\begin{cases}
\lambda,\quad & \text{ for } \lambda \leqslant -1.5 \alpha
\\
\frac{1}{\alpha} \left( -\frac{1}{24\alpha^2} \lambda^4 - \frac{1}{4\alpha}\lambda^3 - \frac{9}{16}\lambda^2+ \frac{7\alpha}{16} \lambda - \frac{27\alpha^2}{128} \right),\quad & \text{ for } -1.5\alpha \leqslant \lambda \leqslant -0.5\alpha 
\\
\frac{1}{\alpha} \left( \frac{1}{12\alpha^2}\lambda^4 - \frac{3}{8}\lambda^2 + \frac{\alpha}{2}\lambda - \frac{13\alpha^2}{64} \right), \quad & \text{ for } -0.5\alpha \leqslant \lambda \leqslant 0.5\alpha
\\
\frac{1}{\alpha} \left( -\frac{1}{24\alpha^2} \lambda^4 + \frac{1}{4\alpha}\lambda^3 - \frac{9}{16}\lambda^2+ \frac{9\alpha}{16} \lambda - \frac{27\alpha^2}{128} \right),\quad & \text{ for } 0.5\alpha \leqslant \lambda \leqslant 1.5\alpha
\\
0, \quad & \text{ for } \lambda \geqslant 1.5\alpha .
\end{cases}
\end{align*}
\end{itemize}

It is remarked that all of the above regularization methods satisfy
\begin{equation}
\lambda = \lambda^+_{\alpha} + \lambda^-_{\alpha} 
\label{consist-1}
\end{equation}
and thus
\begin{equation}
\epsilon = \epsilon_\alpha^+ + \epsilon_\alpha^- .
\label{consist-2}
\end{equation}
Moreover, as can be seen in Fig. \ref{fig:subfig:global}, all of the methods provide a smooth transition
around zero for both positive and negative eigenvalue functions. Furthermore,
the sonic-point and exponential convolution methods
have global effects whereas the smoothed 2-point convolution method only changes the values near
$\lambda = 0$. Finally, for a fixed value of the regularization parameter, $\alpha = 3 \times 10^{-3}$,
in terms of the closeness of the curve of $\lambda_{\alpha}^+$ ($\lambda_{\alpha}^-$) to that of $\lambda^+$
(or $\lambda^-$) the best is the smoothed 2-point convolution and then the exponential convolution
and sonic-point methods.

The effects of these regularization methods on the the convergence of Newton's iteration and
on the numerical solution will be discussed in Section \ref{SEC:effect}.
It is worth mentioning that we have tried a few other regularization methods which only modify
the local behavior of the eigenvalue function near the origin. Since they do not satisfy (\ref{consist-1})
in general and are much less effective in making Newton's iteration convergent than the above methods,
we choose not to report them here to save space.

\subsection{The MMPDE moving mesh method}
\label{SEC:MMPDE}

In this section we describe the MMPDE moving mesh method \cite{HRR94a,HR11}
for generating the new mesh $\mathcal{T}_h^{n+1}$
(as mentioned in Section \ref{SEC:procedure}).
The method takes the $\mathbb{M}$-uniform mesh approach where a nonuniform mesh is viewed
as a uniform one in the metric specified by a tensor $\mathbb{M}$. 
In our computation, we choose the metric tensor as
\begin{equation}
\mathbb{M} = \det(I+|H(d_h)|)^{-\frac{1}{6}} (I+|H(d_h)|),
\label{M-1}
\end{equation}
where $H(d_h)$ is a recovered Hessian of $d_h$ and 
$|H(d_h)| = Q \text{diag}(|\lambda_1|,...,|\lambda_2|) Q^T$,
assuming that the eigen-decomposition of $H(d_h)$ is $Q \text{diag}(\lambda_1,...,\lambda_2) Q^T$.
The recovered Hessian of $d_h$ at a vertex is obtained by twice differentiating a local quadratic
polynomial fitting in the least-squares sense to the nodal values of $d_h$ at the neighboring vertices. 
The form of (\ref{M-1}) is known optimal in terms of the $L^2$ norm of linear interpolation error on
triangular meshes (e.g., see \cite{HS03}).
Notice that $\mathbb{M}$ is symmetric and uniformly positive define on $\Omega$. It will be used to
control the shape, size, and orientation of mesh elements  through the so-called equidistribution
and alignment conditions (to be discussed below). Since (\ref{M-1}) is based on the Hessian
of the phase-field variable $d$, we may expect that the mesh elements are concentrated in the crack regions
where the curvature of $d$ is large.

Let $x_1, ...., x_{N_v}$ be the coordinates of the vertices of $\mathcal{T}_h$.
We choose the very initial physical (simplicial) mesh as the reference computational mesh which is denoted
by $\hat{\mathcal{T}}_{c,h} =  \{\hat{\xi}_1,...,\hat{\xi}_{N_v} \}$. For the purpose of mesh generation,
we need an intermediate simplicial mesh which we refer to as a computational mesh,
$\mathcal{T}_{c,h} = \{\xi_1,...,\xi_{N_v} \}$. We assume that $\mathcal{T}_h$, $\hat{\mathcal{T}}_{c,h}$,
and $\mathcal{T}_{c,h}$ have the same number of elements and vertices and the same connectivity.
As a consequence, there exists a unique element $K_c \in \mathcal{T}_{c,h}$ corresponding
to any element $K \in \mathcal{T}_h$.
The affine mapping between $K_c$ and $K$ is denoted by $F_K$ and its Jacobian matrix by $F_K'$. 
Let the coordinates of the vertices of $K$ and $K_c$ be $x_0^K,\, x_1^K,\, x_2^K$ and
$\xi_0^K,\, \xi_1^K,\, \xi_2^K$, respectively. Then we have
\begin{gather*}
F_K(\xi_i^K) = x_i^K, \quad i = 0, 1, 2 \\
F'_K(\xi_i^K- \xi_0^K) = x_i^K - x_0^K, \quad i = 1, 2.
\end{gather*}
From this we get 
\[
F'_K = E_K \hat{E_K}^{-1}, \quad (F'_K)^{-1} = \hat{E_K} E_K^{-1},
\]
where $E_K$ and $\hat{E}_K$ are the edge matrices of $K$ and $K_c$ defined as
\[
E_K = [x_1^K-x_0^K,x_2^K-x_0^K], \quad \hat{E}_K = [\xi_1^K-\xi_0^K, \xi_2^K-\xi_0^K].
\]

The objective of the MMPDE moving mesh method is to generate an adaptive mesh as a uniform one
in the metric $\mathbb{M}$. Such an $\mathbb{M}$-uniform mesh requires that 
(i) the ratio of the area of $K$ in the metric $\mathbb{M}$ to the area of $K_c$ in the Euclidean metric
stay constant for all elements and
(ii) $K$ measured in the metric $\mathbb{M}$ be similar to $K_c$ measured in the Euclidean metric
for all elements.
These two requirements can be expressed mathematically as the equidistribution and alignment conditions
(e.g., see \cite{Hua06,HR11}),
\begin{align}
& |K| \sqrt{\det(\mathbb{M}_K)} = \frac{|\Omega_h|\; |K_c|}{|\Omega_c|} , \quad \forall{K} \in \mathcal{T}_h \\
& \frac{1}{2} \tr\left( (F'_K)^T \mathbb{M}_K F'_K \right) = \det\left( (F'_K)^T \mathbb{M}_K F'_K \right)^{\frac{1}{2}}, \quad K \in \mathcal{T}_h
\end{align}
where $|K|$ and $|K_c|$ denote the area of $K$ and $K_c$, respectively,
$\mathbb{M}_K$ is the average of $\mathbb{M}$ over $K$, $\det(\cdot)$ and $\tr(\cdot)$ denote the
determinant and trace of a matrix, 
and
\[
|\Omega_h| = \sum_{K \in \mathcal{T}_h} |K| \sqrt{\det(\mathbb{M}_K)},
\quad |\Omega_c| = \sum_{K_c \in \mathcal{T}_{c,h}} |K_c|. 
\]
An energy functional based on these conditions has been proposed by Huang \cite{Hua01} as
\begin{equation}
I_h(\mathcal{T}_h; \mathcal{T}_{c,h}) = \sum_{K \in \mathcal{T}_h} |K|
G(\mathbb{J}_K, \det(\mathbb{J}_K), \mathbb{M}_K),
\end{equation}
where $\mathbb{J}_K = (F'_K)^{-1}$ and
\begin{align*}
G(\mathbb{J}_K, \det(\mathbb{J}_K), \mathbb{M}_K) & = 
\theta \sqrt{\det(\mathbb{M}_K)} \left( \tr(\mathbb{J}_K \mathbb{M}_K \mathbb{J}^T_K) \right)^{p} \\
&  + (1-2 \theta) 2^{p} \sqrt{\det(\mathbb{M}_K)}
\left( \frac{\det(\mathbb{J}_K)}{\sqrt{\det(\mathbb{M}_K)}} \right)^p .
\end{align*}
Here, $0 < \theta \le \frac{1}{2}$ and $p > 1$ are two dimensionless parameters. We use $\theta=\frac{1}{3}$ and
$p = \frac{3}{2}$, which are known experimentally to work well for most problems.

\begin{figure} [!htb]
\centering 
\includegraphics[width=0.5\linewidth]{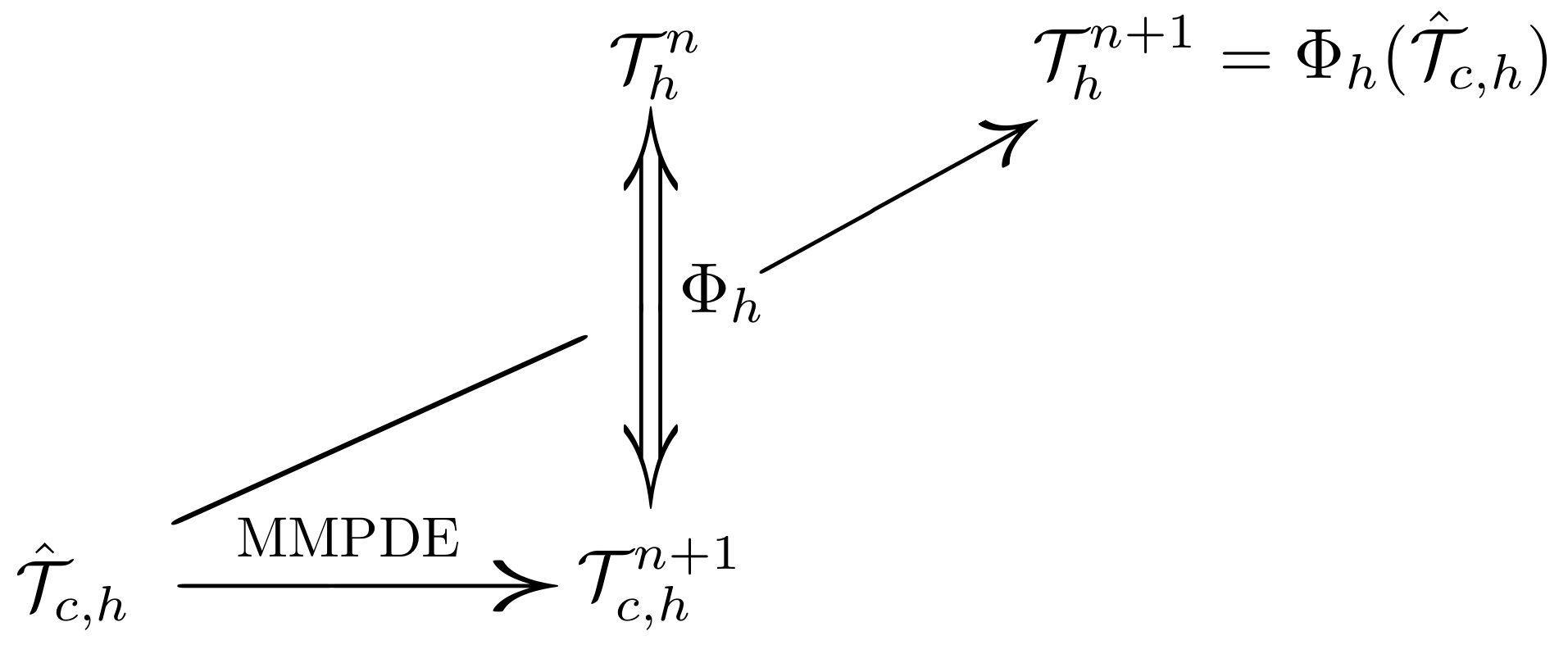}
\caption{A sketch of the procedure for generating the new mesh.}
\label{fig:MMPDE}
\end{figure}

Our goal is to find a new physical mesh $\mathcal{T}_h^{n+1}$ by minimizing $I_h$. We use an indirect approach
with which we take $\mathcal{T}_h=\mathcal{T}_h^n$ and then minimize $I_h$ with respect to $\mathcal{T}_{c,h}$.
The minimization is carried out by integrating the MMPDE (see (\ref{MMPDE1}) below)
from $t^n$ to $t^{n+1}$ with the reference computational mesh $\hat{\mathcal{T}}_{c,h}$ as the initial mesh.
The obtained computational mesh is denoted as $\mathcal{T}_{c,h}^{n+1}$.
Notice that $\mathcal{T}_h^n$ is kept unchanged during the integration and 
$\mathcal{T}_{c,h}^{n+1}$ and $\mathcal{T}_h^n$ form a correspondence, i.e.,
$\mathcal{T}_h^n = \Phi_h(\mathcal{T}_{c,h}^{n+1})$.
The new physical mesh $\mathcal{T}_h^{n+1}$ is then defined as
\[
\mathcal{T}_h^{n+1} = \Phi_h(\hat{\mathcal{T}}_{c,h}),
\]
which can be computed readily using linear interpolation. This procedure is explained in Fig. \ref{fig:MMPDE}.

We now describe the MMPDE moving mesh method \cite{HRR94a, HRR94b} for minimizing $I_h$.
The MMPDE is defined as a gradient system of $I_h$, i.e., 
\begin{equation}
\frac{d\xi_j}{dt} = -\frac{P_j}{\tau} \left( \frac{\partial I_h}{\partial \xi_j} \right)^T, \quad j = 1,...,N_v
\label{MMPDE1}
\end{equation}
where $\tau > 0$ is a parameter used to adjust the time scale of mesh movement
and $P_j = \det(\mathbb{M}(x_j))^{\frac{p-1}{2}}$ is chosen such that (\ref{MMPDE1}) is invariant under
the scaling transformation of $\mathbb{M}$.
The derivative of $I_h$ with respect to $\xi_j$ is considered as a row vector and can be found analytically by the notion of scale-by-matrix differentiation; see \cite{HK15a}. Using the analytical formula, we can rewrite (\ref{MMPDE1}) as
\begin{equation}
\frac{d\xi_j}{dt} = \frac{P_j}{\tau}\sum_{K\in\omega_j} |K| {v}_{j_K}^K, 
\label{MMPDE2}
\end{equation}
where $\omega_j$ is the element patch associated with the $j$-th vertex, $j_K$ is its local index of the vertex
on $K$, and the local velocities ${v}_{j_K}^K$ are given by
\begin{displaymath}
\left[ \begin{array}{ccc}
({v}_1^K)^T \\
({v}_2^K)^T
\end{array} \right] = 
-E_K^{-1} \frac{\partial G}{\partial \mathbb{J}} - \frac{\partial G}{\partial \det(\mathbb{J})} \frac{\det(\hat{E}_K)}{\det(E_K)} \hat{E}_K^{-1},
\quad {v}_0^K = - \sum_{i=1}^2 {v}_j^K,
\end{displaymath}
where the derivatives of the function $G$ are given by
\begin{align*}
& \frac{\partial G}{\partial \mathbb{J}} = 2 p\theta \sqrt{\det(\mathbb{M})}\left( \tr(\mathbb{J}\mathbb{M}^{-1}\mathbb{J}^T) \right)^{p-1}\mathbb{M}^{-1}\mathbb{J}^T, \\
& \frac{\partial G}{\partial \det(\mathbb{J})} = p(1-2\theta) 2^{p} \det (\mathbb{M})^{\frac{1-p}{2}} \det(\mathbb{J})^{p-1}.
\end{align*}
MMPDE (\ref{MMPDE2}), with proper modifications for the boundary vertices to allow them only to slide
on the boundary, is integrated from $t^n$ to $t^{n+1}$ with the initial mesh $\hat{\mathcal{T}}_{c,h}$.
Notice that the integration of the MMPDE is equivalent to performing steepest descent for minimizing $I_h$.
In our computation, we use the Matlab\textsuperscript \textregistered\, function {\em ode15s},
a Numerical Differentiation Formula based integrator, for the integration.

\section{Numerical results}
\label{SEC:numerics}

In this section we present numerical results obtained with the moving mesh finite element method described
in the previous section for three examples. The first two examples are benchmark problems commonly
used in the existing literature to examine mathematical models for brittle fracture and related
numerical algorithms. The last example is chosen to test the ability of our method to handle multiple and complex cracks.
Special attention is paid to the demonstration of the effectiveness of the MMPDE method to track
crack propagation and the effects of the regularization methods on the convergence of Newton's iteration.
In the results presented in this section, an adaptive mesh of size $N = 6,400$ and
the sonic-point regularization method with $\alpha = 1 \times 10^{-3}$ are used, unless stated otherwise.

\subsection{Example 1. Single edge notched tension test}

\begin{figure} [!htb]
\centering 
\subfigure[tension test]{\label{fig:subfig:tension}
\includegraphics[width=0.4\linewidth]{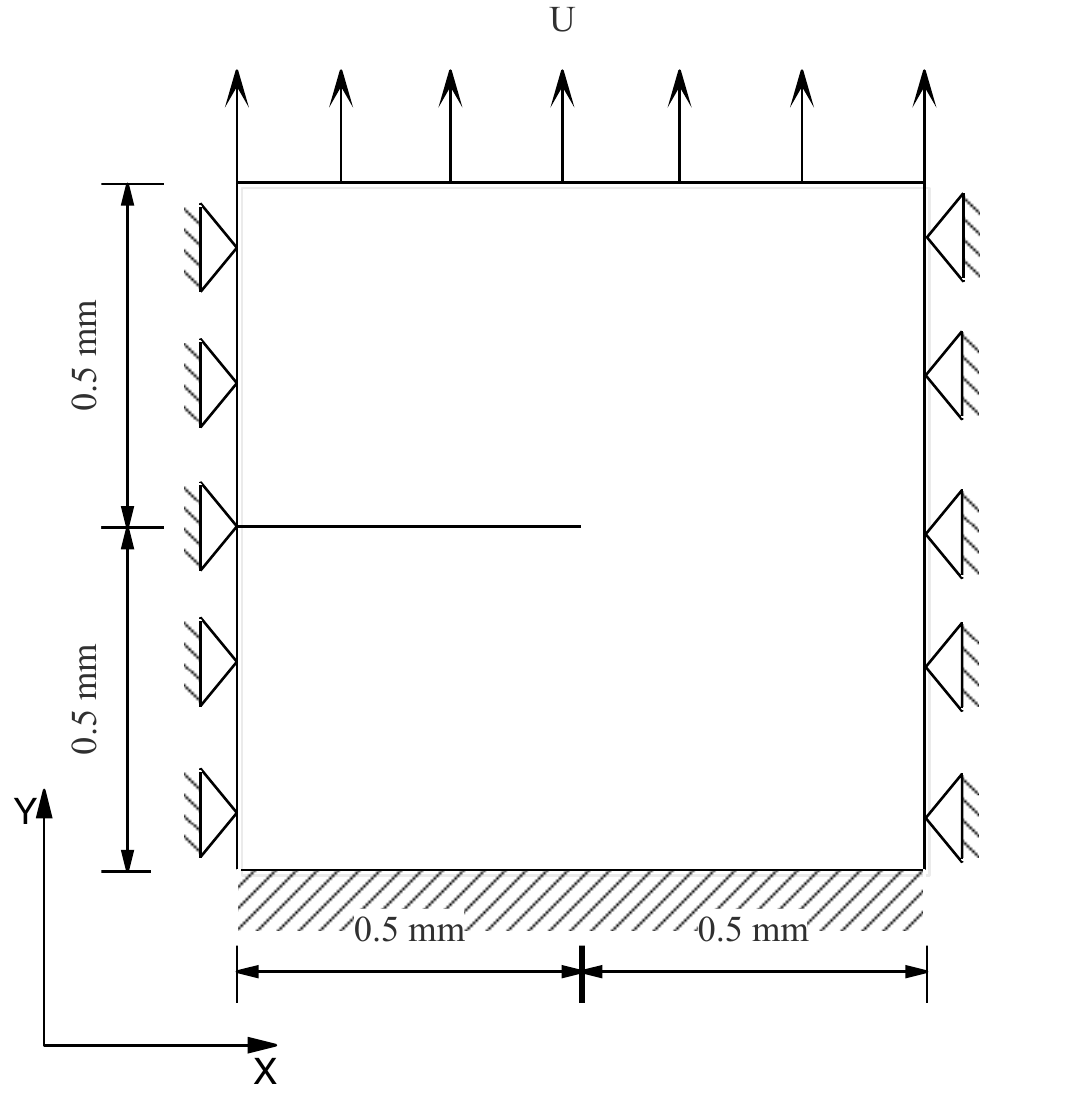}}
\subfigure[shear test]{\label{fig:subfig:shear}
\includegraphics[width=0.4\linewidth]{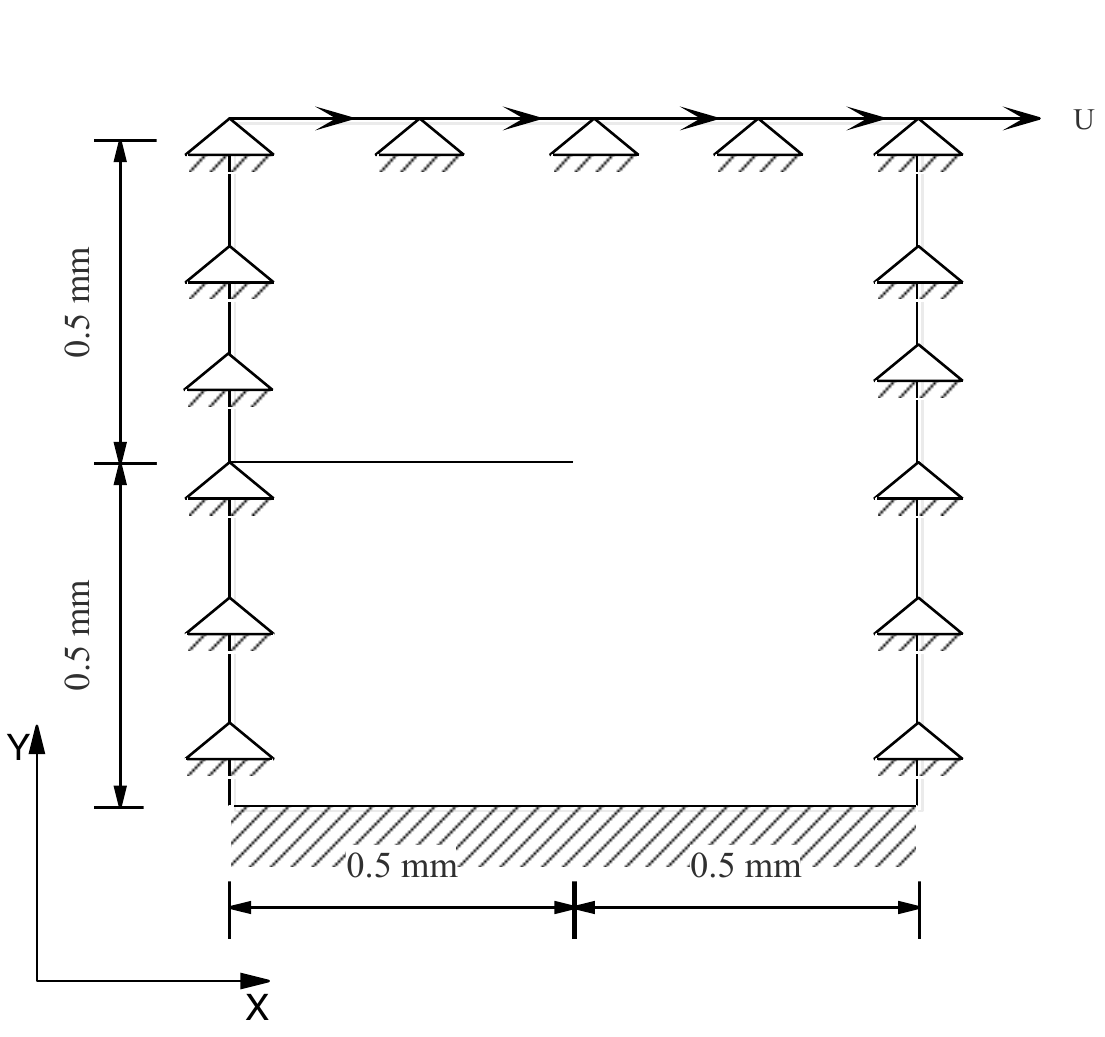}}
\caption{Domain and boundary conditions for single edge notched tests, (a) tension test for Example 1
and (b) shear test for Example 2.}
\label{fig:TS specimen}
\end{figure}

We first consider a single edge notched tension test from Miehe et al. \cite{MHW10}, with the domain and boundary conditions shown in Fig. \ref{fig:subfig:tension}.  For the boundary conditions, the bottom edge of the domain is fixed and the top edge is fixed along $x$-direction while a uniform $y$-displacement $U$ is increased with time to drive the crack propagation. 

The solid is assumed to be homogeneous isotropic with elastic bulk modulus $\lambda = 121.15$~kN/mm$^{2}$
and shear modulus $\mu = 80.77$~kN/mm$^{2}$. The fracture toughness is $g_c = 2.7 \times 10^{-3}$~kN/mm.
Two displacement increments have been prescribed for the computation, $\Delta U = 1 \times 10^{-5}$~mm is chosen
for the first 500 time steps and $\Delta U = 1 \times 10^{-6}$~mm afterwards. We consider two values of $l$,
$0.00375$~mm and $0.0075$~mm. An initial triangular mesh is constructed from a rectangular mesh
by subdividing each rectangle into four triangles along the diagonal directions. 

\begin{figure} 
\centering 
\subfigure[$U = 1.0 \times 10^{-3}$ mm]{\label{fig:subfig:TM_l75_u1_0}
\includegraphics[width=0.25\linewidth]{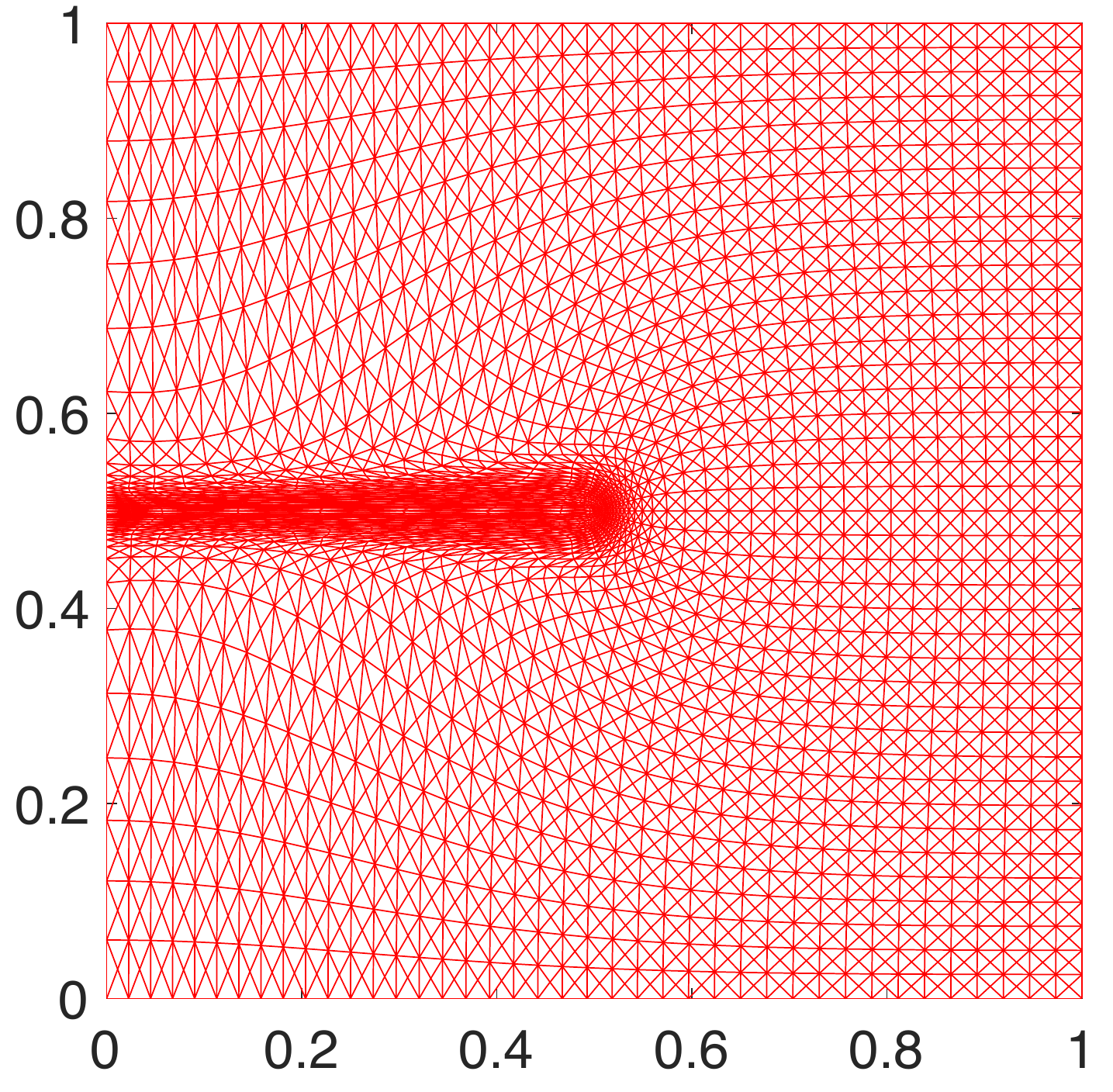}}
\subfigure[$U = 5.2 \times 10^{-3}$ mm]{\label{fig:subfig:TM_l75_u5_2}
\includegraphics[width=0.25\linewidth]{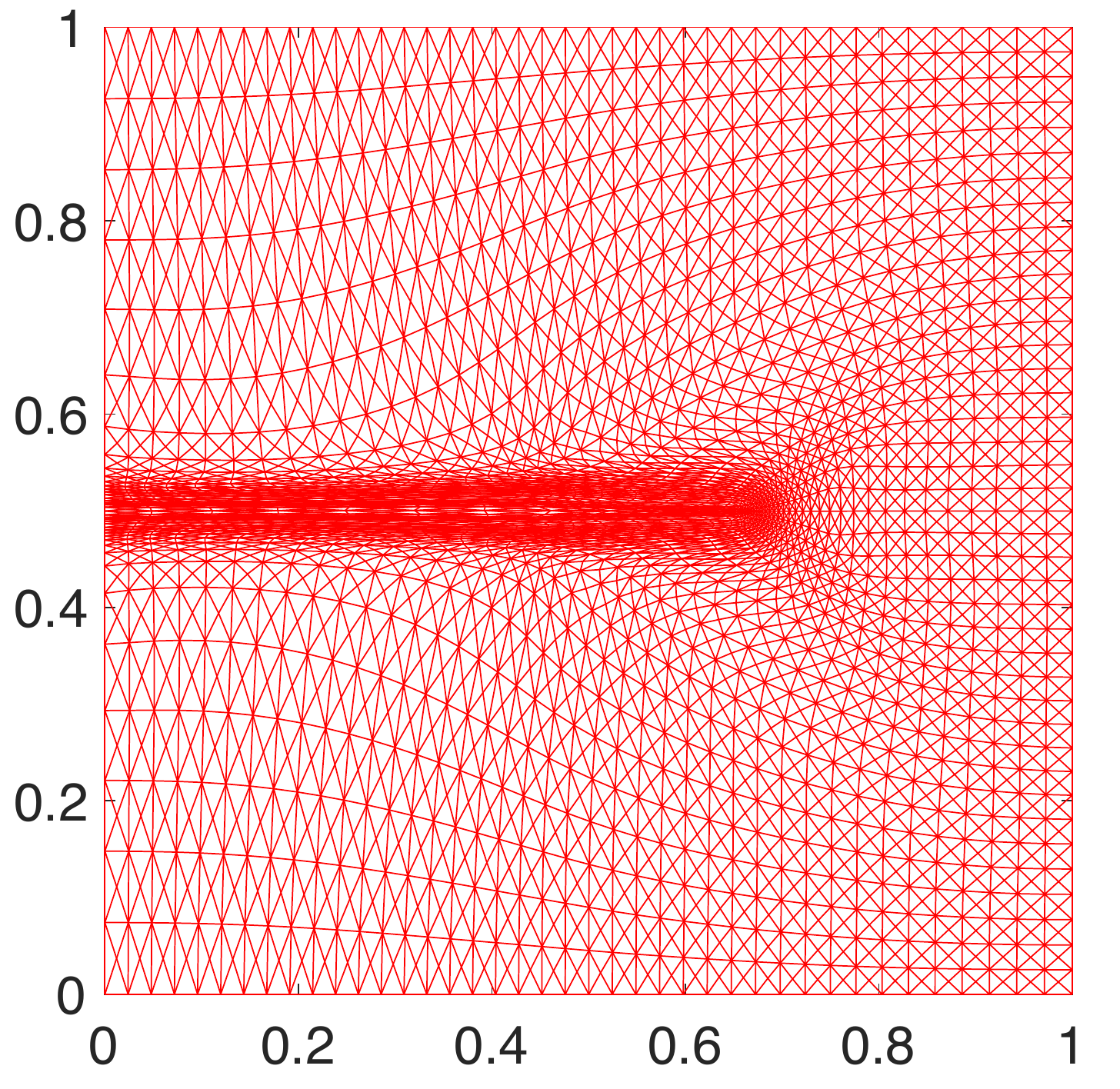}}
\subfigure[$U = 5.5 \times 10^{-3}$ mm]{\label{fig:subfig:TM_l75_u5_5}
\includegraphics[width=0.25\linewidth]{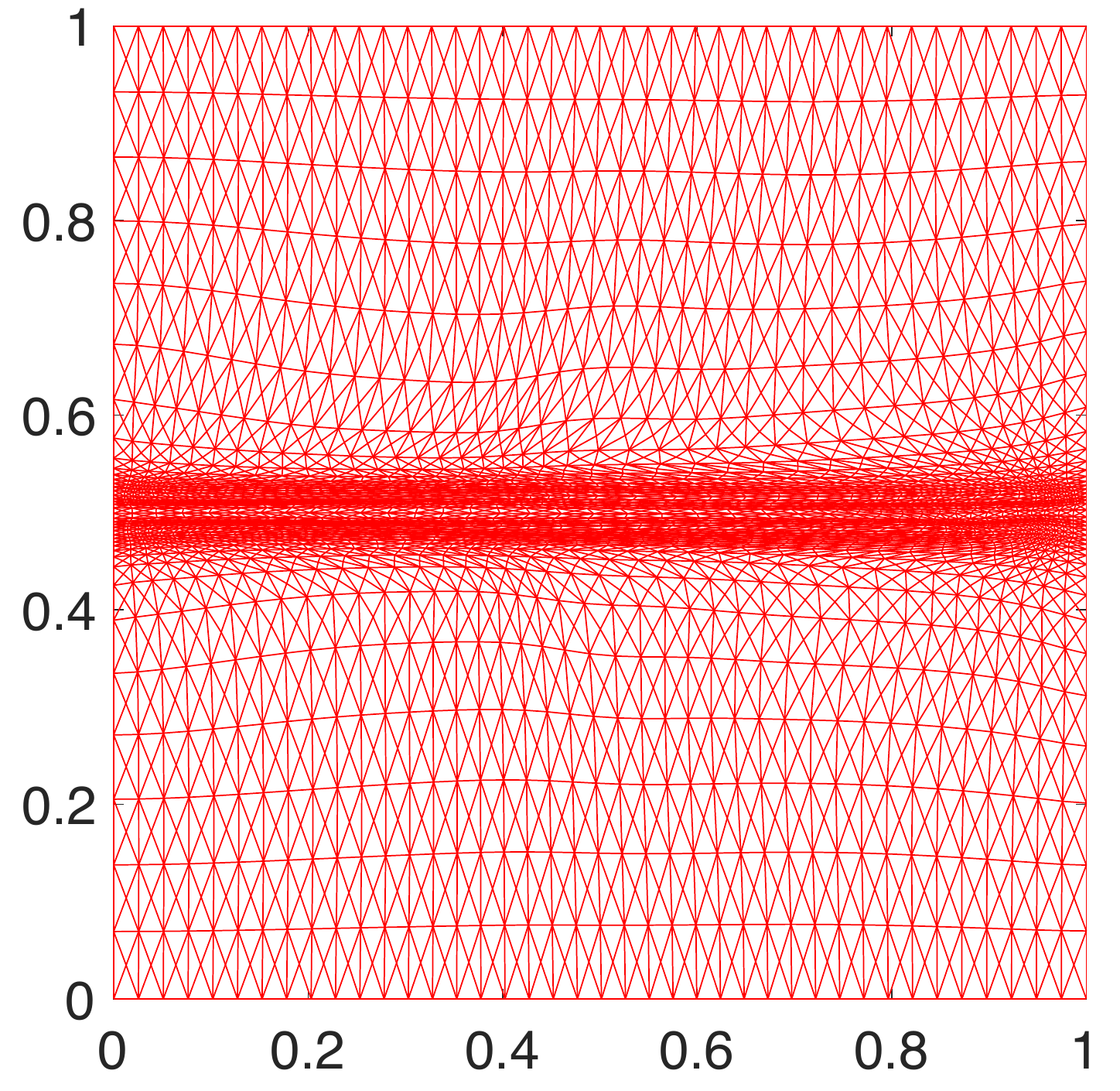}}
\vfill
\subfigure[$U = 1.0 \times 10^{-3}$ mm]{\label{fig:subfig:TD_l75_u1_0}
\includegraphics[width=0.25\linewidth]{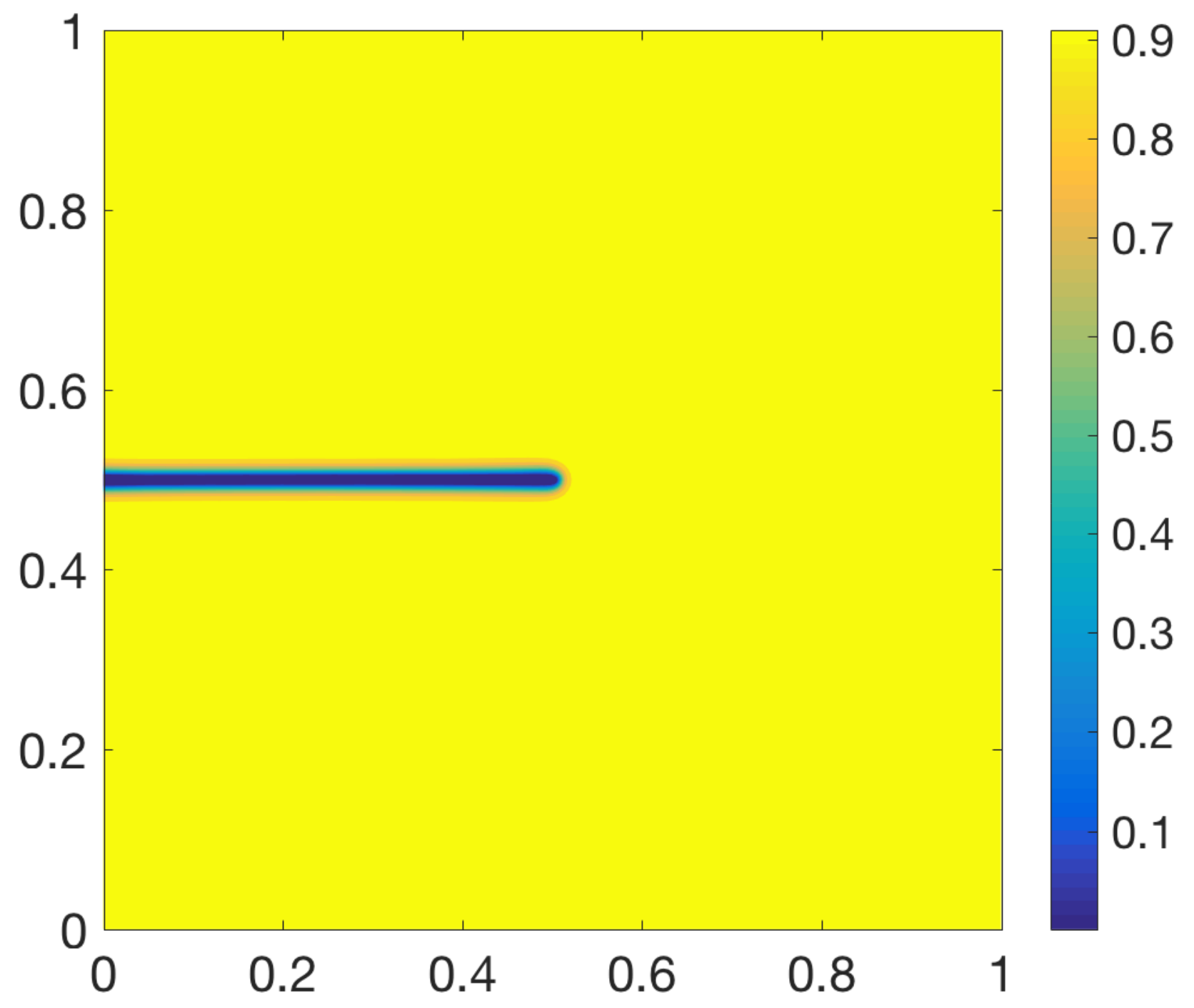}}
\subfigure[$U = 5.2 \times 10^{-3}$ mm]{\label{fig:subfig:TD_l75_u5_2}
\includegraphics[width=0.25\linewidth]{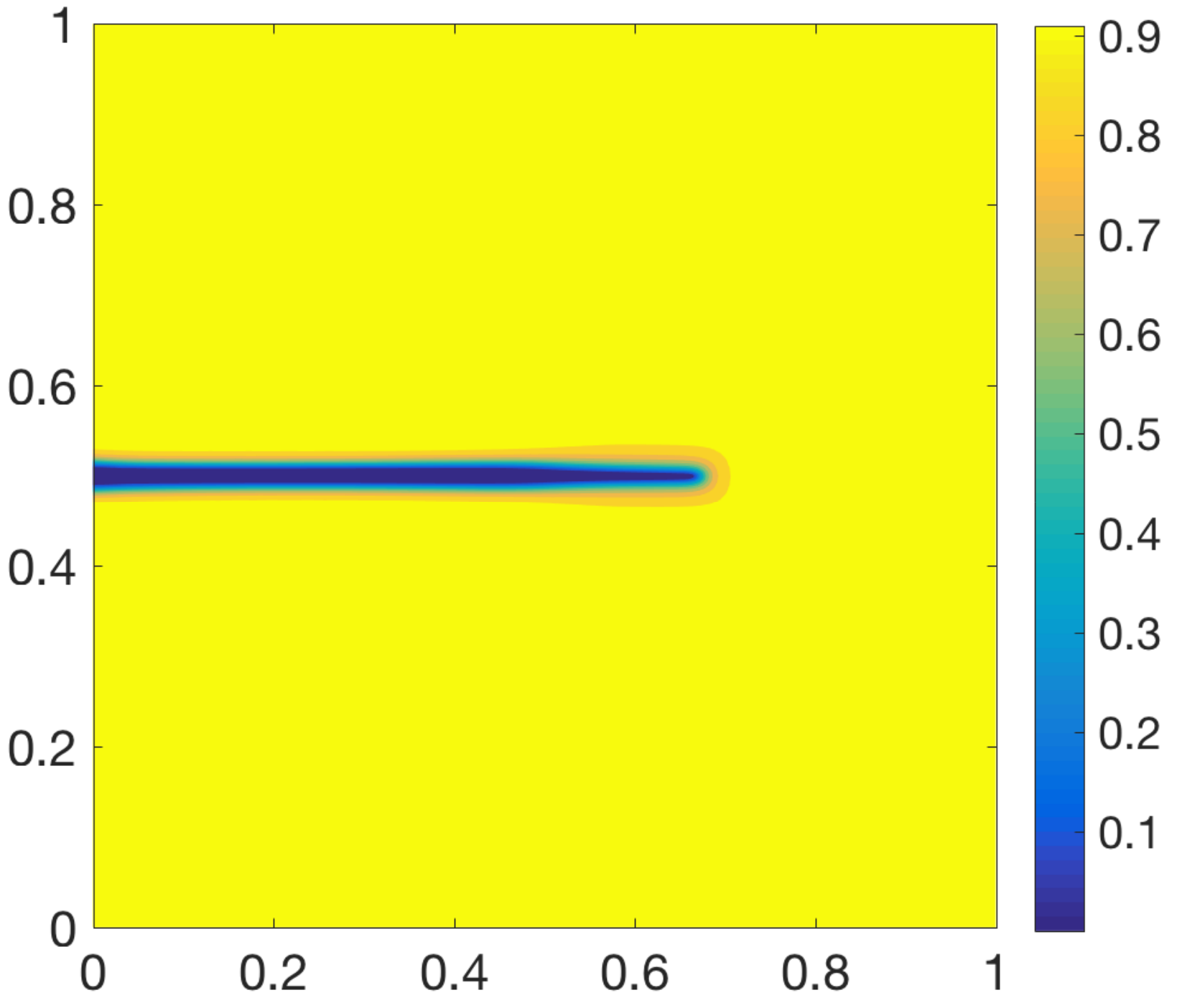}}
\subfigure[$U = 5.5 \times 10^{-3}$ mm]{\label{fig:subfig:TD_l75_u5_5}
\includegraphics[width=0.25\linewidth]{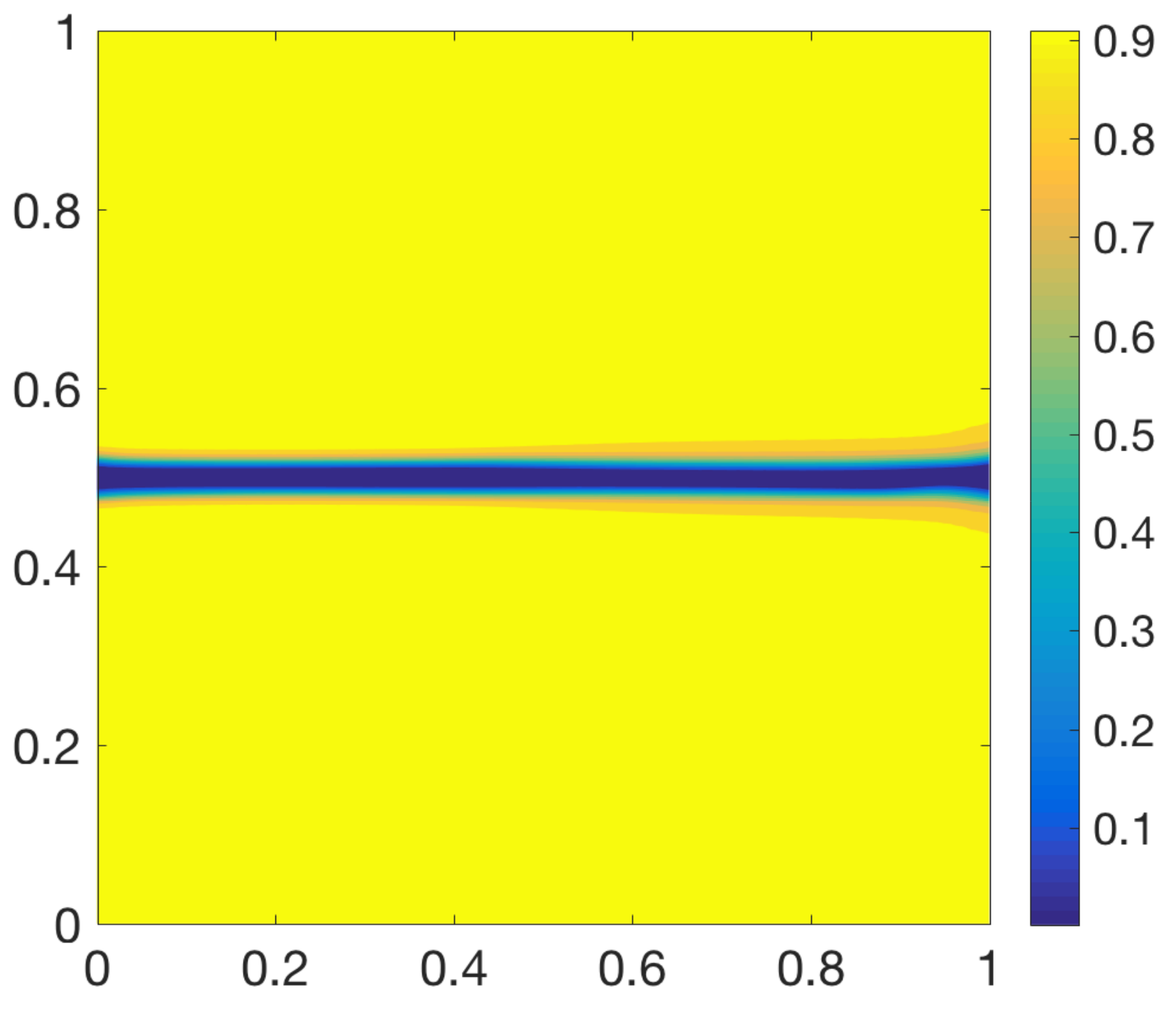}}
\caption{Example 1. The mesh and contours of the phase-field distribution during crack evolution
for the tension test with $l = 0.00375$~mm. ($N = 6,400$)}
\label{fig:l75_tension}
\end{figure}

\begin{figure} 
\centering 
\subfigure[$U = 1.0 \times 10^{-3}$ mm]{\label{fig:subfig:TM_l15_u1_0}
\includegraphics[width=0.25\linewidth]{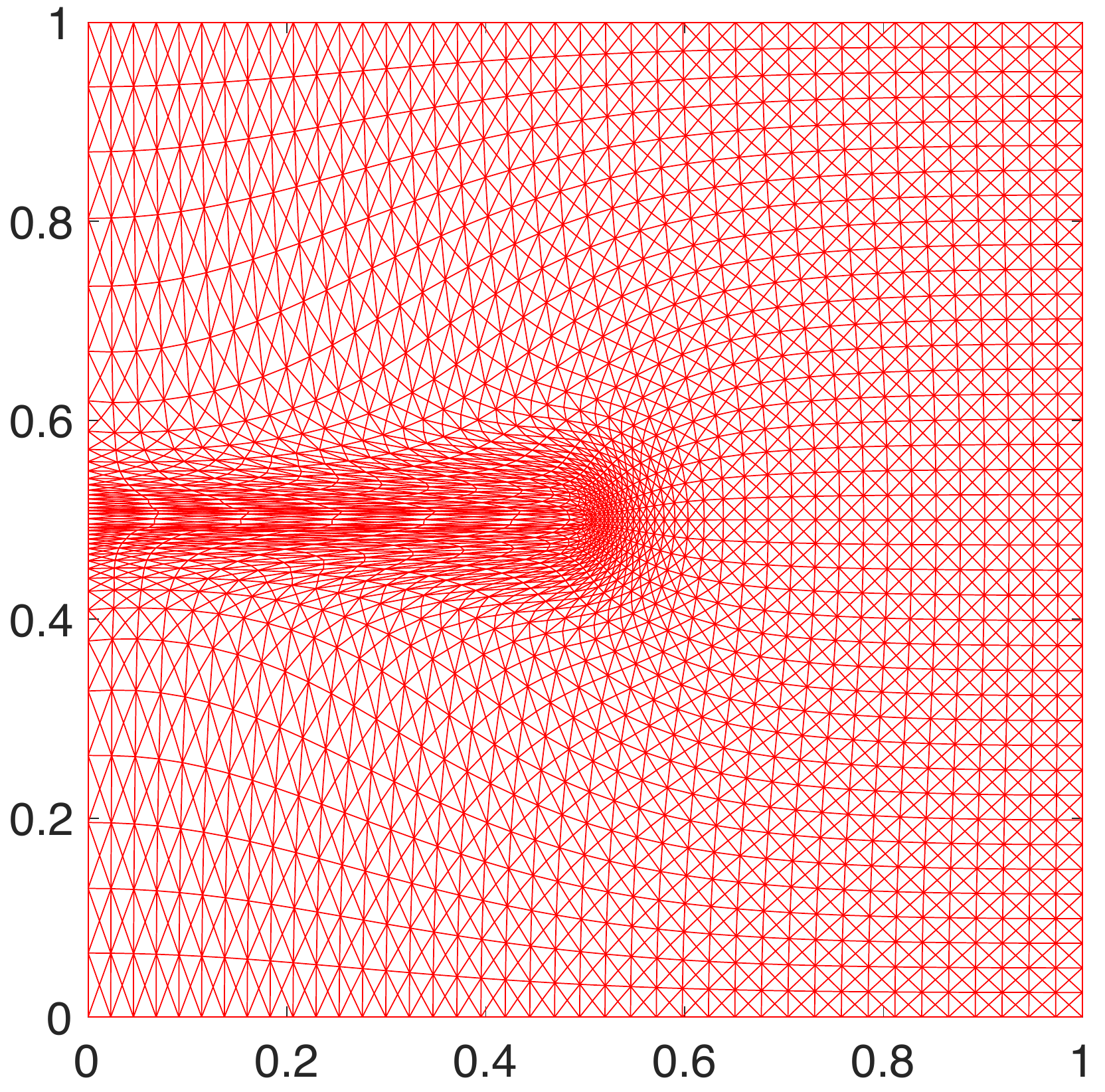}}
\subfigure[$U = 5.2 \times 10^{-3}$ mm]{\label{fig:subfig:TM_l15_u5_2}
\includegraphics[width=0.25\linewidth]{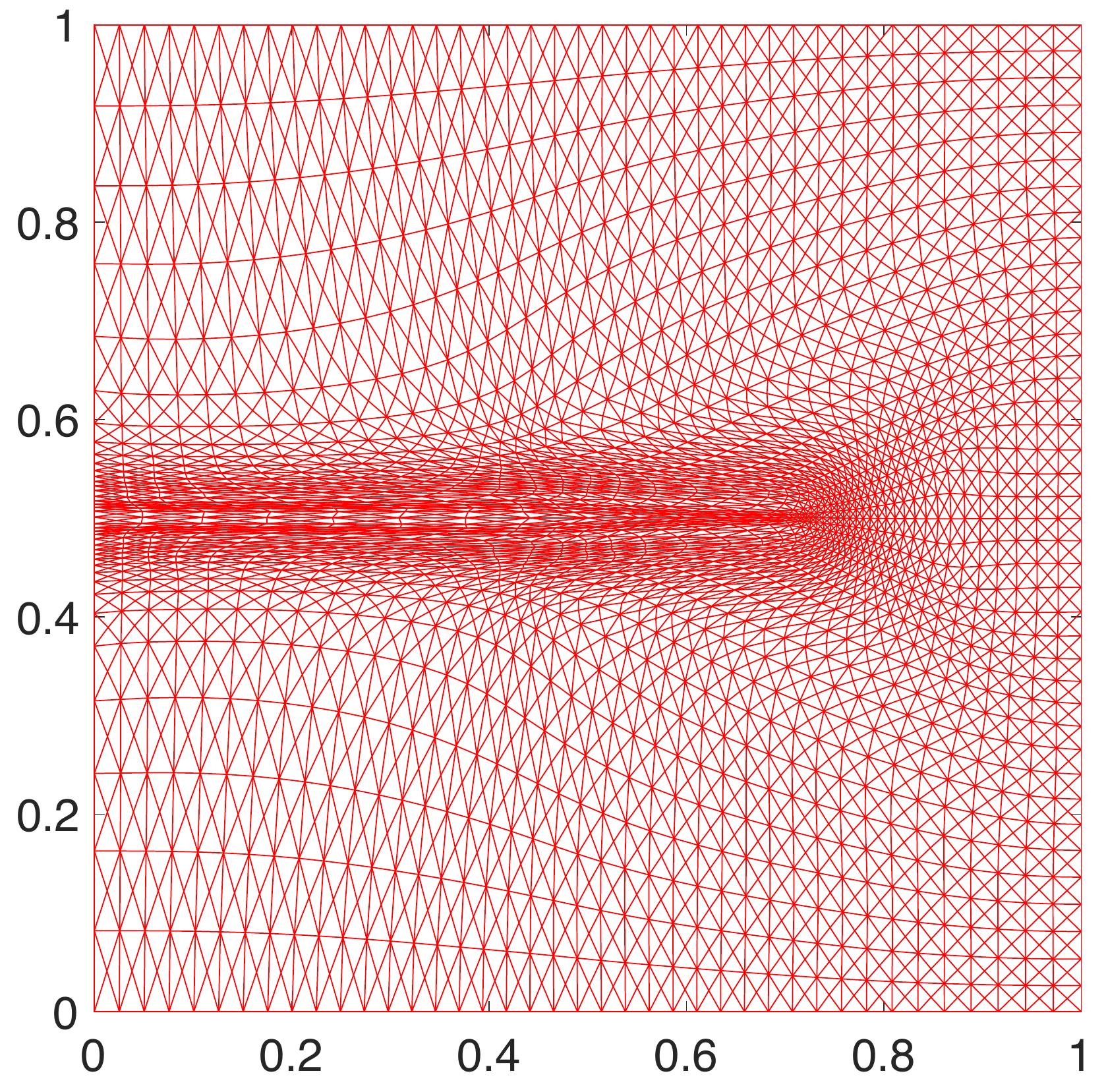}}
\subfigure[$U = 5.5 \times 10^{-3}$ mm]{\label{fig:subfig:TM_l15_u5_5}
\includegraphics[width=0.25\linewidth]{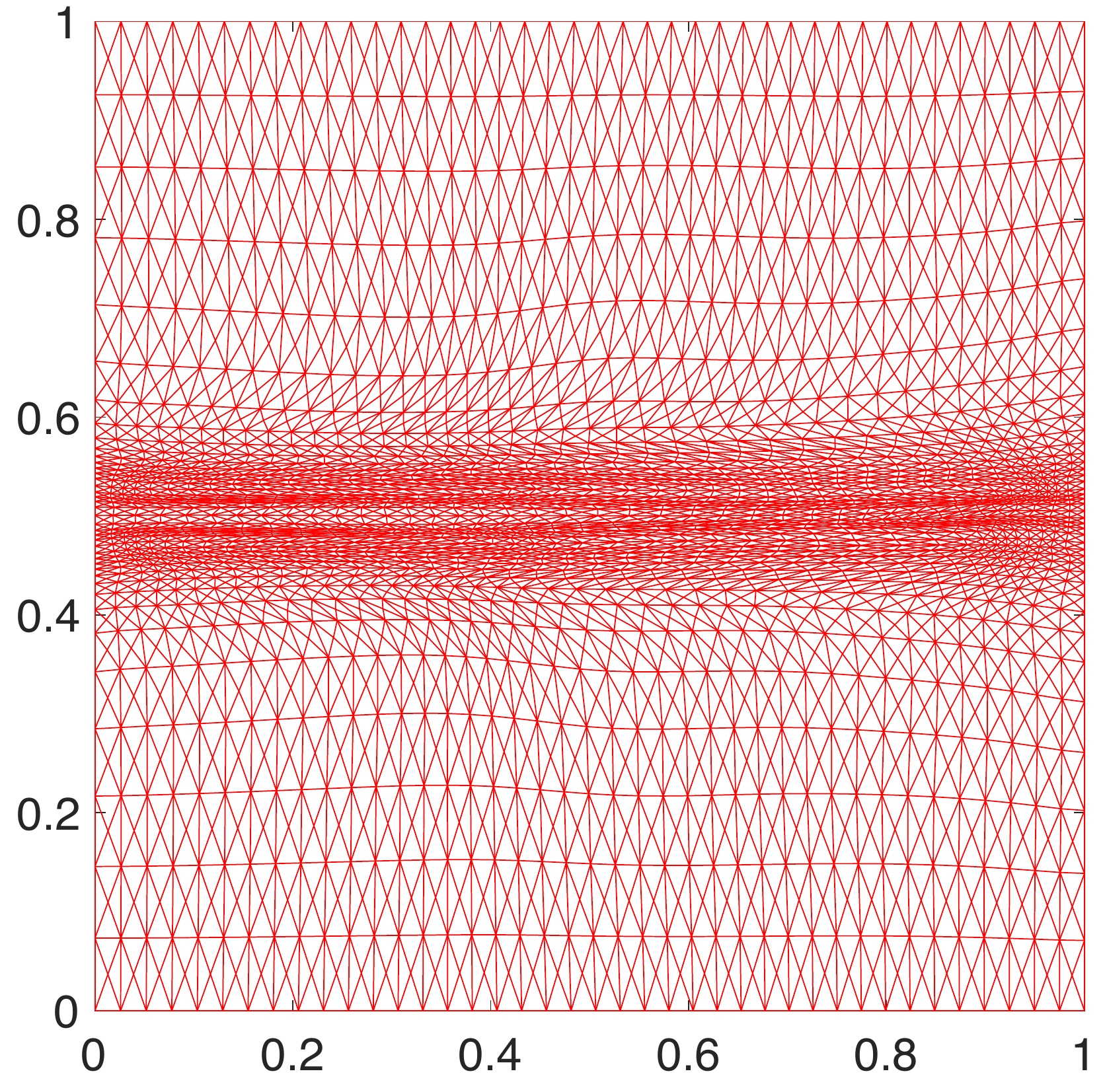}}
\vfill
\subfigure[$U = 1.0 \times 10^{-3}$ mm]{\label{fig:subfig:TD_l15_u1_0}
\includegraphics[width=0.25\linewidth]{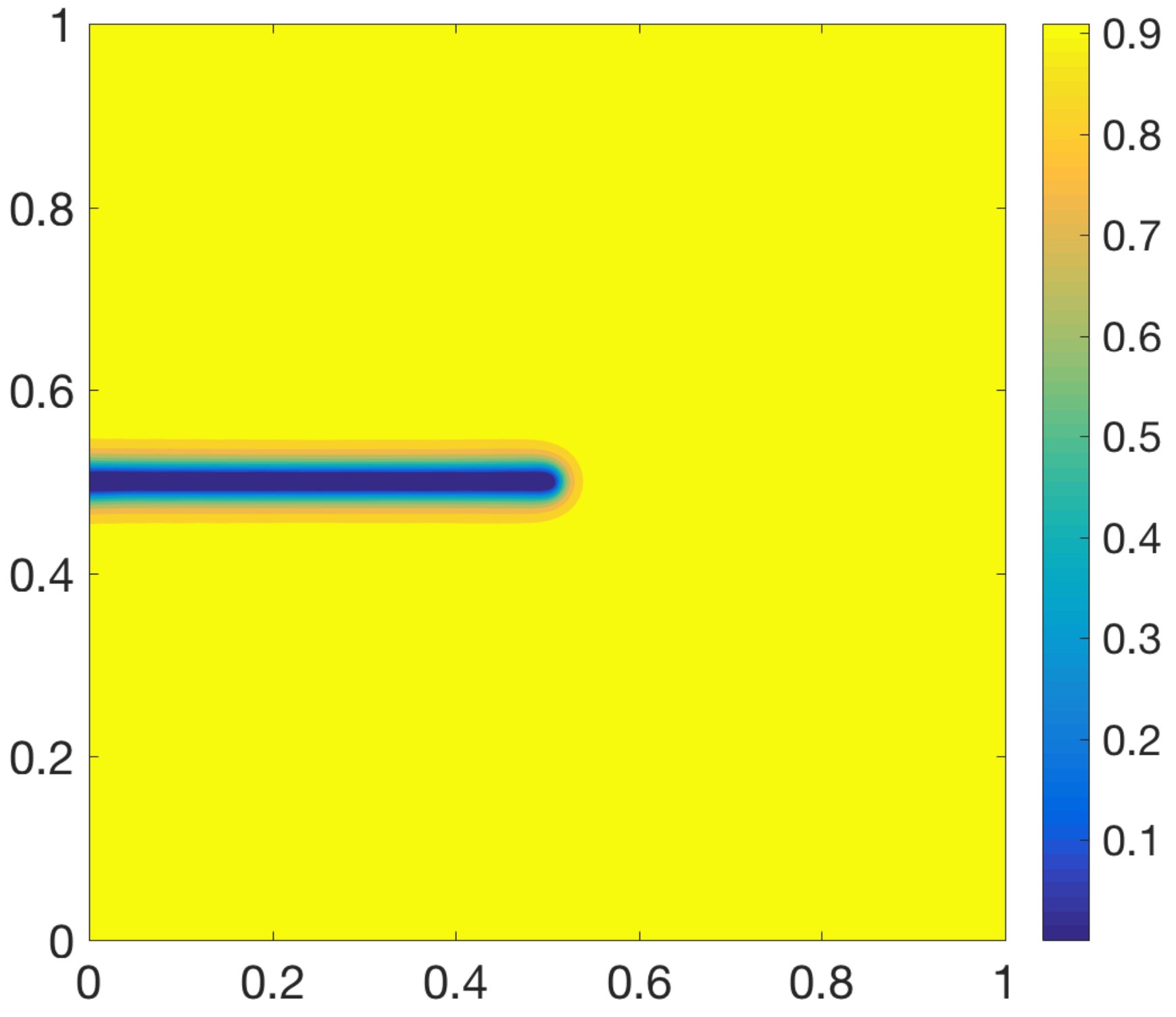}}
\subfigure[$U = 5.2 \times 10^{-3}$ mm]{\label{fig:subfig:TD_l15_u5_2}
\includegraphics[width=0.25\linewidth]{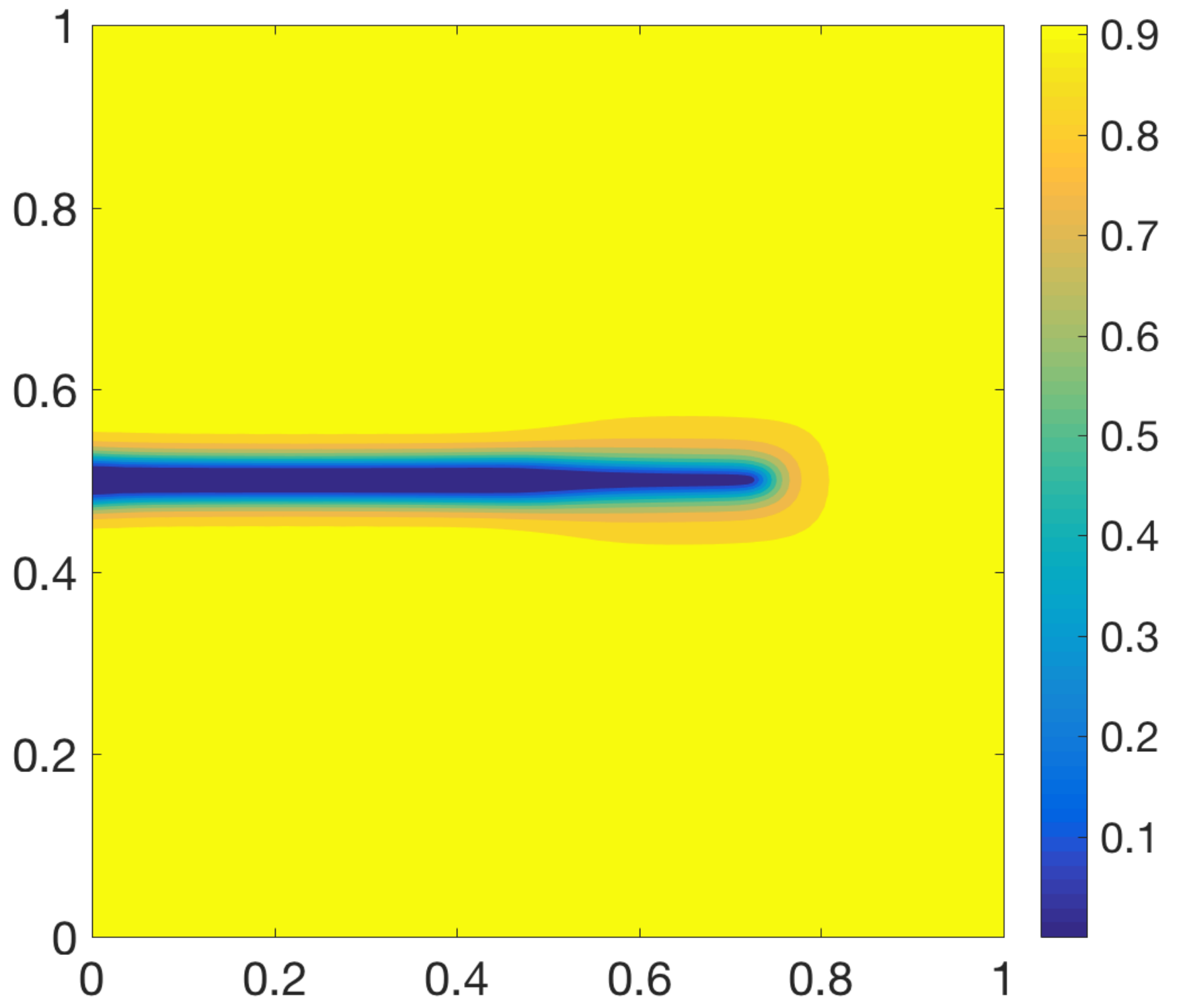}}
\subfigure[$U = 5.5 \times 10^{-3}$ mm]{\label{fig:subfig:TD_l15_u5_5}
\includegraphics[width=0.25\linewidth]{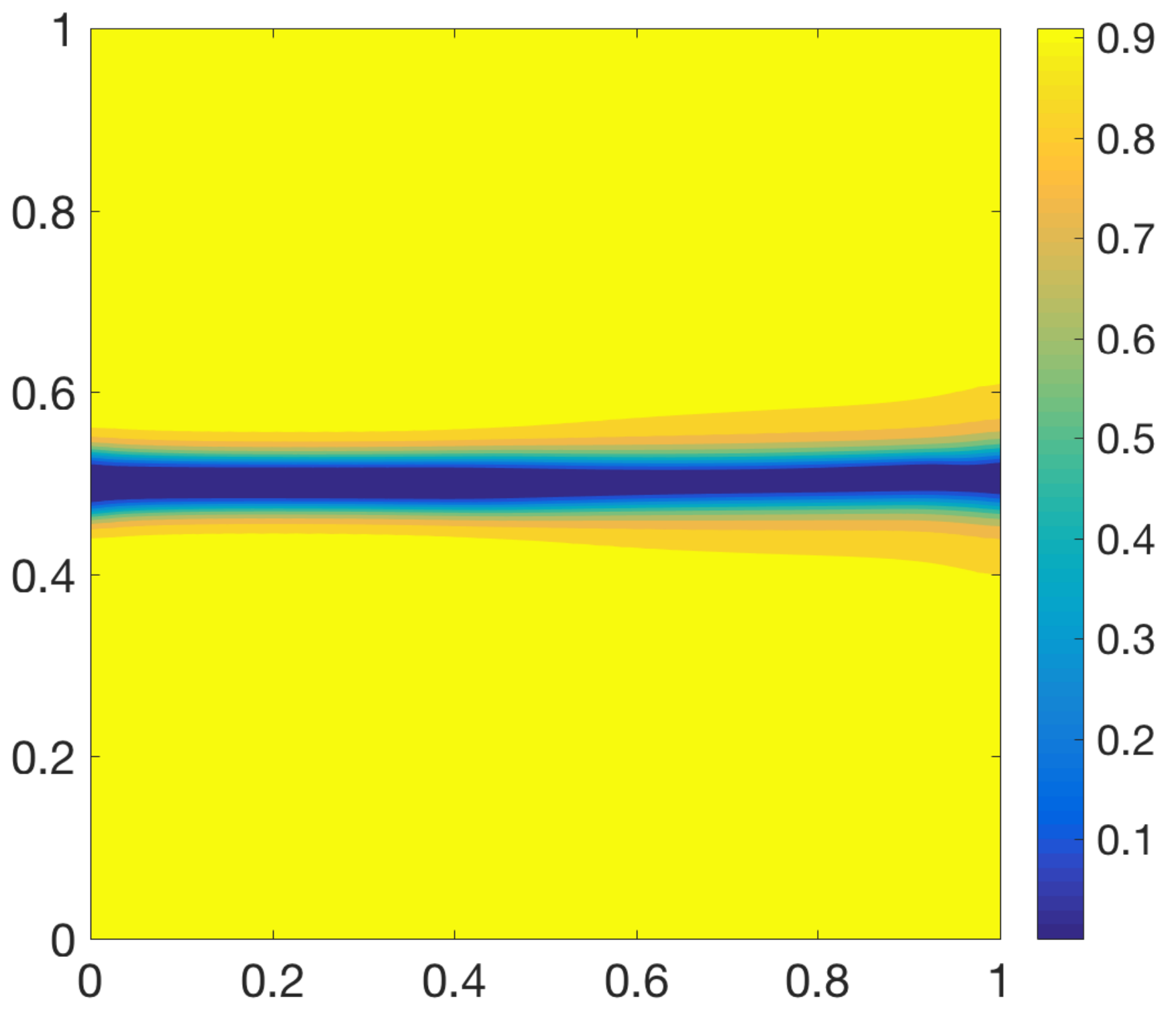}}
\caption{Example 1. The mesh and contours of the phase-field distribution during
crack evolution for the tension test with $l = 0.0075$~mm. ($N = 6,400$)}
\label{fig:l15_tension}
\end{figure}

\begin{figure} [!htb]
\centering 
\subfigure[adaptive mesh]{\label{fig:subfig:TM_l75_u1_0_2}
\includegraphics[width=0.28\linewidth]{./TM_l75_u1_0.pdf}}
\subfigure[close view around the crack]{\label{fig:subfig:Tcrack_around}
\includegraphics[width=0.29\linewidth]{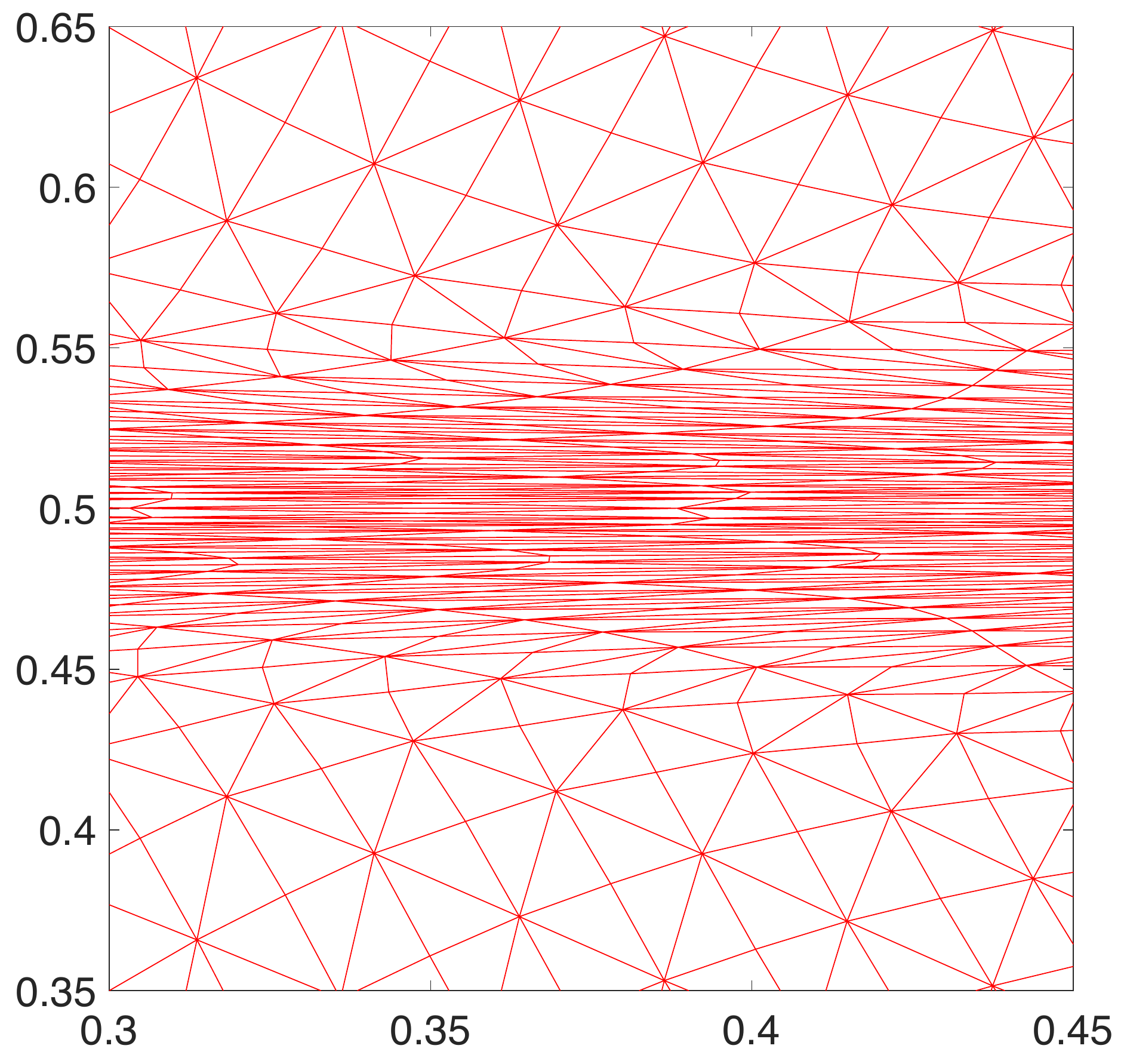}}
\subfigure[close view near the crack tip]{\label{fig:subfig:Tcrack_tip}
\includegraphics[width=0.29\linewidth]{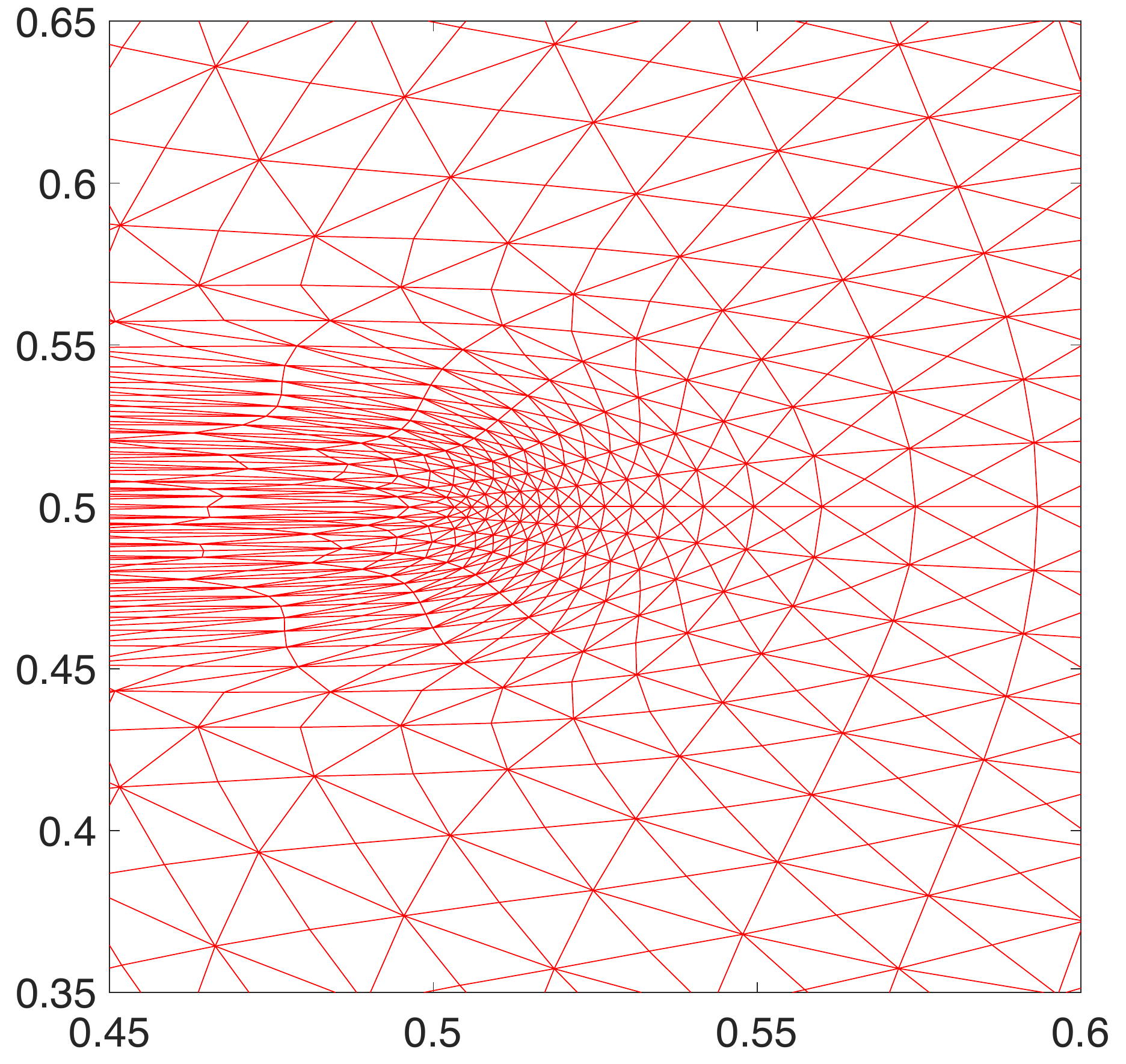}}
\caption{Example 1. An adaptive mesh of $N = 6,400$ and its close views near the crack
and crack tip ($l = 0.00375$~mm).}
\label{fig:initial-mesh}
\end{figure}

\begin{figure} [!htb]
\centering 
\subfigure[$l = 0.00375$ mm]{\label{fig:subfig:Tld_l75}
\includegraphics[width=0.4\linewidth]{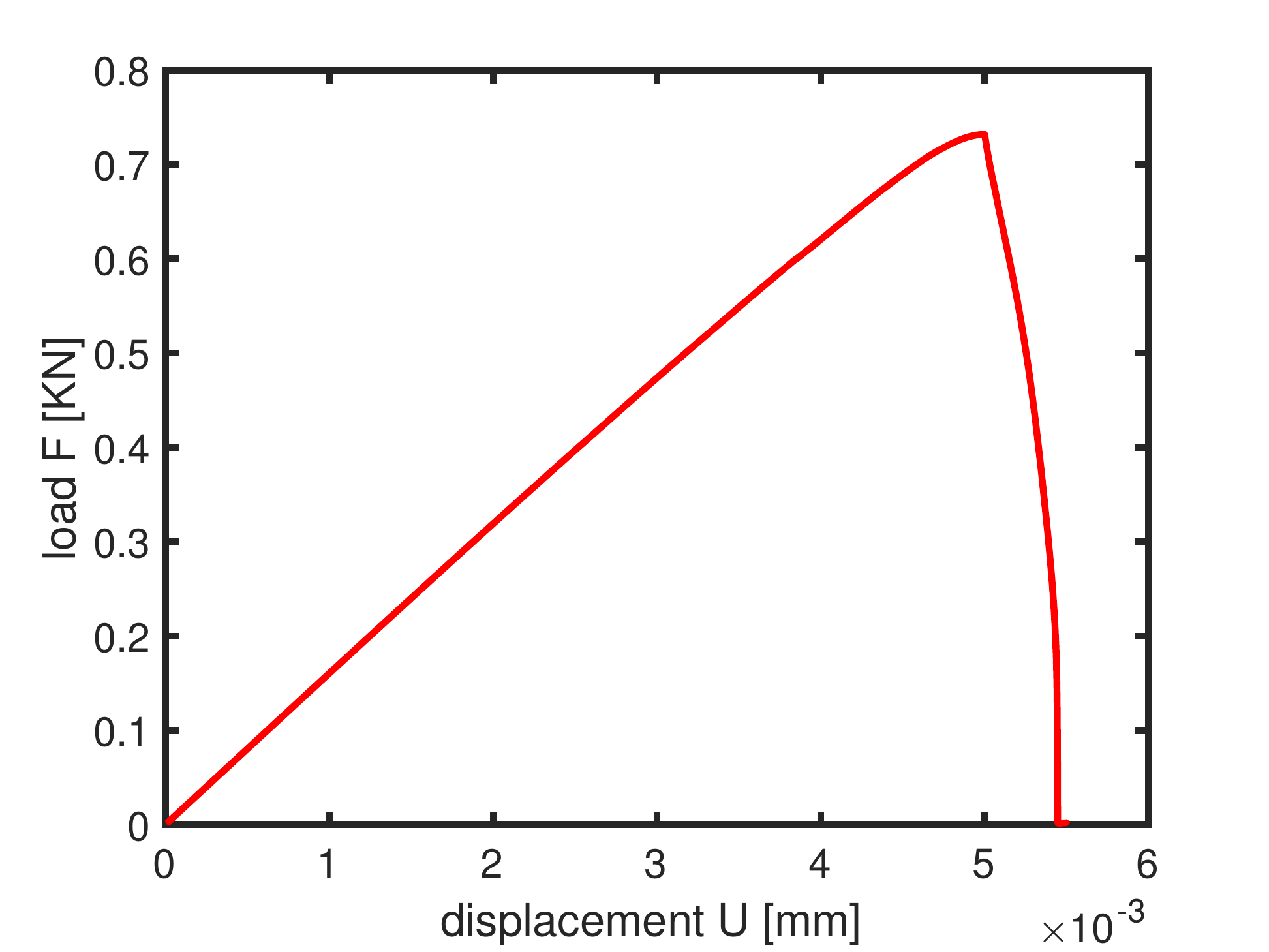}}
\subfigure[$l = 0.0075$ mm]{\label{fig:subfig:Tld_l15}
\includegraphics[width=0.4\linewidth]{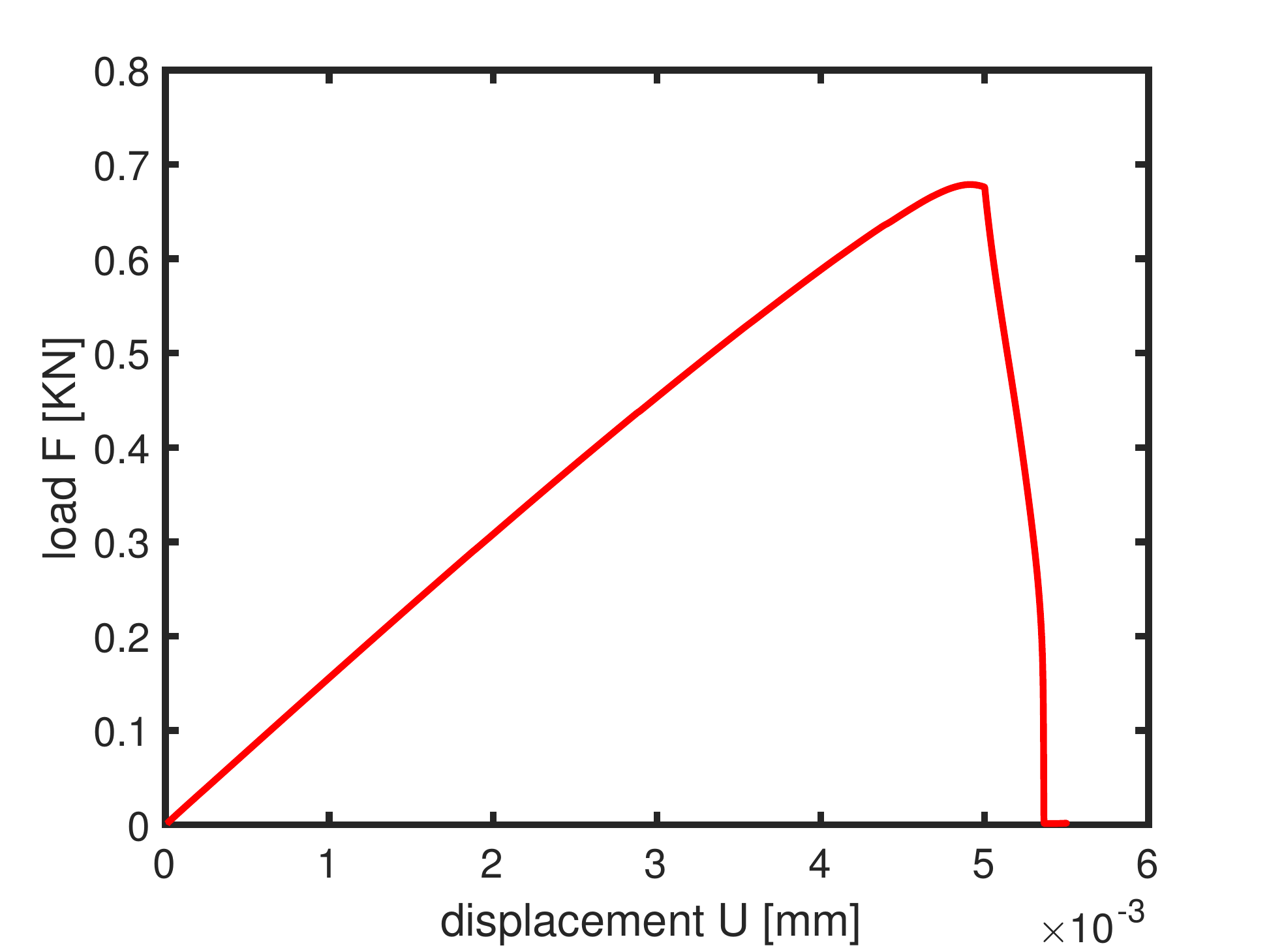}}
\caption{Example 1. The load-deflection curves for the tension test with different values of $l$. ($N = 6,400$)}
\label{fig:Tld_Diffl}
\end{figure}

Typical adaptive meshes and contours of the phase-field variable during crack evolution
for $l = 0.00375$~mm and $0.0075$~mm are shown in Fig. \ref{fig:l75_tension} and Fig. \ref{fig:l15_tension},
respectively. As can be seen, the mesh points concentrate in the region around the crack, which demonstrates
the effectiveness of the mesh adaptation strategy for this tension test. Closer views around the crack and crack
tip are shown in Fig. \ref{fig:subfig:Tcrack_around} and Fig. \ref{fig:subfig:Tcrack_tip}. It can be observed that
the mesh stays symmetric near the crack for all time.
For comparison purpose, we compute the surface load vector on the top edge as 
\[
F = (F_x,F_y):= \int_{\text{top edge}} \sigma(u) \cdot n \, d l ,
\]
where $n$ is the unit outward normal to the top edge.
Particularly, we are interested in $F_y$ for tension test and $F_x$ for shear test.
The load-deflection curves are shown in Fig. \ref{fig:Tld_Diffl}.

\subsubsection{Effects of the regularization methods} 
\label{SEC:effect}

\begin{figure} [!htb]
\centering 
\subfigure[No regularization at step 1]{\label{fig:subfig:non_regularization}
\includegraphics[width=0.34\linewidth]{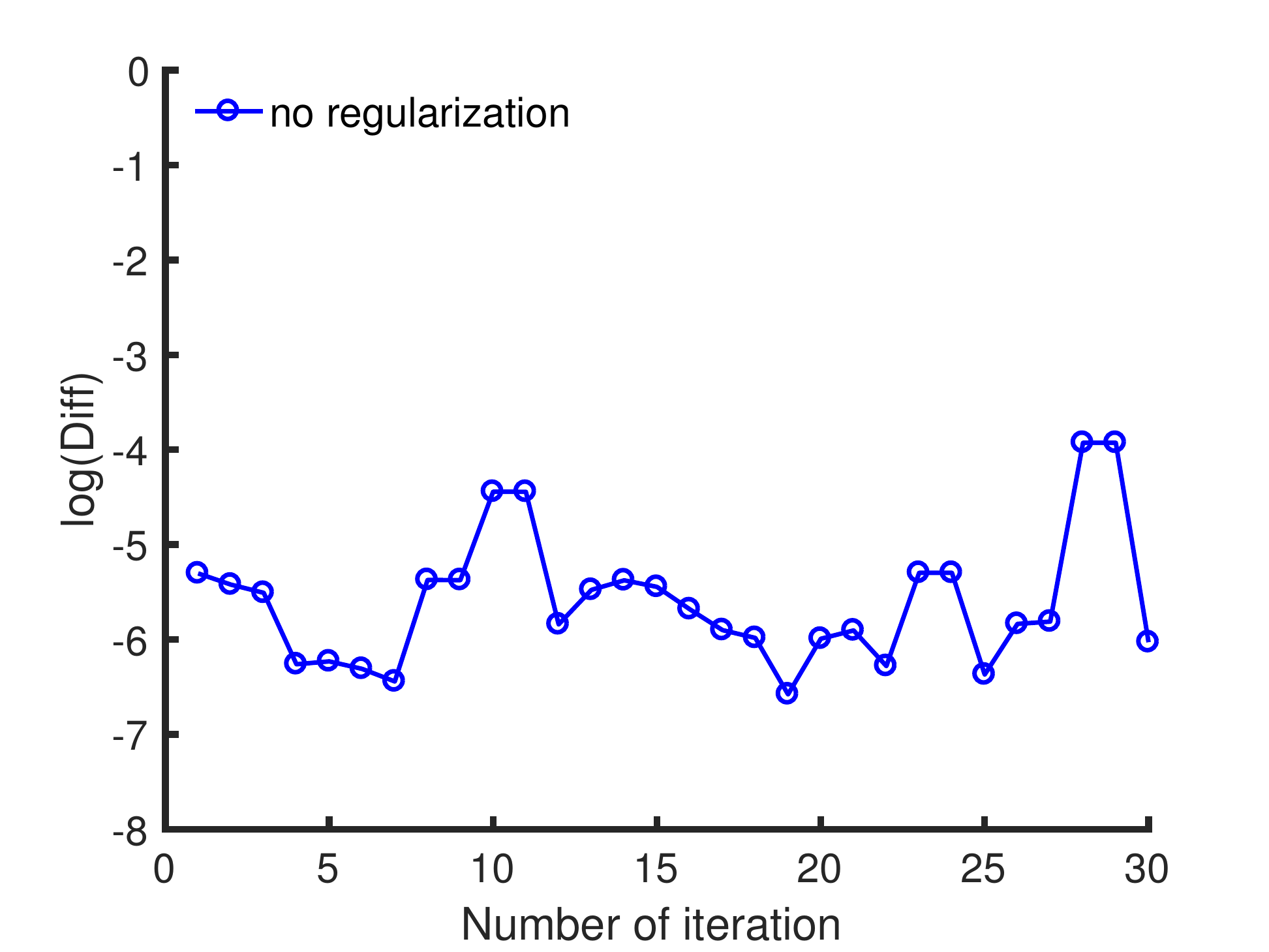}}
\subfigure[No regularization at step 2]{\label{fig:subfig:non_regularization_step2}
\includegraphics[width=0.34\linewidth]{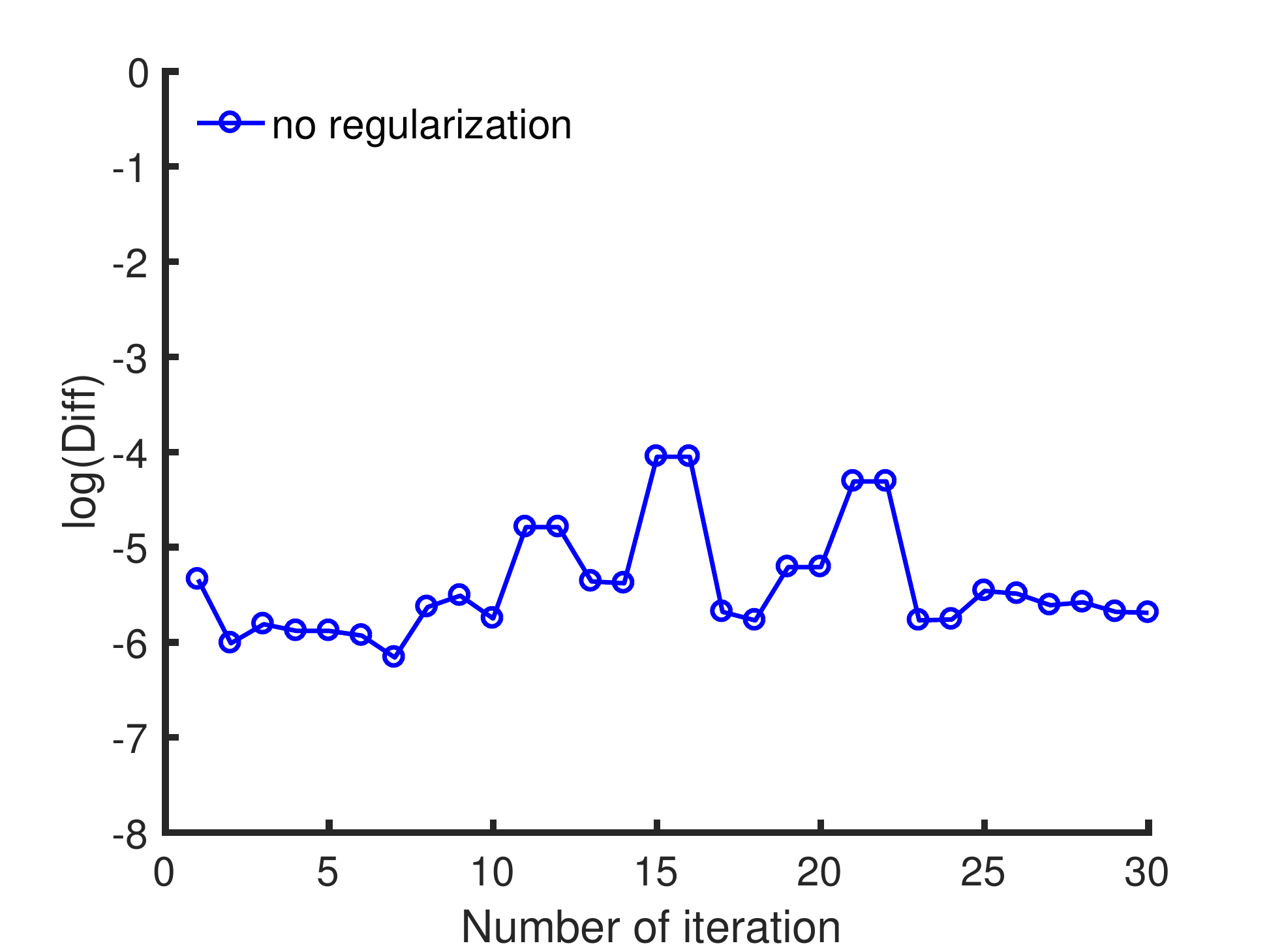}}
\vfill
\subfigure[Sonic-point at step 1]{\label{fig:subfig:sonic_point}
\includegraphics[width=0.34\linewidth]{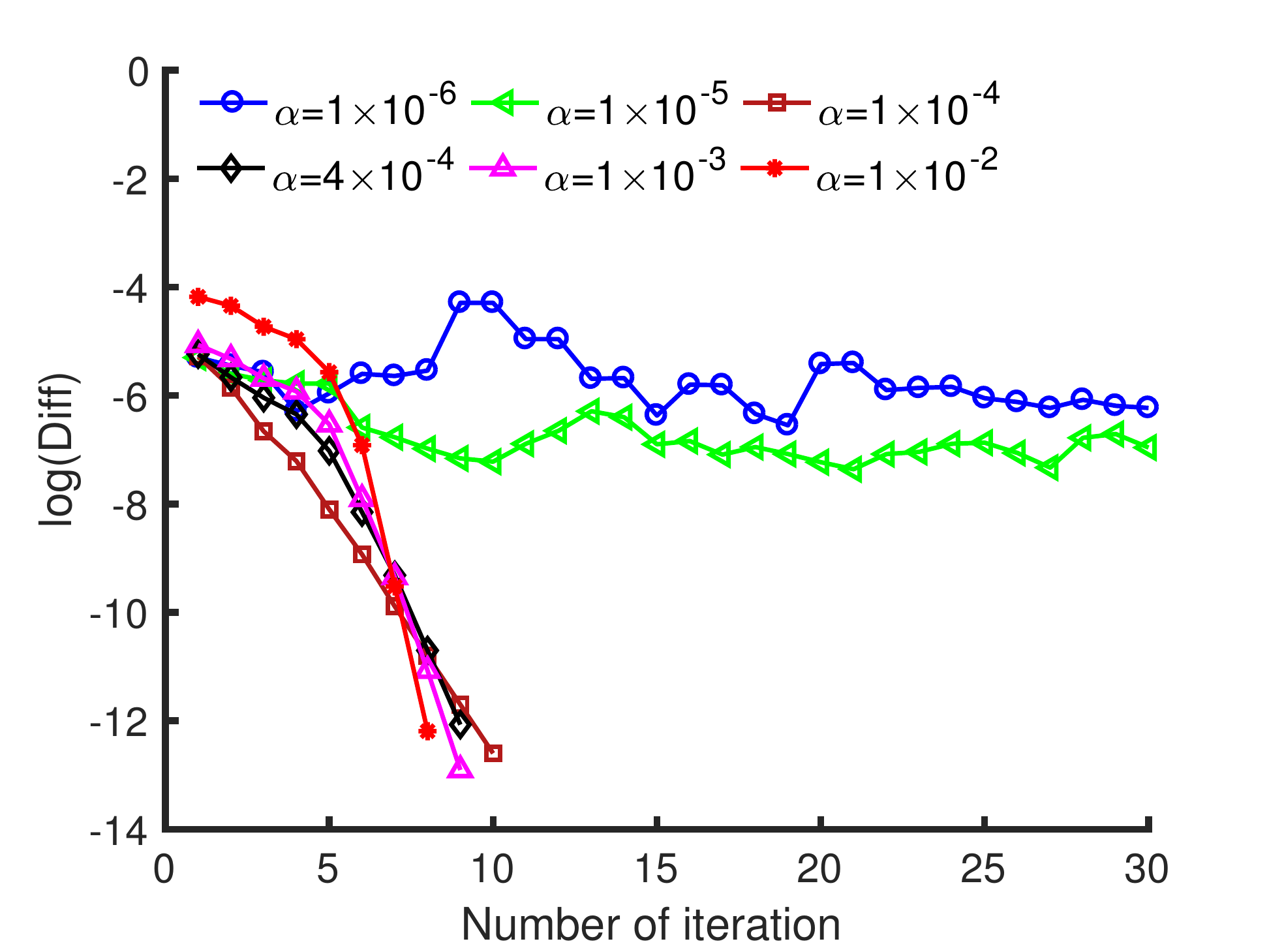}}
\subfigure[Sonic-point at step 2]{\label{fig:subfig:sonic_point_step2}
\includegraphics[width=0.34\linewidth]{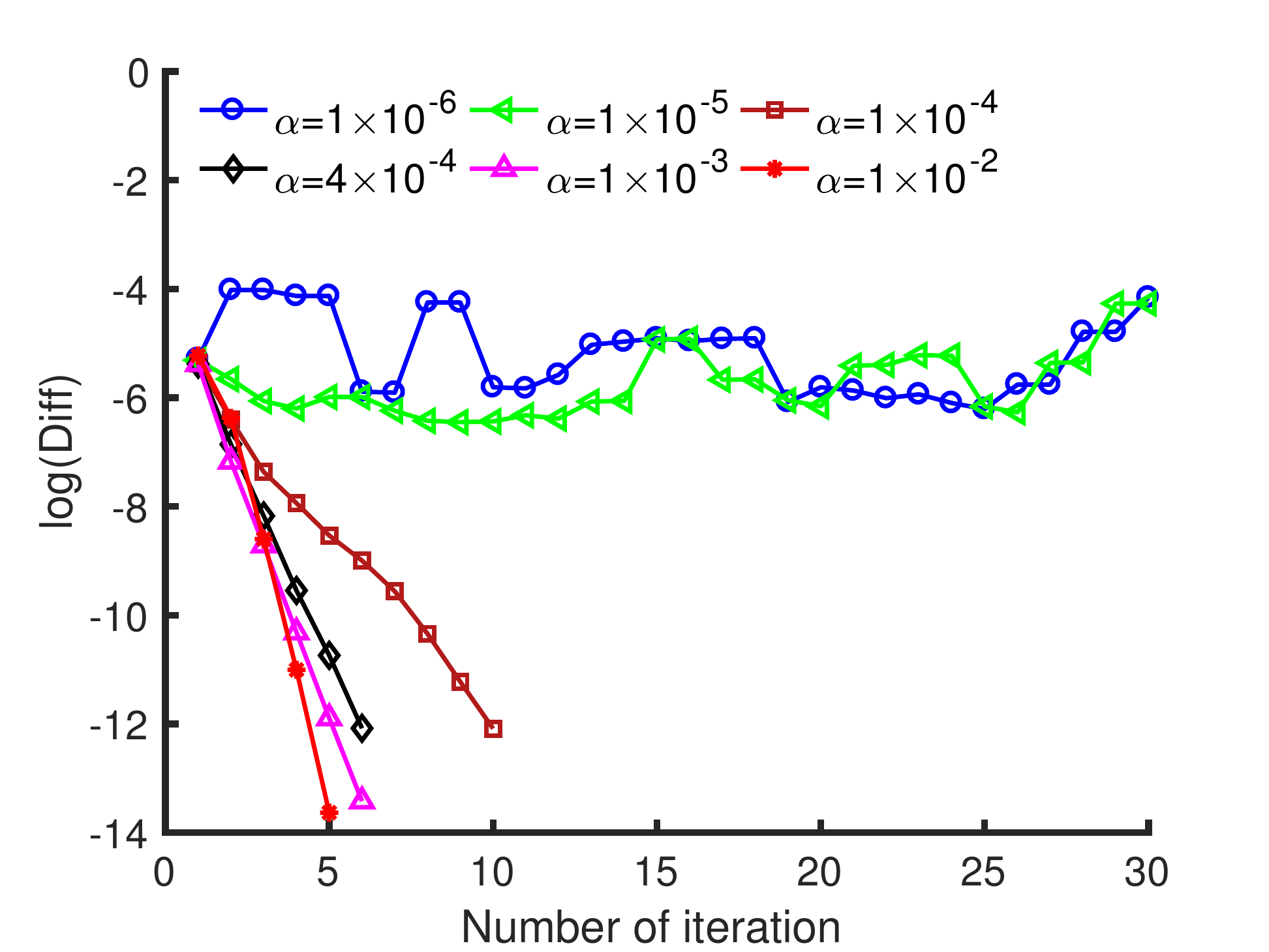}}
\vfill
\subfigure[Exponential convolution at step 1]{\label{fig:subfig:exponential_convolution}
\includegraphics[width=0.34\linewidth]{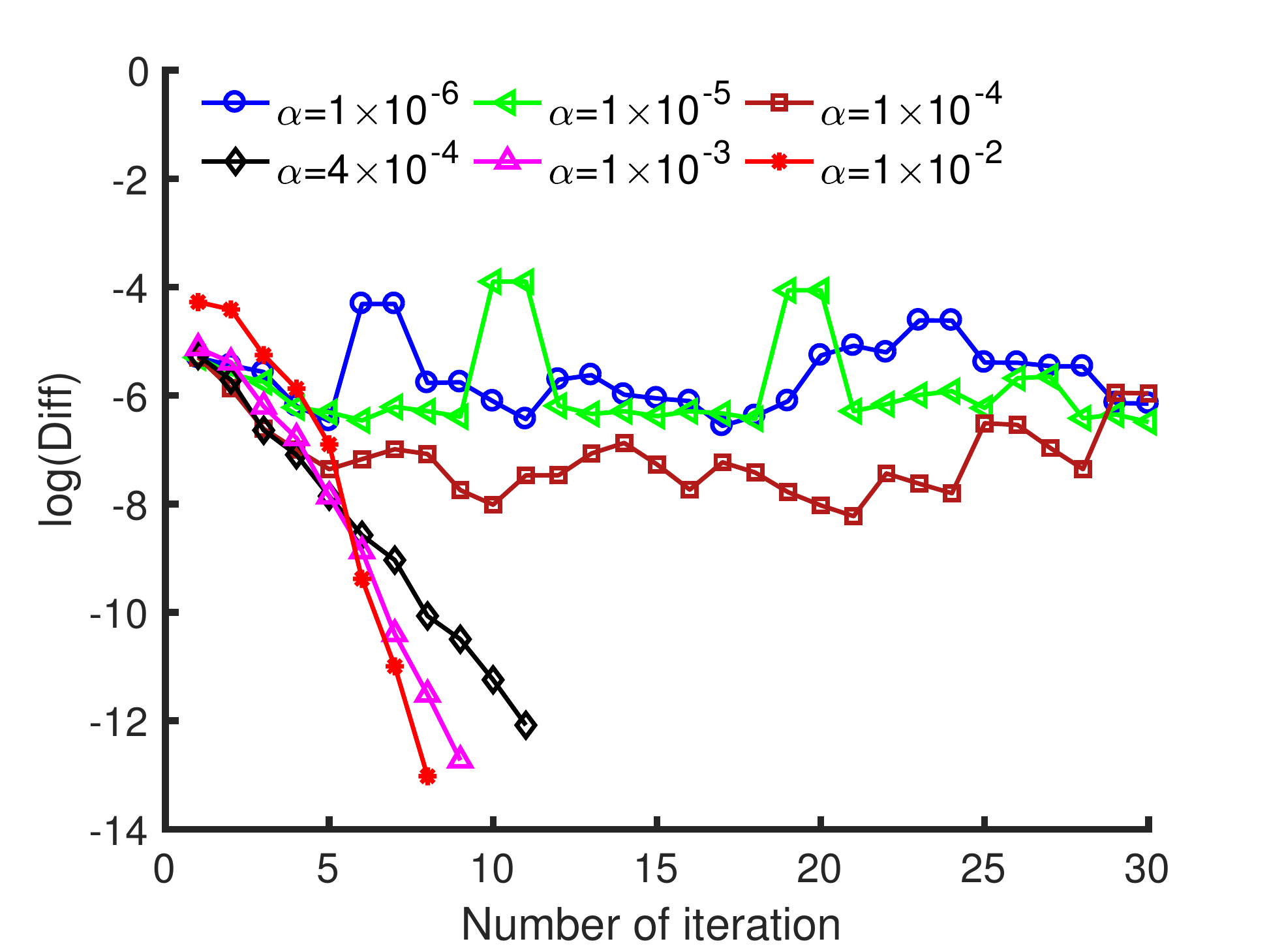}}
\subfigure[Exponential convolution at step 2]{\label{fig:subfig:exponential_convolution_step2}
\includegraphics[width=0.34\linewidth]{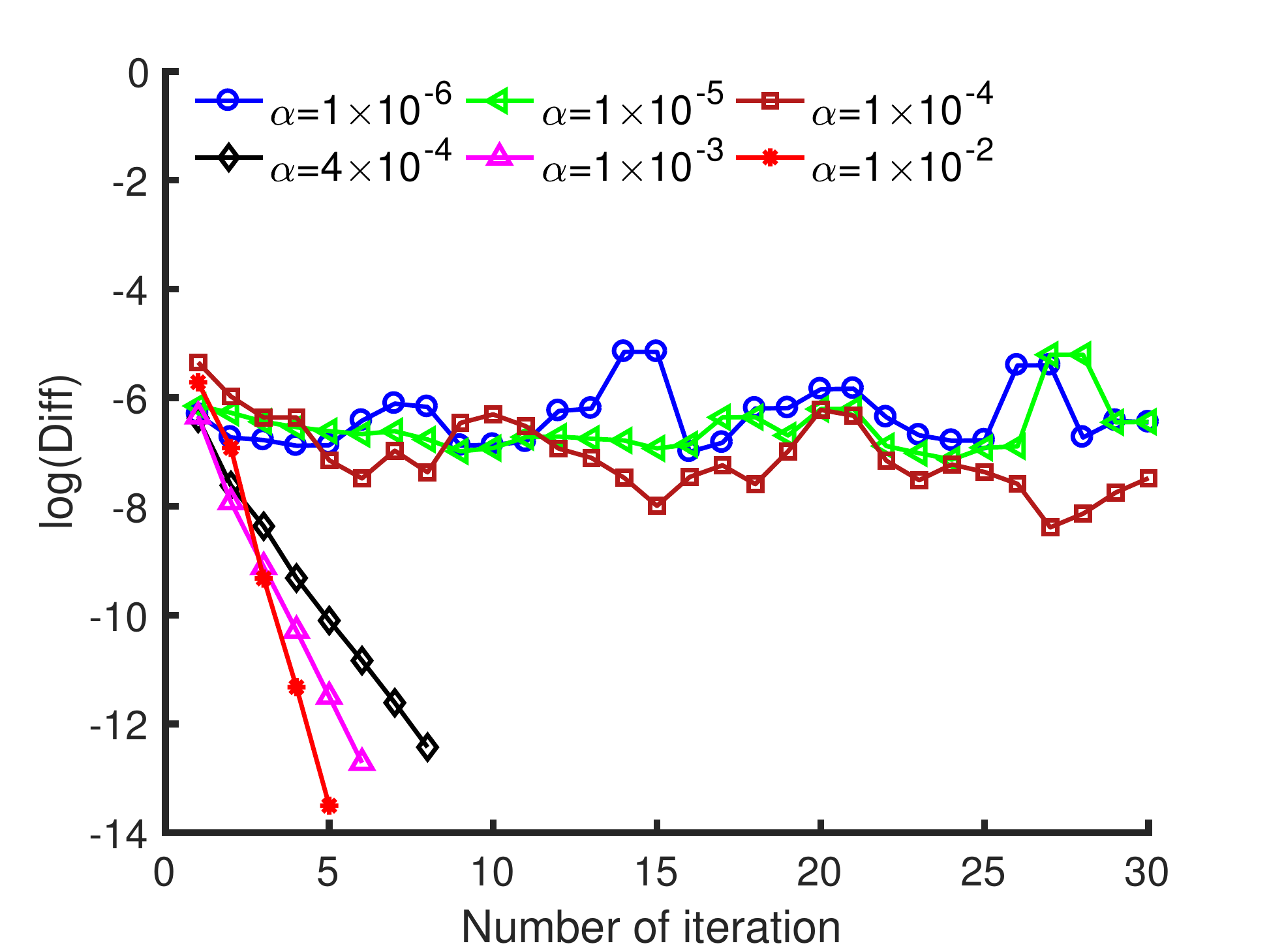}}
\vfill
\subfigure[Smoothed 2-point convolution at step 1]{\label{fig:subfig:smoothed_2_point_convolution}
\includegraphics[width=0.34\linewidth]{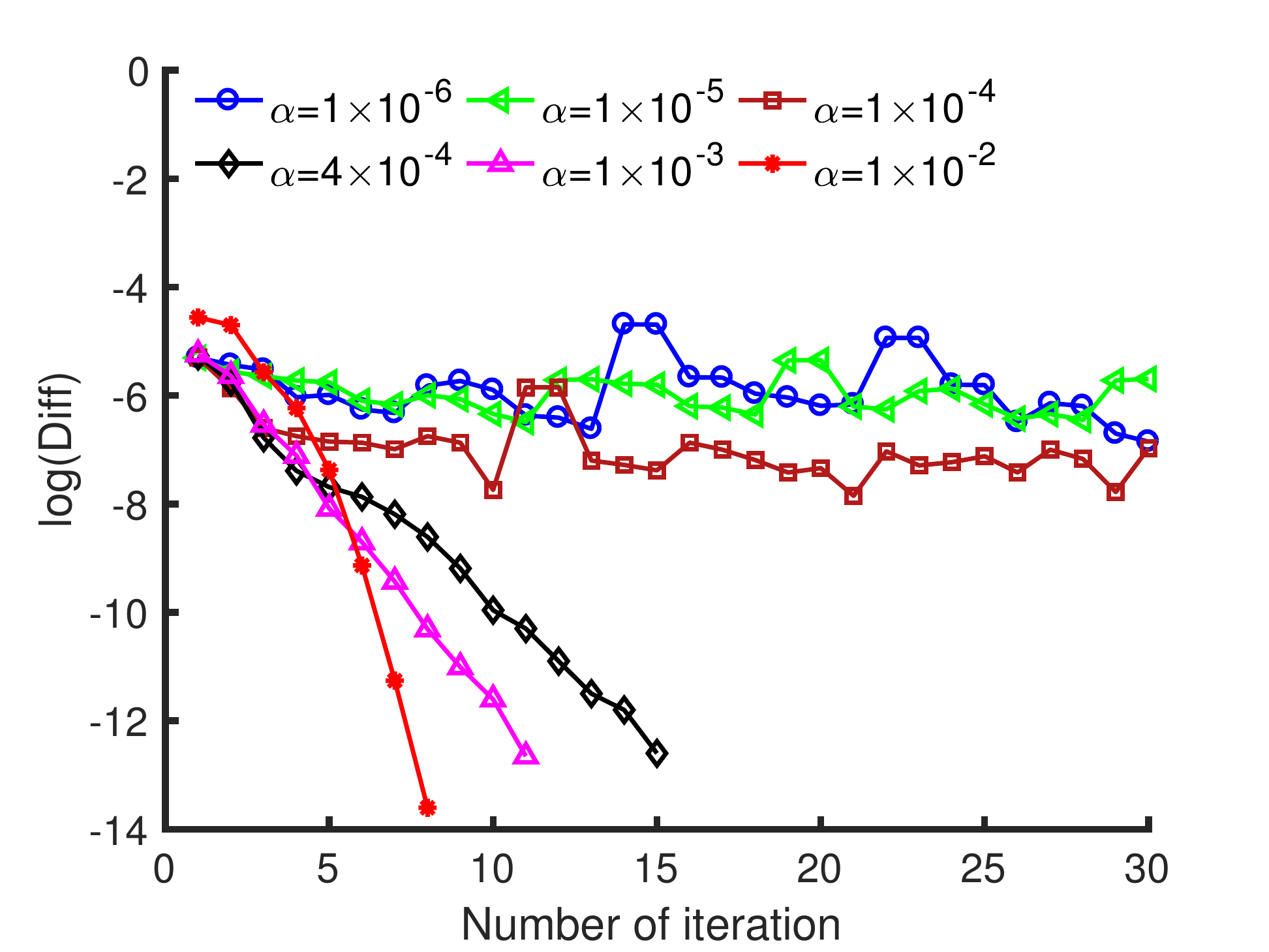}}
\subfigure[Smoothed 2-point convolution at step 2]{\label{fig:subfig:smoothed_2_point_convolution_step2}
\includegraphics[width=0.34\linewidth]{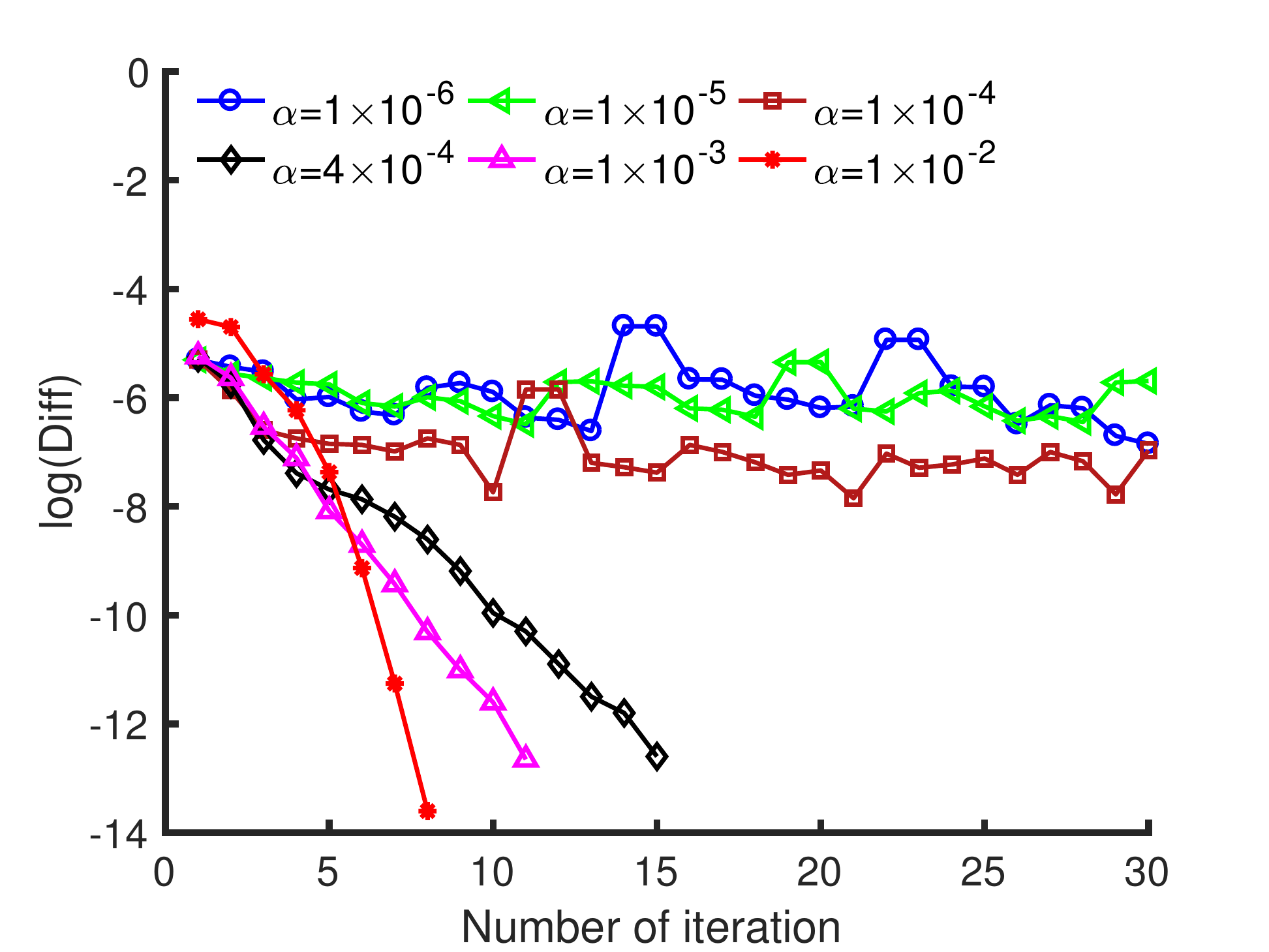}}
\caption{Example 1. Convergence of Newton's iteration using different regularization methods for the tension
test where Diff denotes the $L^2$ norm of the difference between two consecutive approximations at steps 1
and 2. ($N = 6,400$ and $l = 0.0075$~mm)
\label{fig:tension_convergence}}
\end{figure}

We now investigate the effects of the regularization methods (cf. Section \ref{SEC:convergence})
on the convergence of Newton's iteration and load-deflection curves.
For convergence, we consider the first displacement increment (where $\Delta U = 1\times10^{-5}$~mm
is applied on the top edge) with $k_l = 0$. The convergence history of Newton's iteration is shown
in Fig. \ref{fig:tension_convergence}. As can be seen, Newton's iteration fails to converge 
without regularization on the decomposition of the strain tensor.
On the other hand, convergence is reached using the sonic-point regularization with $\alpha \ge 1\times 10^{-4}$
and exponential convolution and smoothed 2-point convolution methods with $\alpha \ge 4\times 10^{-4}$.
Moreover, Newton's iteration converges faster for larger $\alpha$ for all of the three methods.

\begin{figure} [!htb]
\centering 
\subfigure[$\alpha = 10^{-3}$ for three regularization methods,
without regularization $\alpha = 0$ ($k_l =10^{-3})$]{\label{fig:subfig:Tld_Diffmeth}
\includegraphics[width=0.4\linewidth]{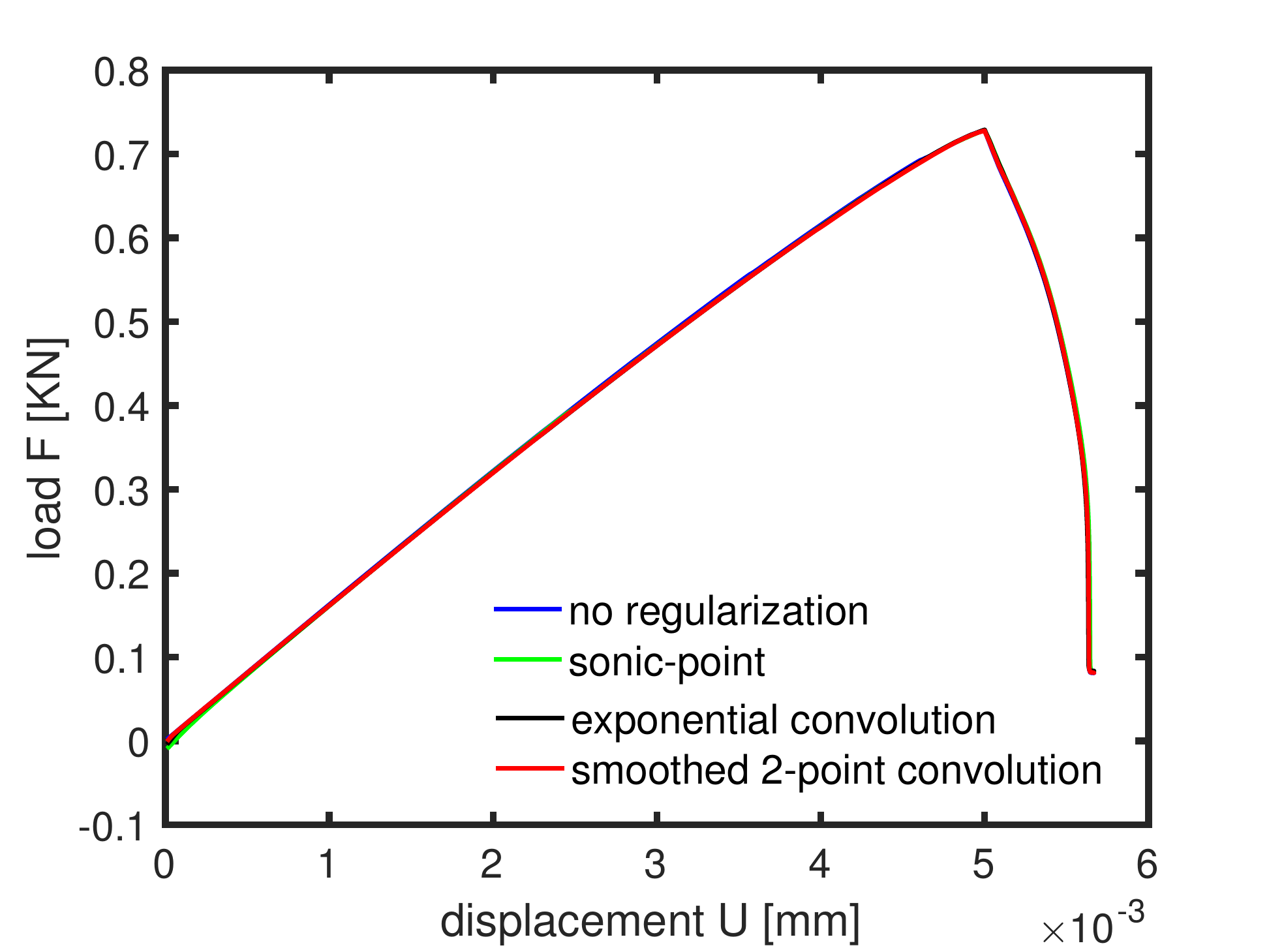}}
\subfigure[sonic-point ($k_l = 0$)]{\label{fig:subfig:Tld_job1}
\includegraphics[width=0.4\linewidth]{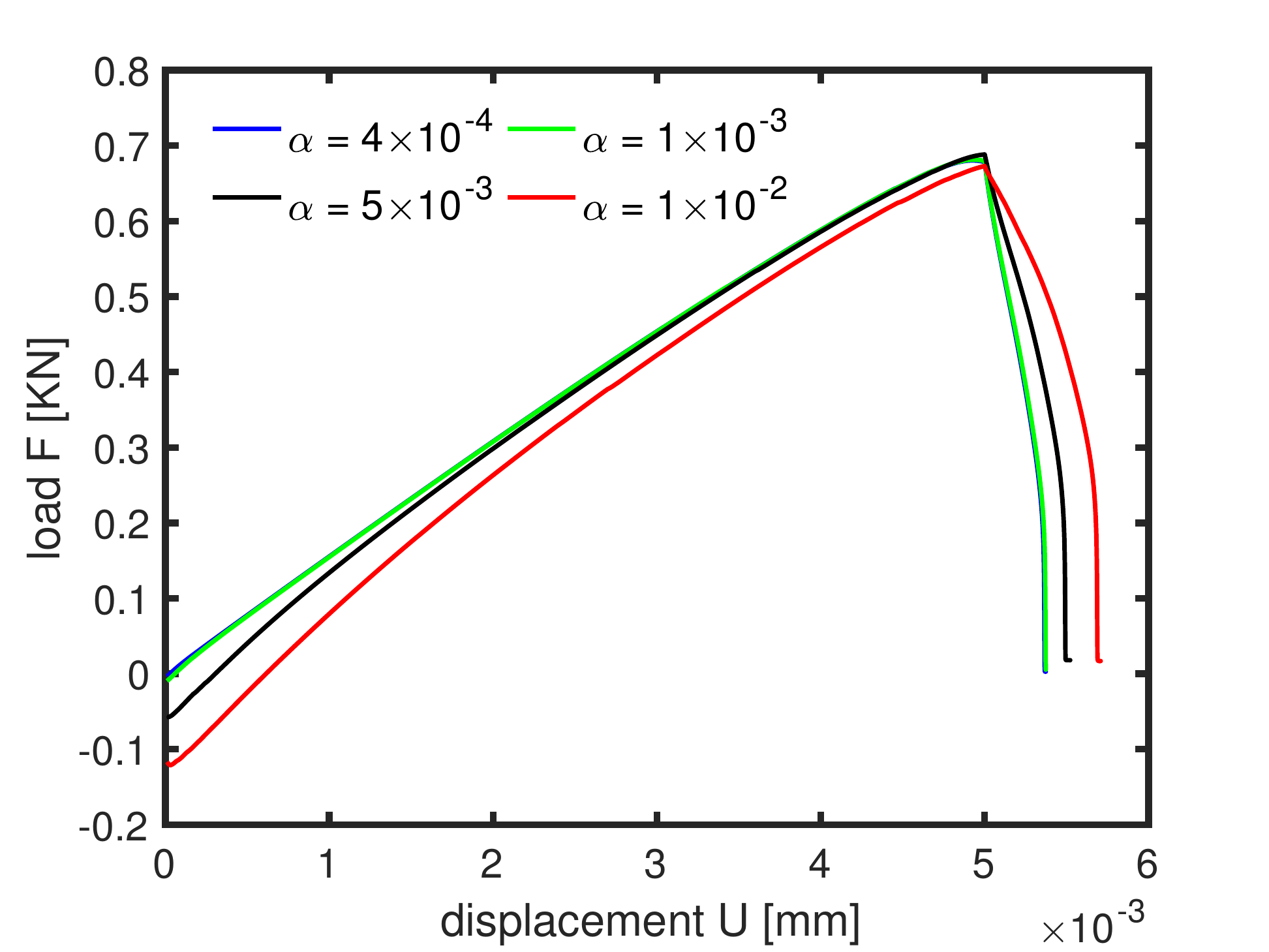}}
\subfigure[exponential convolution ($k_l = 0$)]{\label{fig:subfig:Tld_job2}
\includegraphics[width=0.4\linewidth]{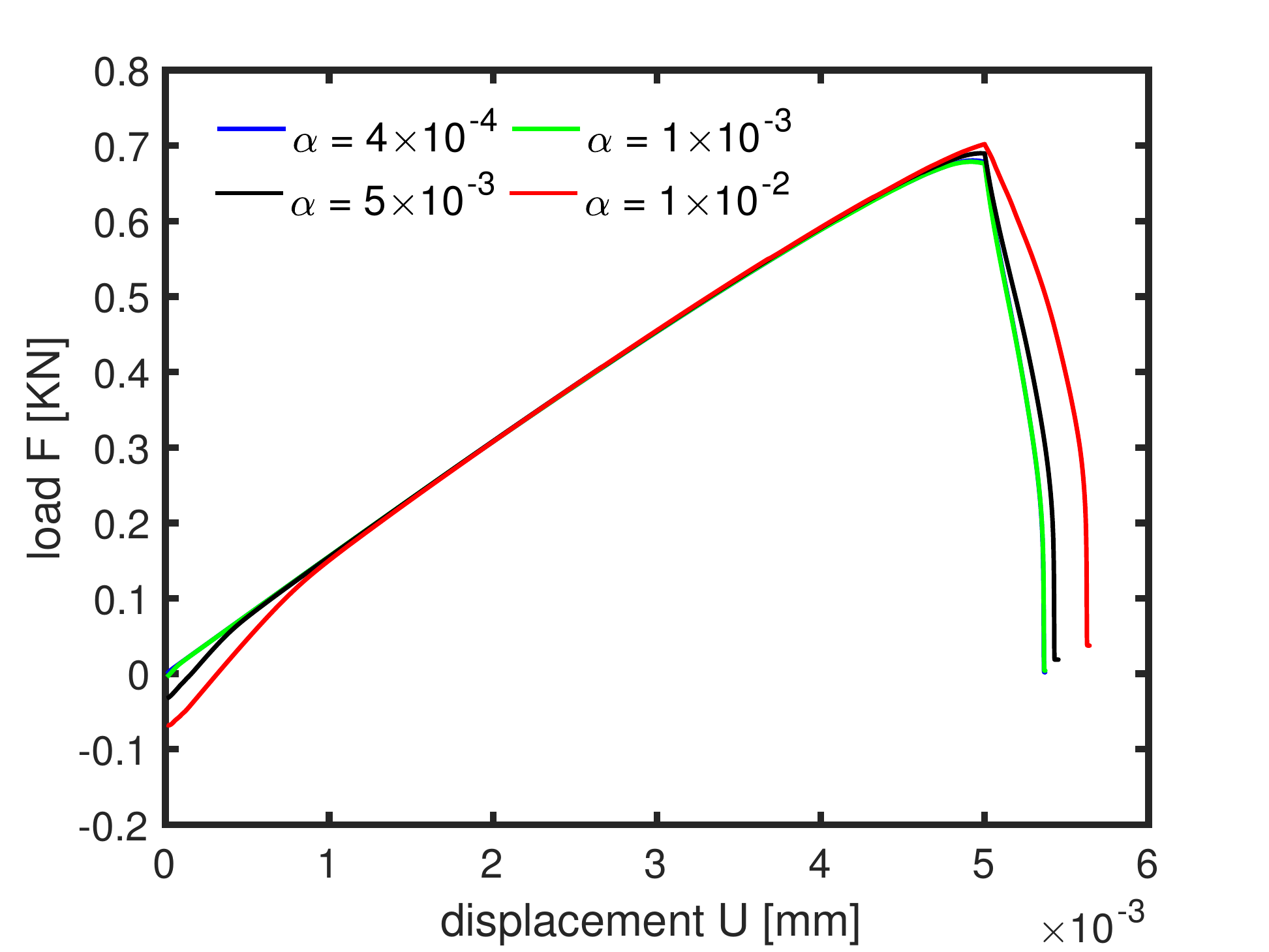}}
\subfigure[smoothed 2-point convolution ($k_l = 0$)]{\label{fig:subfig:Tld_job3}
\includegraphics[width=0.4\linewidth]{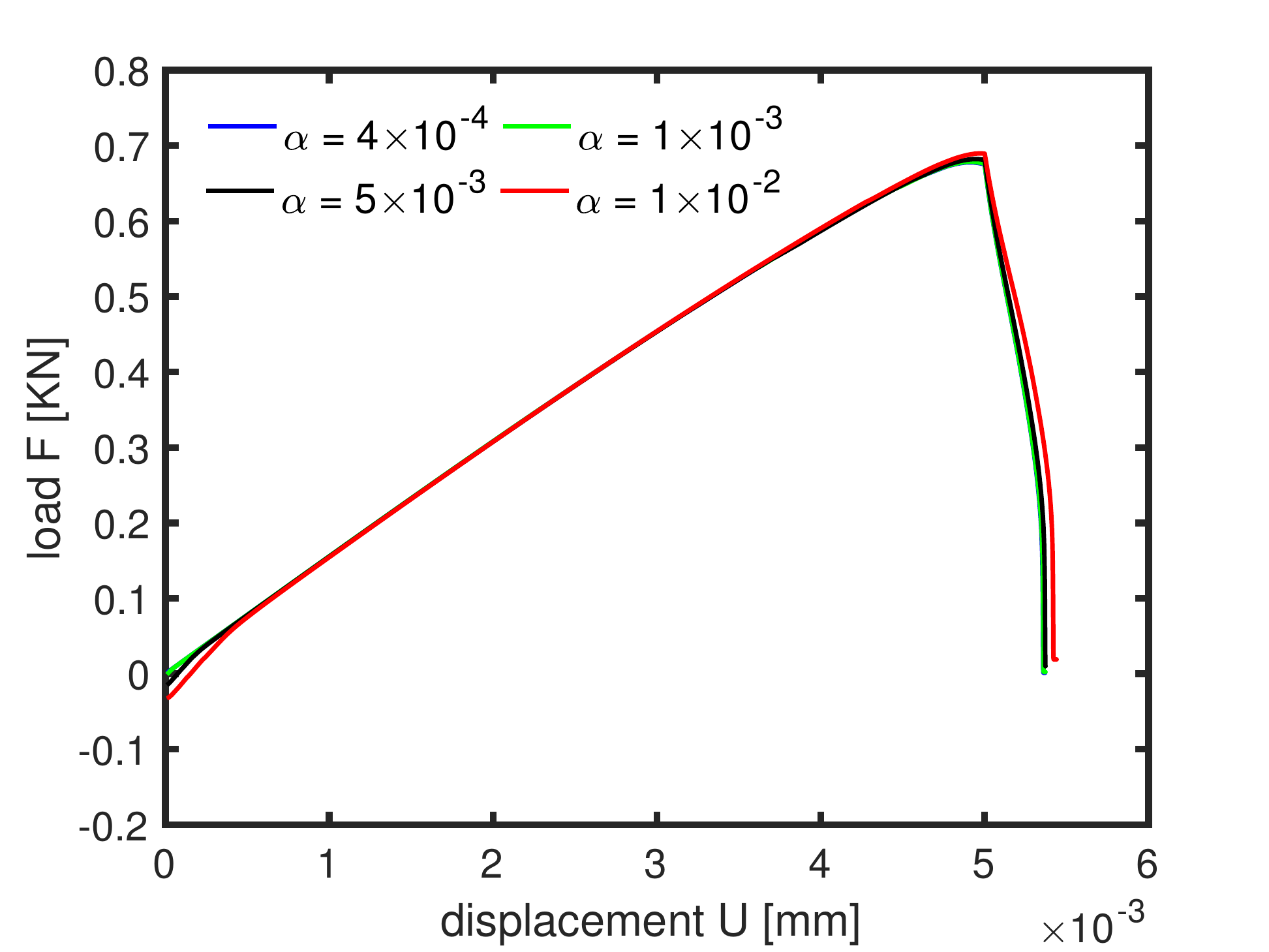}}
\caption{Example 1. The load-deflection curves for the regularization methods for
the tension test with $l=0.0075$ mm and $N=6,400$.
(a) Different regularization methods; (b) Sonic-point method with various $\alpha$;
(c) Exponential convolution method with various $\alpha$;
(d) Smoothed 2-point convolution method with various $\alpha$.}
\label{fig:tension effects of the regularization method}
\end{figure}

Next, we examine the effects of the regularization methods on load-deflection curves. We first compare the three
regularization methods ($\alpha = 1\times10^{-3}$, $l=0.0075$~mm and $k_l = 1\times10^{-3}$).
For the computation without regularization ($\alpha = 0$), $k_l = 1\times10^{-3}$ is used to obtain convergent results.
The results are shown in Fig. \ref{fig:subfig:Tld_Diffmeth}. As we can see, there is no significant difference
between the regularization methods with $\alpha = 1\times10^{-3}$ and without regularization,
except at the beginning of the load time.

We then investigate the effects of the regularization parameter. The load-deflection curves
for the three regularization methods with various values of
$\alpha$ ($4\times10^{-4}, 1\times10^{-3}, 5\times10^{-3}$, and $1\times10^{-2}$) are shown
in Fig. \ref{fig:subfig:Tld_job1}, \ref{fig:subfig:Tld_job2}, \ref{fig:subfig:Tld_job3}, respectively.
We can see that, for all three regularization methods, when $\alpha$ is small ($\alpha \le 1\times10^{-3}$),
influence on the numerical solution is almost invisible. For large $\alpha$ values ($\alpha \ge 5\times10^{-3}$),
all three regularization methods underestimate the load before the crack starts propagating and
overestimates afterwards. Moreover, the effects using the smoothed 2-point convolution method are less
significant than those with the sonic-point method. The effects of the exponential convolution method
lie between those of the smoothed 2-point convolution and sonic-point methods. 
Since the sonic-point method has the simplest form and is more economic to compute than the other two methods,
and since all three methods give almost identical results when $\alpha$ is small,
we use the sonic-point method with $\alpha = 1\times10^{-3}$ for later computations.

\subsubsection{Effects of mesh adaptivity} 

\begin{figure} [!htb]
\centering 
\subfigure[uniform mesh]{\label{fig:subfig:Tld_hrefine_uniform}
\includegraphics[width=0.4\linewidth]{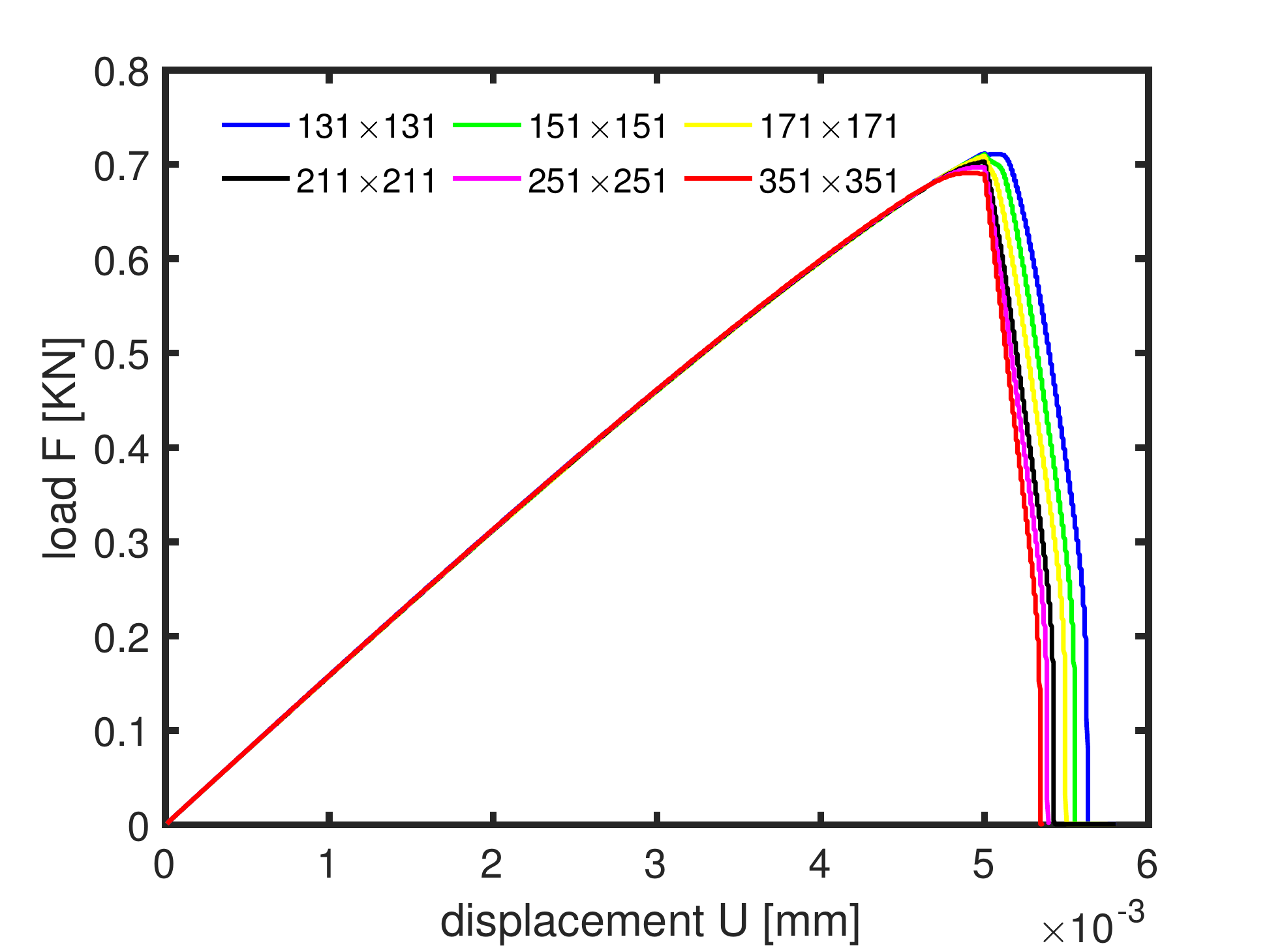}}
\subfigure[adaptive moving mesh]{\label{fig:subfig:Tld_hrefine_moving}
\includegraphics[width=0.4\linewidth]{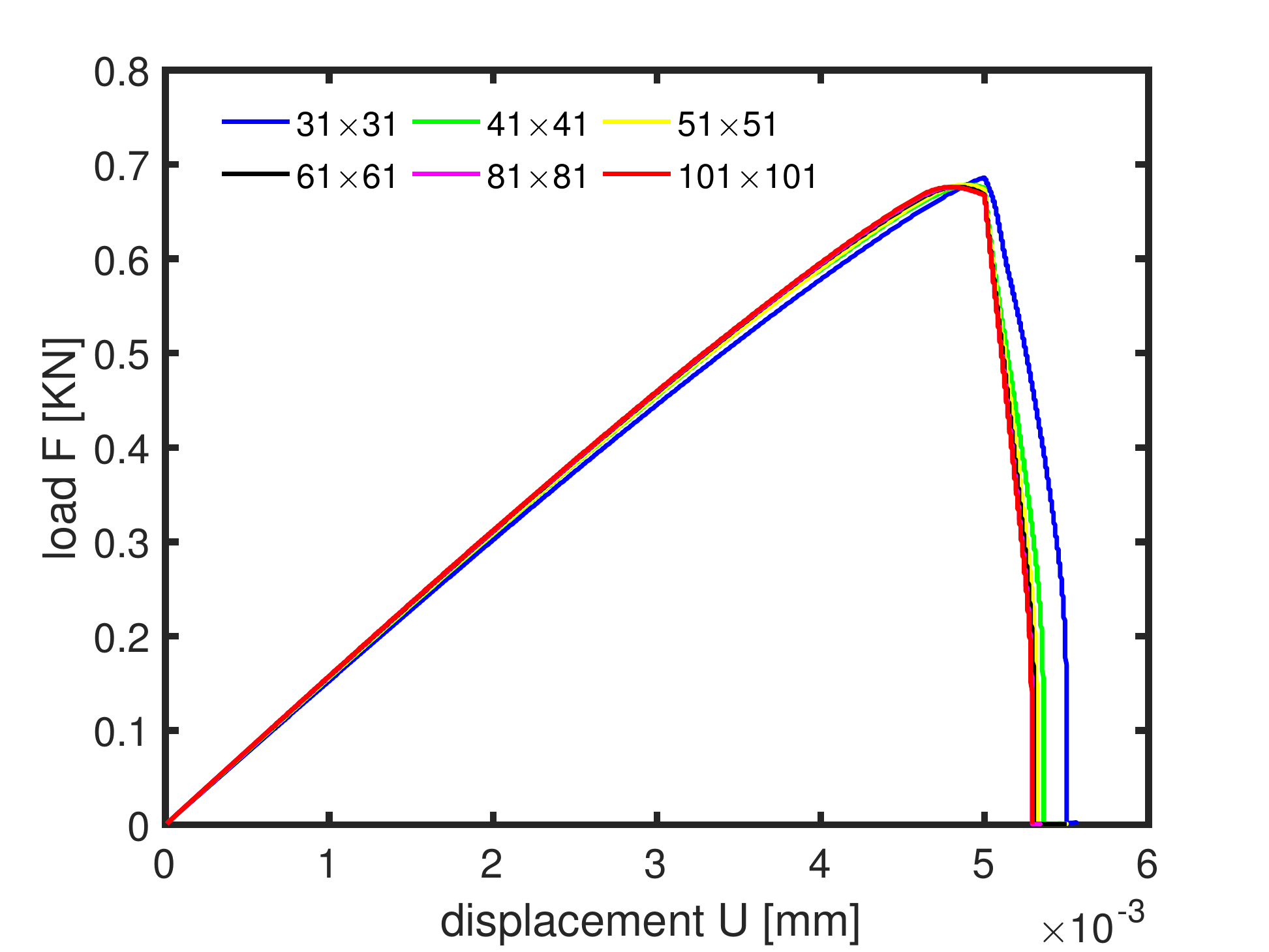}}
\caption{Example 1. The load-deflection curves for $l=0.0075$ mm for meshes of different size:
(a) uniform mesh; (b) adaptive moving mesh.}
\label{fig:mesh convergence}
\end{figure}

\begin{figure} [!htb]
\centering 
\includegraphics[width=0.4\linewidth]{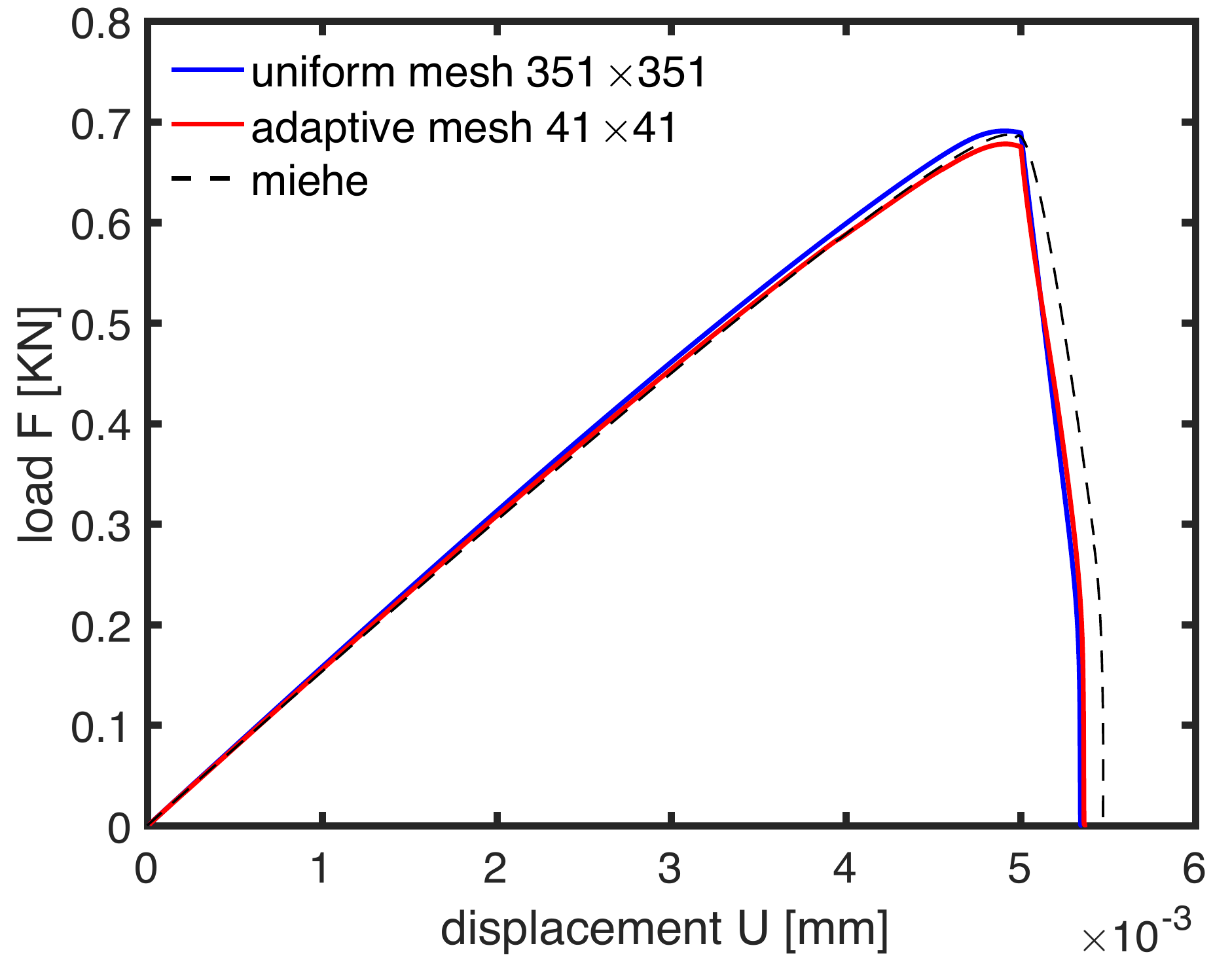}
\caption{Example 1. The load-deflection curves for different meshes for the tension test with $l=0.0075$~mm.}
\label{fig:mesh-compare}
\end{figure}

We now study the effects of mesh adaptivity on the numerical solution. We take $l=0.0075$~mm in the computation
and compare the results with uniform meshes. As mentioned before, the mesh should be sufficiently fine
($h \ll l$) to resolve the cracks. A uniform triangular mesh is obtained
by first partitioning the domain into a uniform rectangular mesh and then dividing each sub-rectangle
into four triangles using the diagonal lines and is represented by the number of sub-rectangles
in $x-$ and $y$-directions. We consider different mesh sizes as $131 \times 131$ ($N = 67,600$),
$151\times151$ ($N = 90,000$), $171\times171$ ($N = 115,600$), $211\times211$ ($N = 176,400$),
$251\times251$ ($N = 250,000$), and $351\times351$ ($N = 490,000$). For adaptive meshes, we start
with uniform meshes of size $31\times31$ ($N = 3,600$), $41\times41$ ($N = 6,400$), and $51\times51$ ($N = 10,000$). 
Fig. \ref{fig:mesh convergence} shows the load-deflection curves for different meshes. It clearly demonstrates
the spatial convergence in terms of the number of mesh elements for both uniform and adaptive meshes.
Moreover, moving meshes use significantly fewer elements to achieve the same accuracy,
as can be seen from Fig. \ref{fig:subfig:Tld_hrefine_moving} and Fig. \ref{fig:mesh-compare}.

To verify the accuracy of the numerical solutions, the load-deflection curves with a uniform mesh of
$351\times351$ ($N = 490,000$) and an adaptive mesh of $41\times41$ ($N = 6,400$) are plotted
in Fig. \ref{fig:mesh-compare}. They are comparable with each other and agree very well with the result
obtained by Miehe~et~al.~\cite{MHW10} where a mesh pre-refined around the regions of the crack
and its expected propagation path is used. The computation (with a Matlab\textsuperscript \textregistered\,
implementation of the method)
takes about 1,339.7 seconds of CPU time (for one time step)
on the machine with a single AMD Opteron 6386 SE $2.80$~GHz
processor for the uniform mesh and about 141.07 seconds of CPU time for the adaptive mesh.
It is clear that the moving mesh method improves computational efficiency significantly.
The relative cost between moving mesh and solving $u$ and $d$ can be clearly seen in Table \ref{tab:CPU}.
The CPU time for mesh movement takes a large proportion ($73.58\%$) when the moving mesh method
described in the previous section was used. This can further be improved by employing
a two-level mesh movement strategy; e.g., see Huang \cite{Hua01b}.

\begin{table}[!hbp]
\centering
\caption{Average CPU time (in seconds) for one time step for solving for $d$, $u$ and mesh.
\label{tab:CPU}}
\vspace{3mm}
\begin{tabular}{|c|c|c|c|c|c|c|c|}
\hline
 \multirow{2}{*}{Mesh size} &
 \multicolumn{2}{c|}{CPU for $d$} &
 \multicolumn{2}{c|}{CPU for $u$} &
 \multicolumn{2}{c|}{CPU for mesh} &
 \multirow{2}{*}{Total CPU time} \\
 \cline{2-7}
 & time & \% & time & \% & time & \% & \multicolumn{1}{c|}{} \\
 \hline
 Moving mesh $41\times41$ & 5.39 & 3.82 & 31.89 & 22.61 & 103.79 & 73.58 & 141.07 \\
 \hline
 Uniform mesh $351\times351$ & 47.9 & 3.58 & 1291.8 & 96.4 & \multicolumn{2}{c|}{} & 1339.7 \\
 \hline
\end{tabular}
\end{table}

\subsection{Example 2. Single edge notched shear test}

In this example, we consider a single edge notched shear test, with the domain and boundary conditions shown
in Fig. \ref{fig:subfig:shear}. The bottom edge of the domain is fixed and the top edge is fixed along $y$-direction while a uniform $x$-displacement $U$ is increased with time to drive the crack propagation. The elastic bulk modulus
is $\lambda = 121.15$~kN/mm$^{2}$, the shear modulus is $\mu = 80.77$~kN/mm$^{2}$, and the fracture toughness
is $g_c = 2.7 \times 10^{-3} $~kN/mm. The displacement increment is chosen as $\Delta U = 1\times10^{-5}$~mm
for the computation. 

\begin{figure} 
\centering 
\subfigure[$U = 1.0 \times 10^{-2}$ mm]{\label{fig:subfig:sm1.0}
\includegraphics[width=0.25\linewidth]{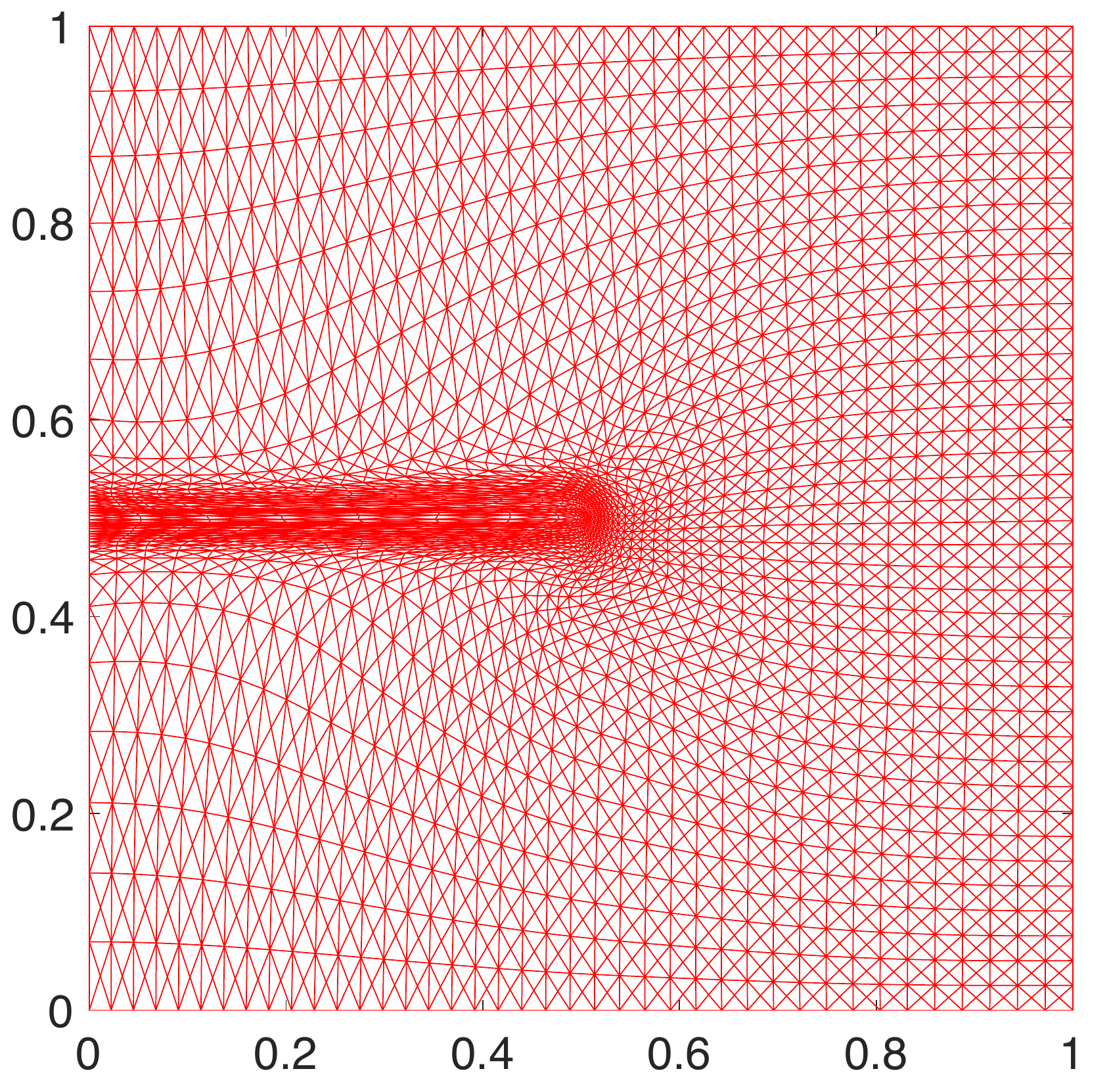}}
\subfigure[$U = 1.3 \times 10^{-2}$ mm]{\label{fig:subfig:sm1.3}
\includegraphics[width=0.25\linewidth]{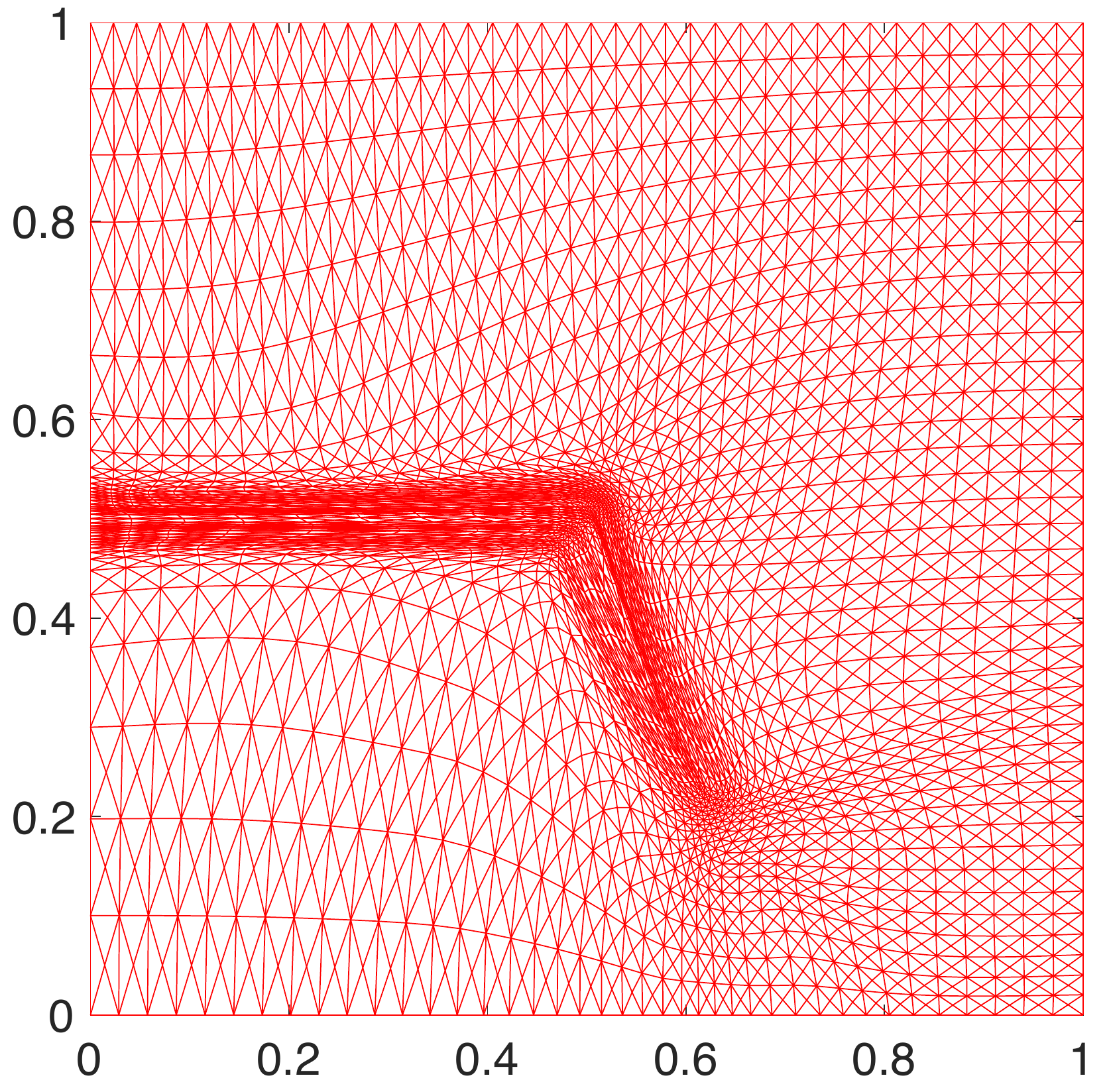}}
\subfigure[$U = 1.45 \times 10^{-2}$ mm]{\label{fig:subfig:sm1.45}
\includegraphics[width=0.25\linewidth]{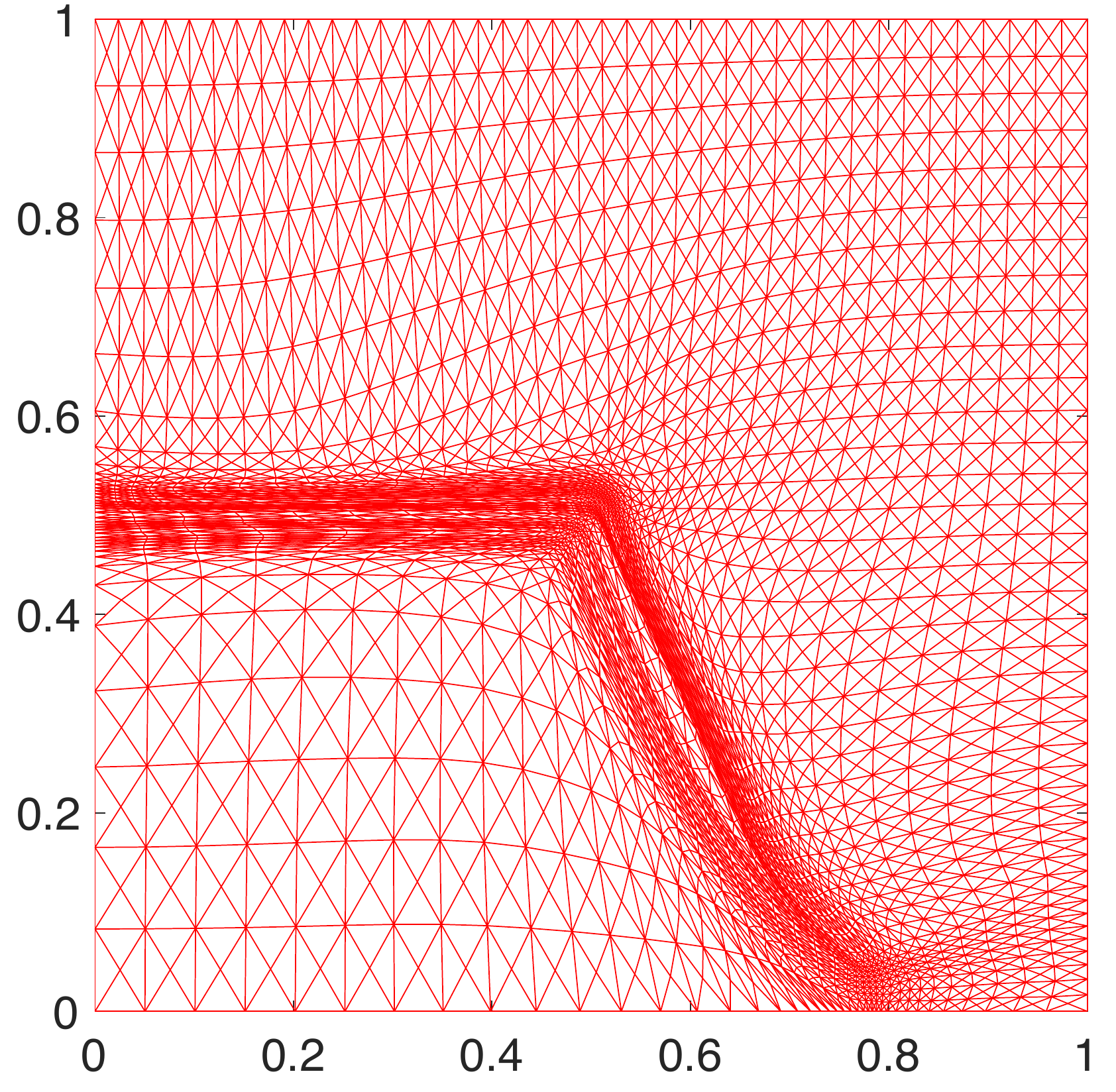}}
\vfill
\subfigure[$U = 1.0 \times 10^{-2}$ mm]{\label{fig:subfig:sd1.0}
\includegraphics[width=0.25\linewidth]{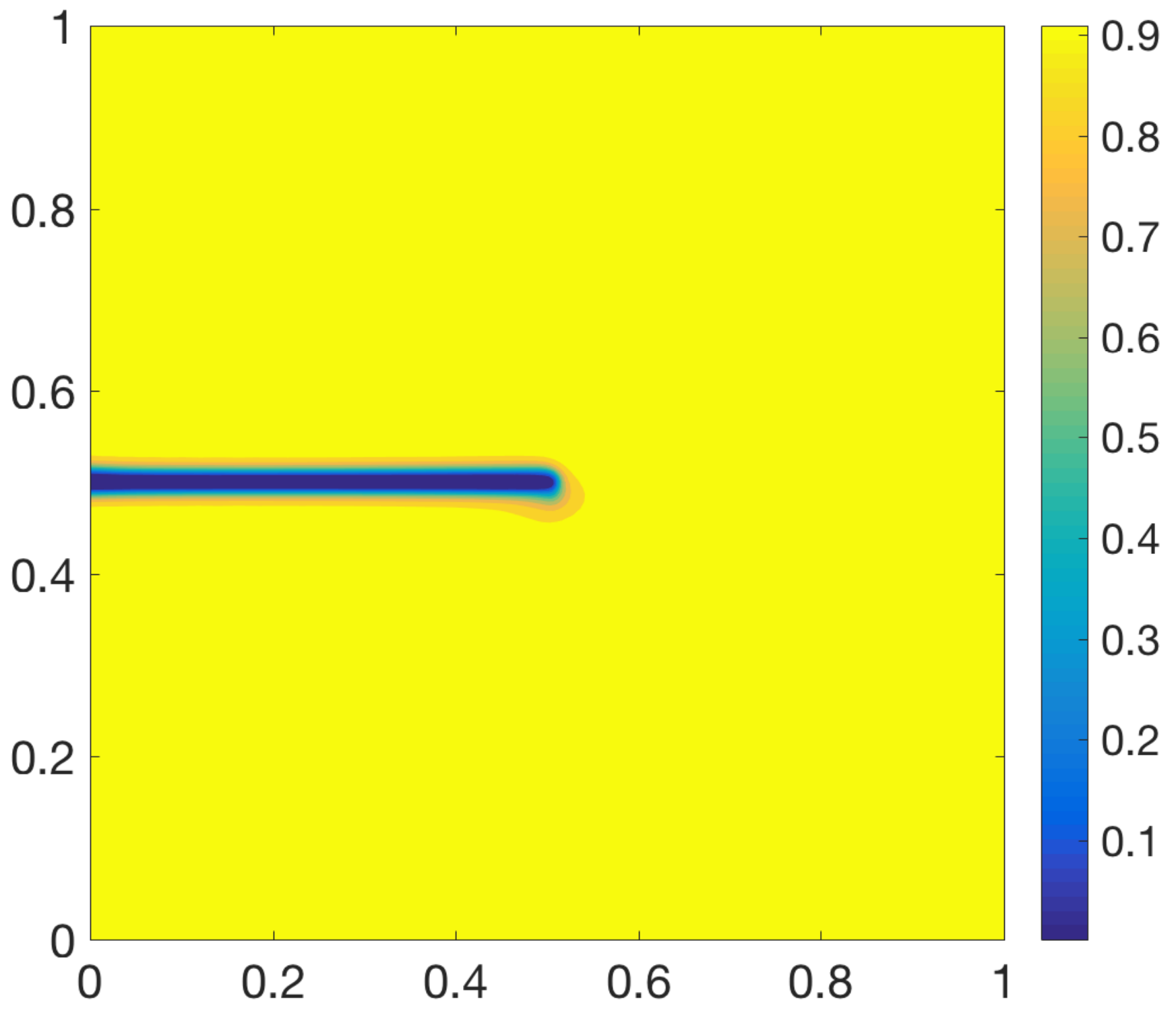}}
\subfigure[$U = 1.3 \times 10^{-2}$ mm]{\label{fig:subfig:sd1.3}
\includegraphics[width=0.25\linewidth]{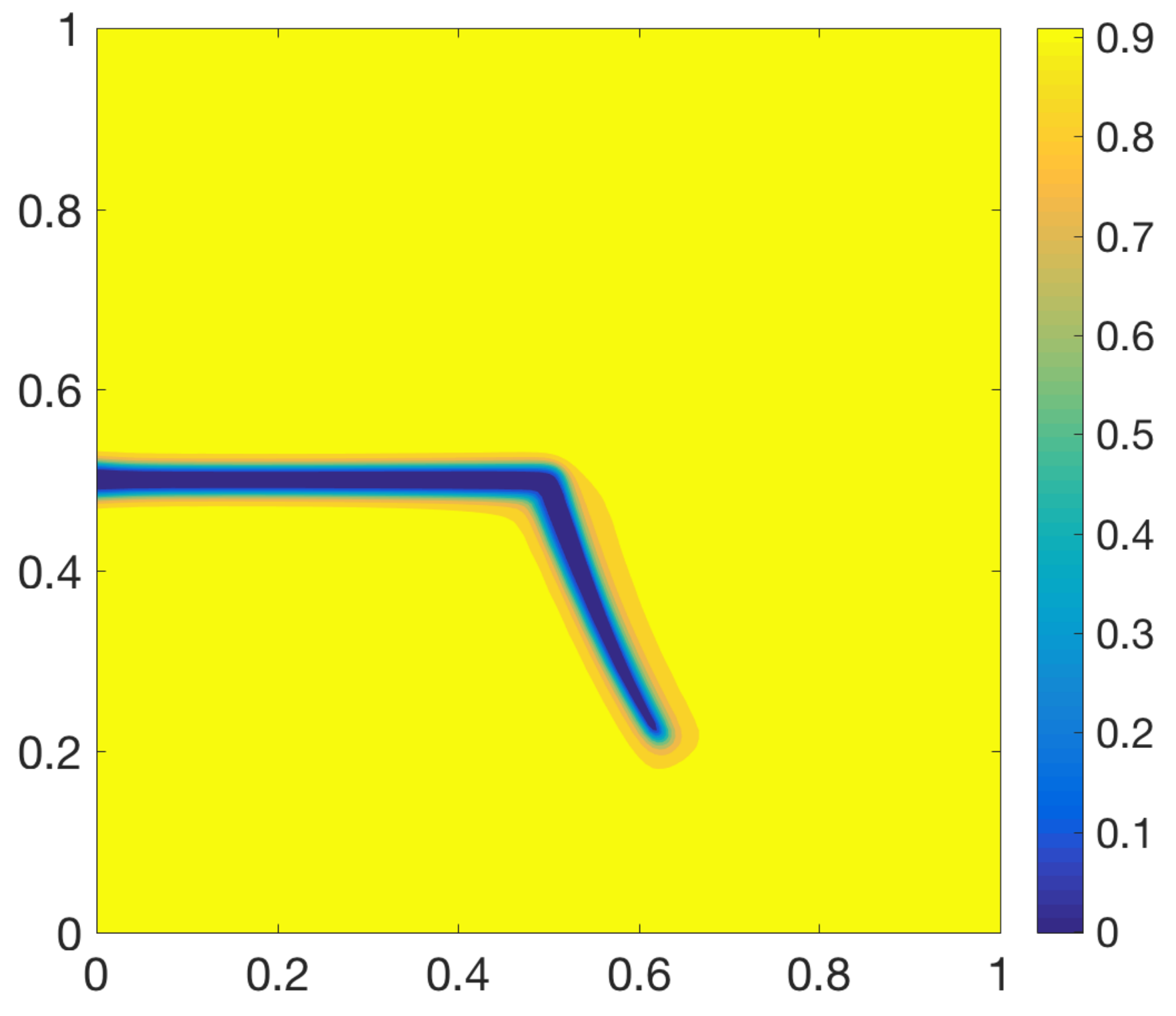}}
\subfigure[$U = 1.45 \times 10^{-2}$ mm]{\label{fig:subfig:sd1.45}
\includegraphics[width=0.25\linewidth]{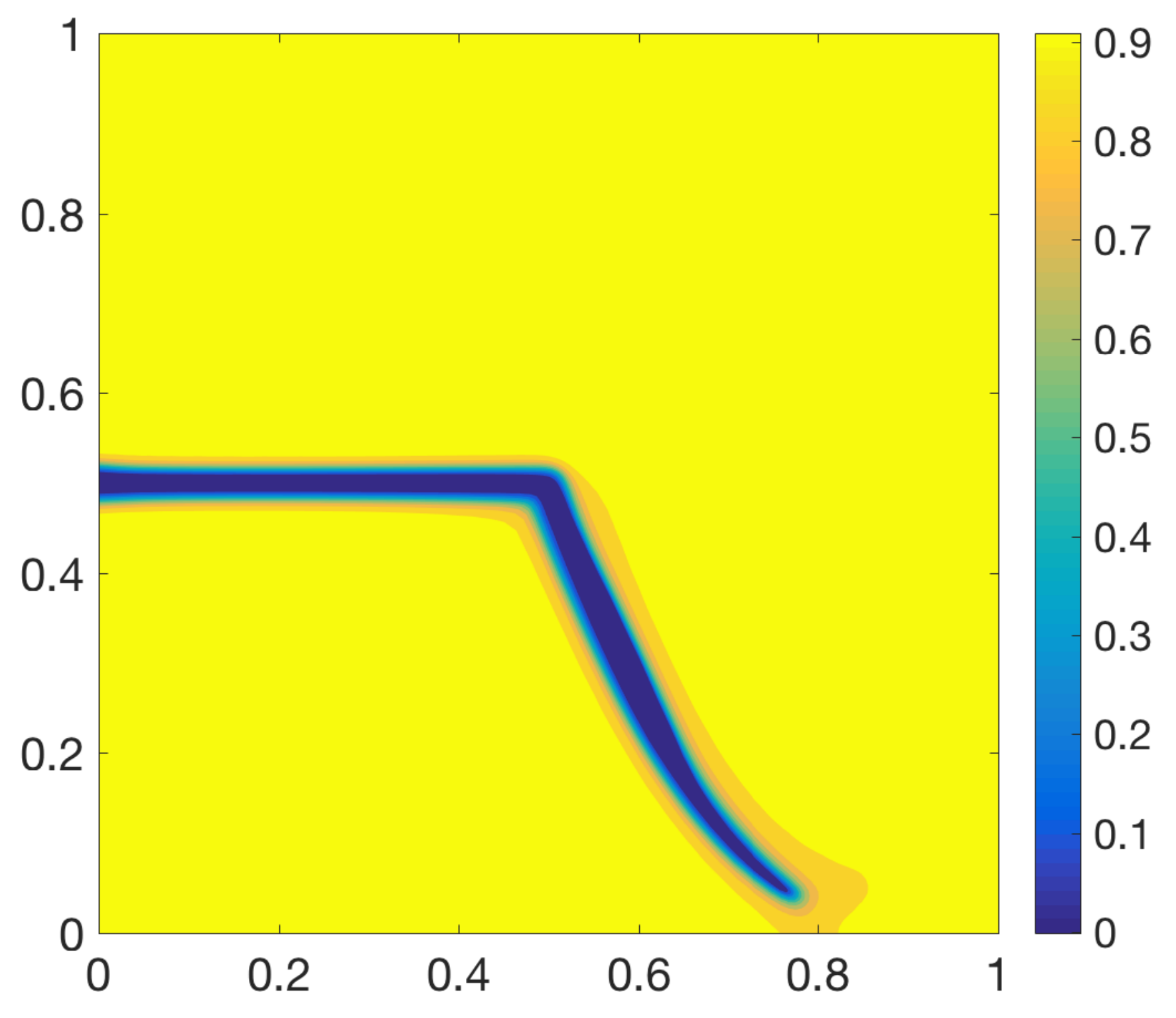}}
\caption{Example 2. The mesh and contours of the phase-field variable during crack evolution for the shear test
with $l = 0.00375$~mm and $N = 6,400$.}
\label{fig:l75_shear}
\end{figure}

\begin{figure} 
\centering 
\subfigure[$U = 1.0 \times 10^{-2}$ mm]{\label{fig:subfig:sm1.0_b}
\includegraphics[width=0.25\linewidth]{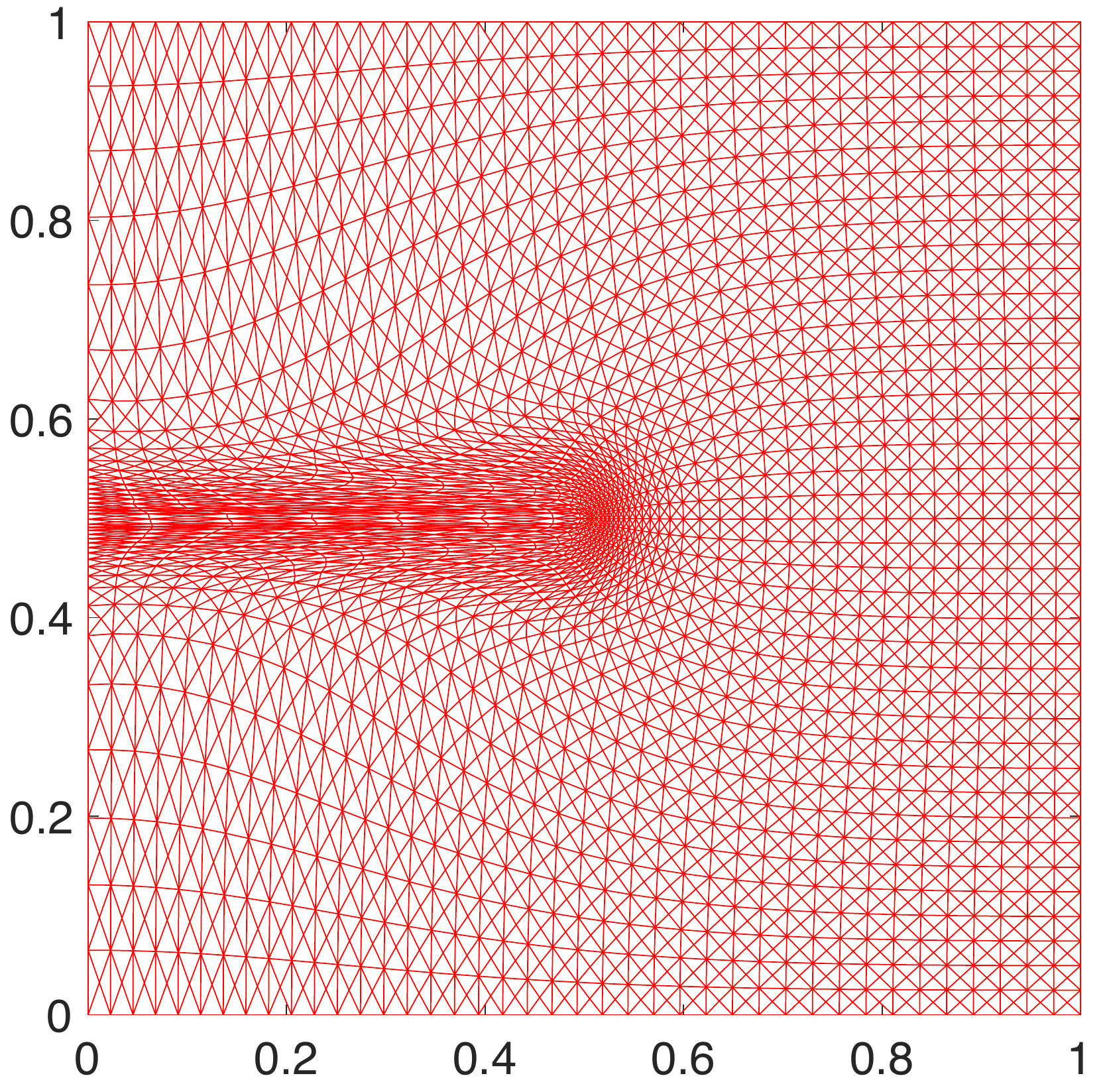}}
\subfigure[$U = 1.3 \times 10^{-2}$ mm]{\label{fig:subfig:sm1.3_b}
\includegraphics[width=0.25\linewidth]{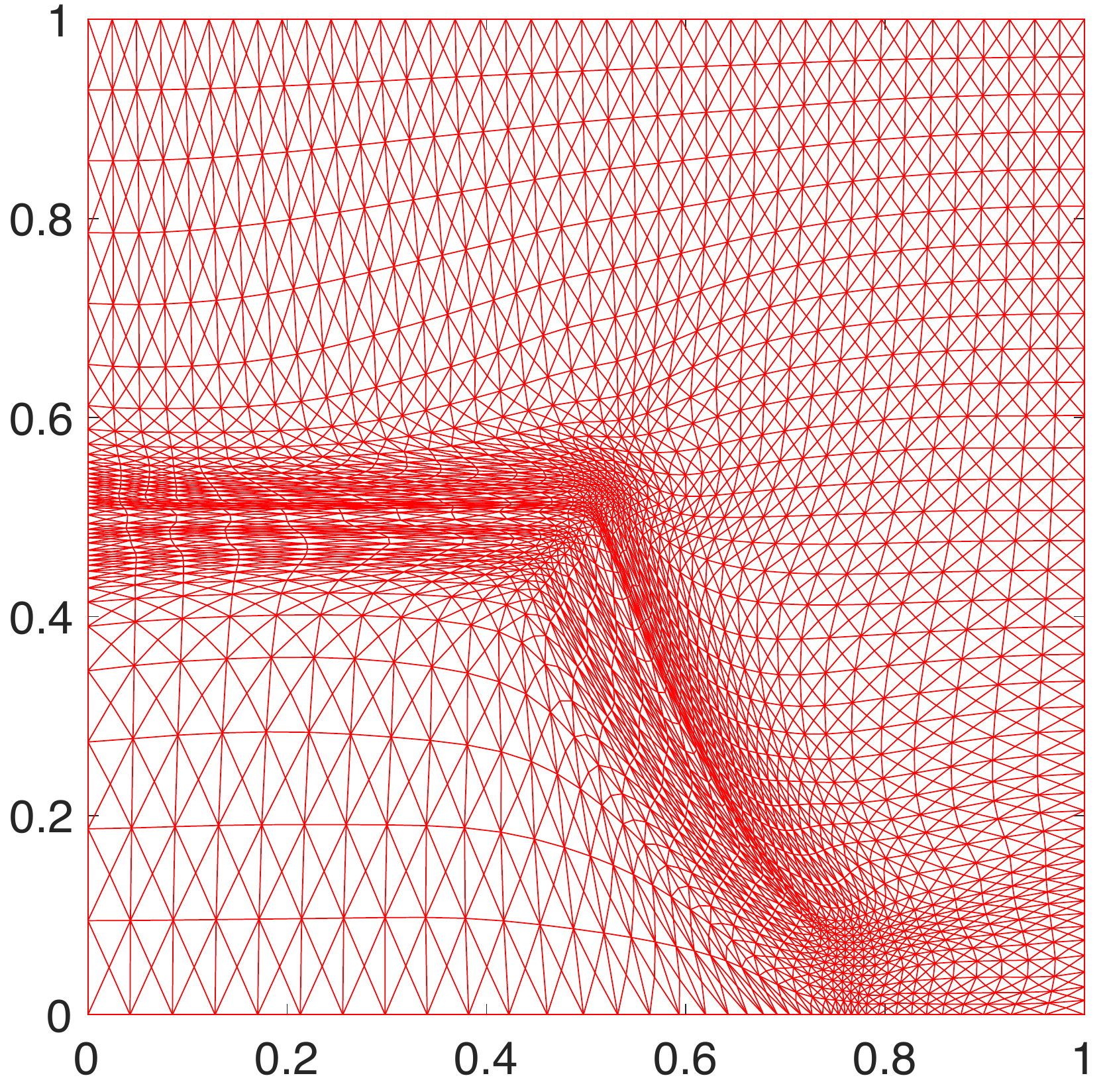}}
\subfigure[$U = 1.45 \times 10^{-2}$ mm]{\label{fig:subfig:sm1.45_b}
\includegraphics[width=0.25\linewidth]{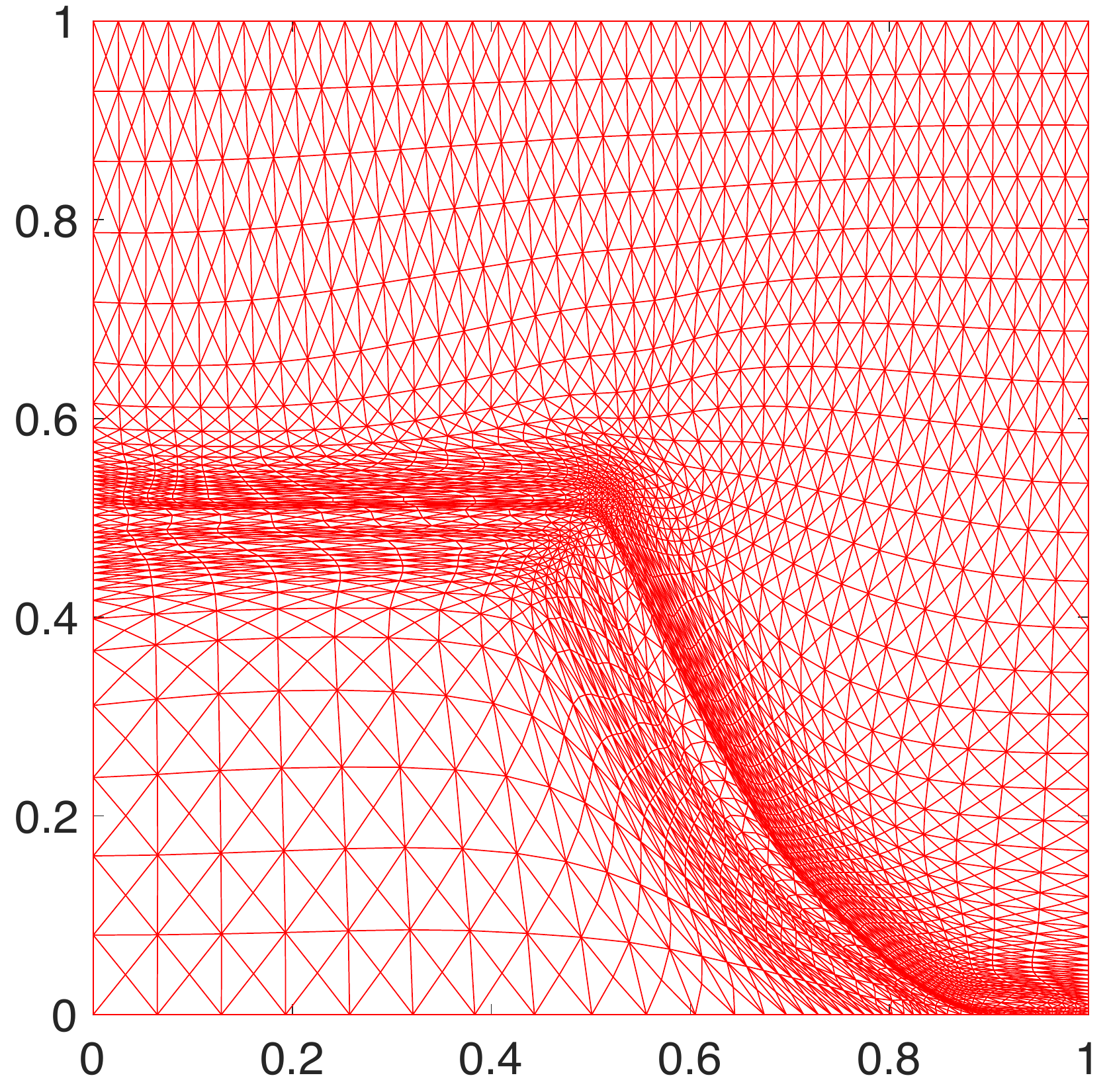}}
\vfill
\subfigure[$U = 1.0 \times 10^{-2}$ mm]{\label{fig:subfig:sd1.0_b}
\includegraphics[width=0.25\linewidth]{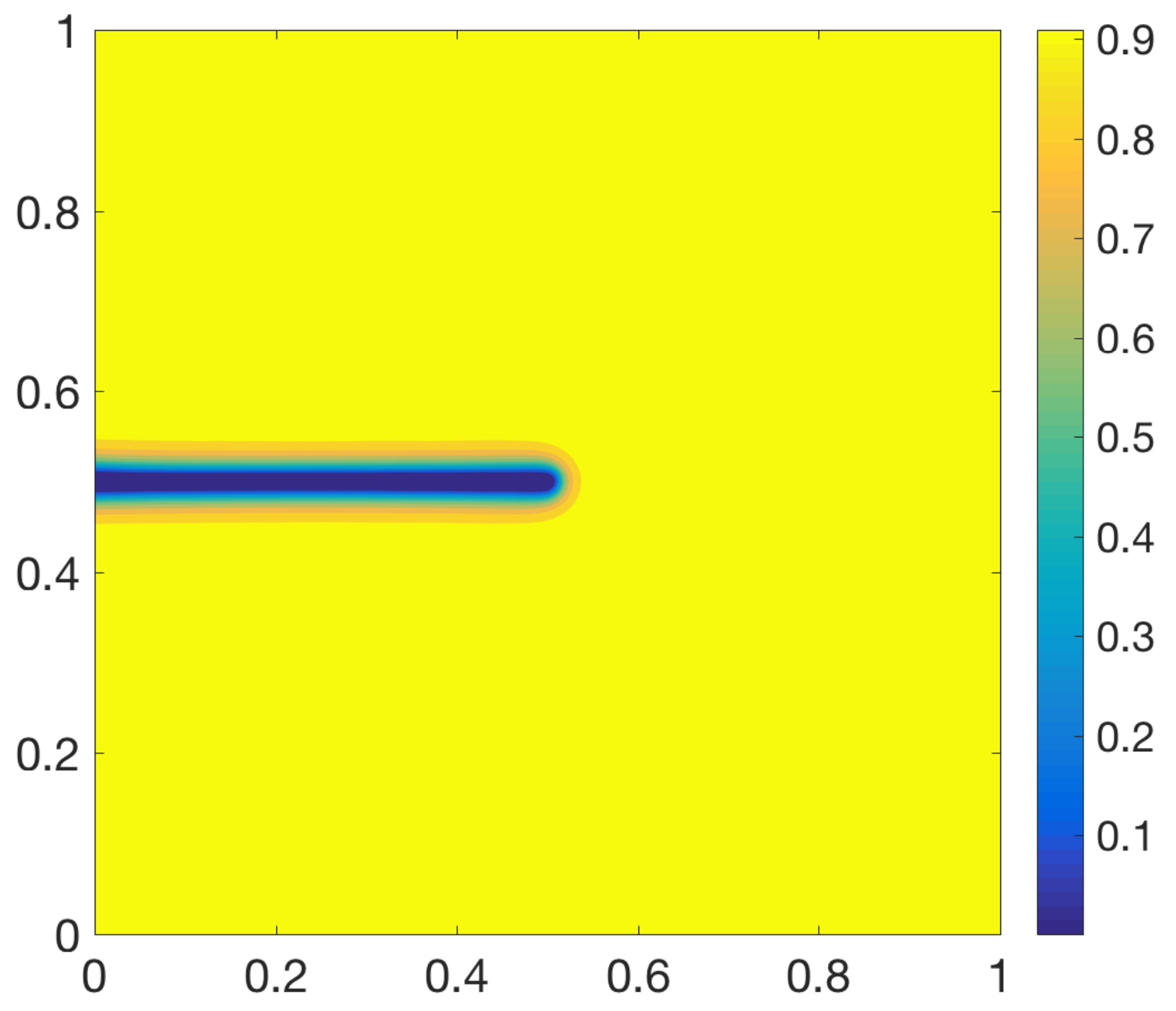}}
\subfigure[$U = 1.3 \times 10^{-2}$ mm]{\label{fig:subfig:sd1.3_b}
\includegraphics[width=0.25\linewidth]{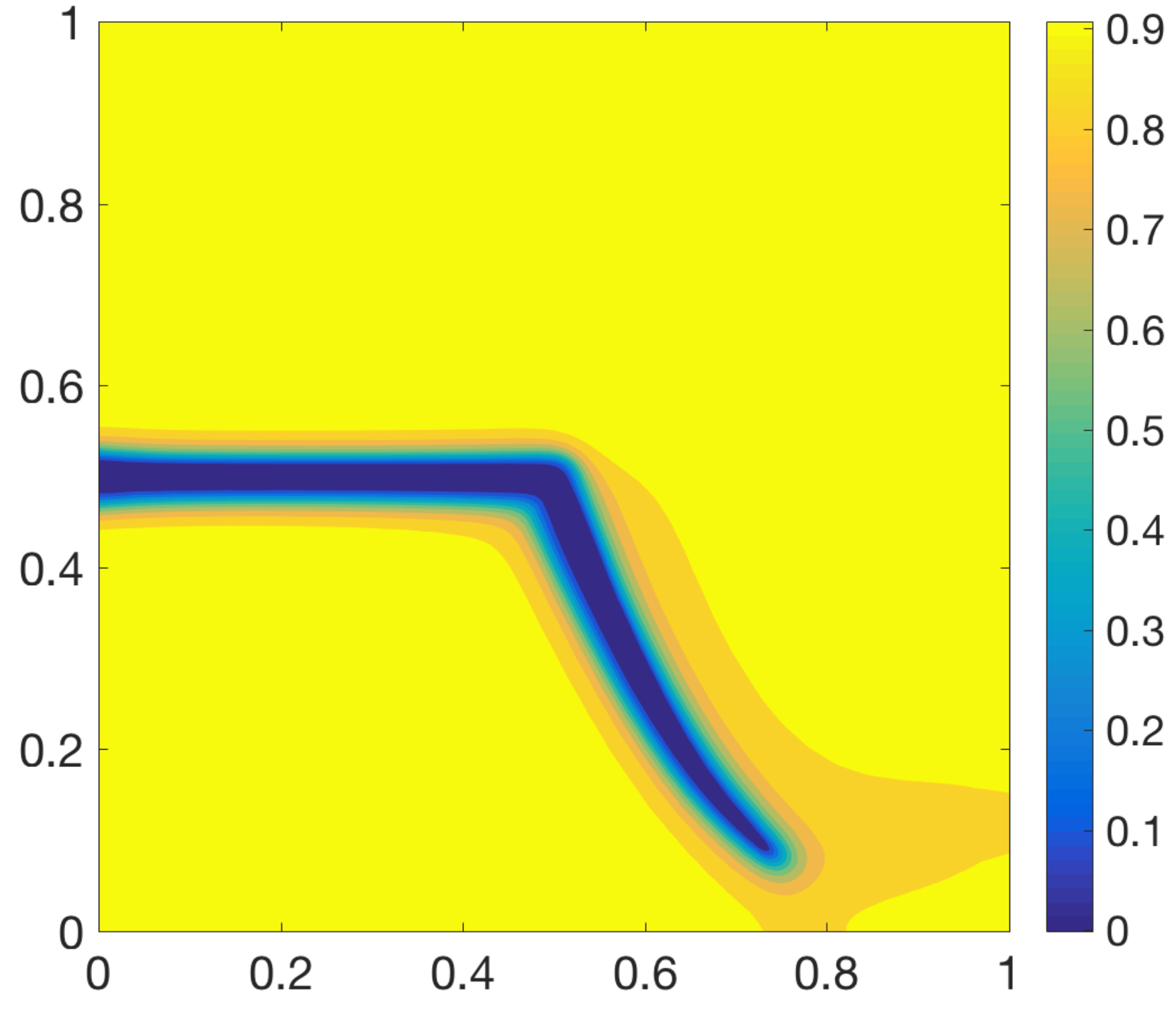}}
\subfigure[$U = 1.45 \times 10^{-2}$ mm]{\label{fig:subfig:sd1.45_b}
\includegraphics[width=0.25\linewidth]{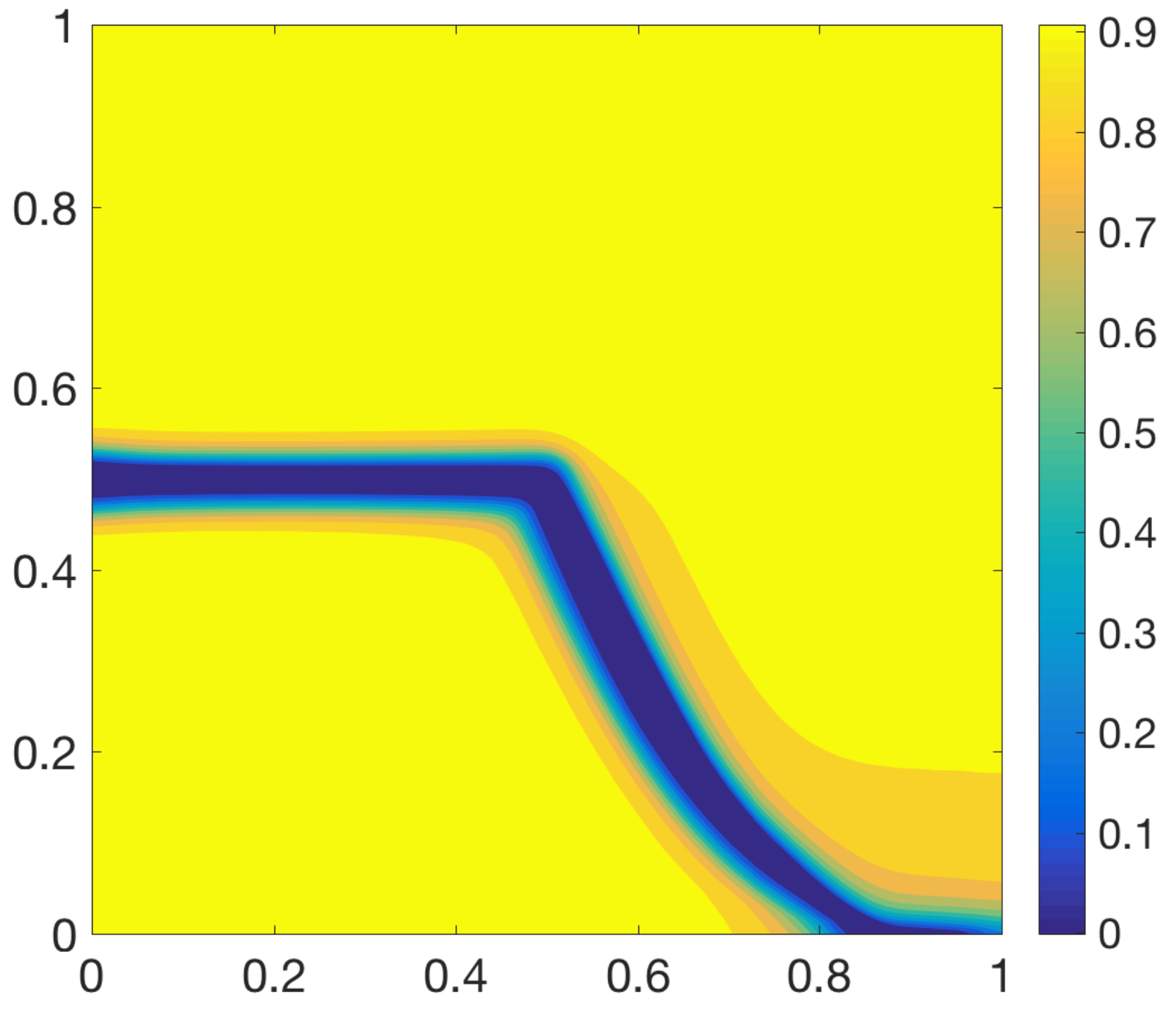}}
\caption{Example 2. The mesh and contours of the phase-field variable during crack evolution for the shear test
with $l = 0.0075$~mm and $N = 6,400$.}
\label{fig:l15_shear}
\end{figure}

\begin{figure} [htb]
\centering 
\subfigure[adaptive mesh]{\label{fig:subfig:sm1_0}
\includegraphics[width=0.29\linewidth]{./SM_l75_u1_3.pdf}}
\subfigure[close view around the crack]{\label{fig:subfig:Scrack_around}
\includegraphics[width=0.3\linewidth]{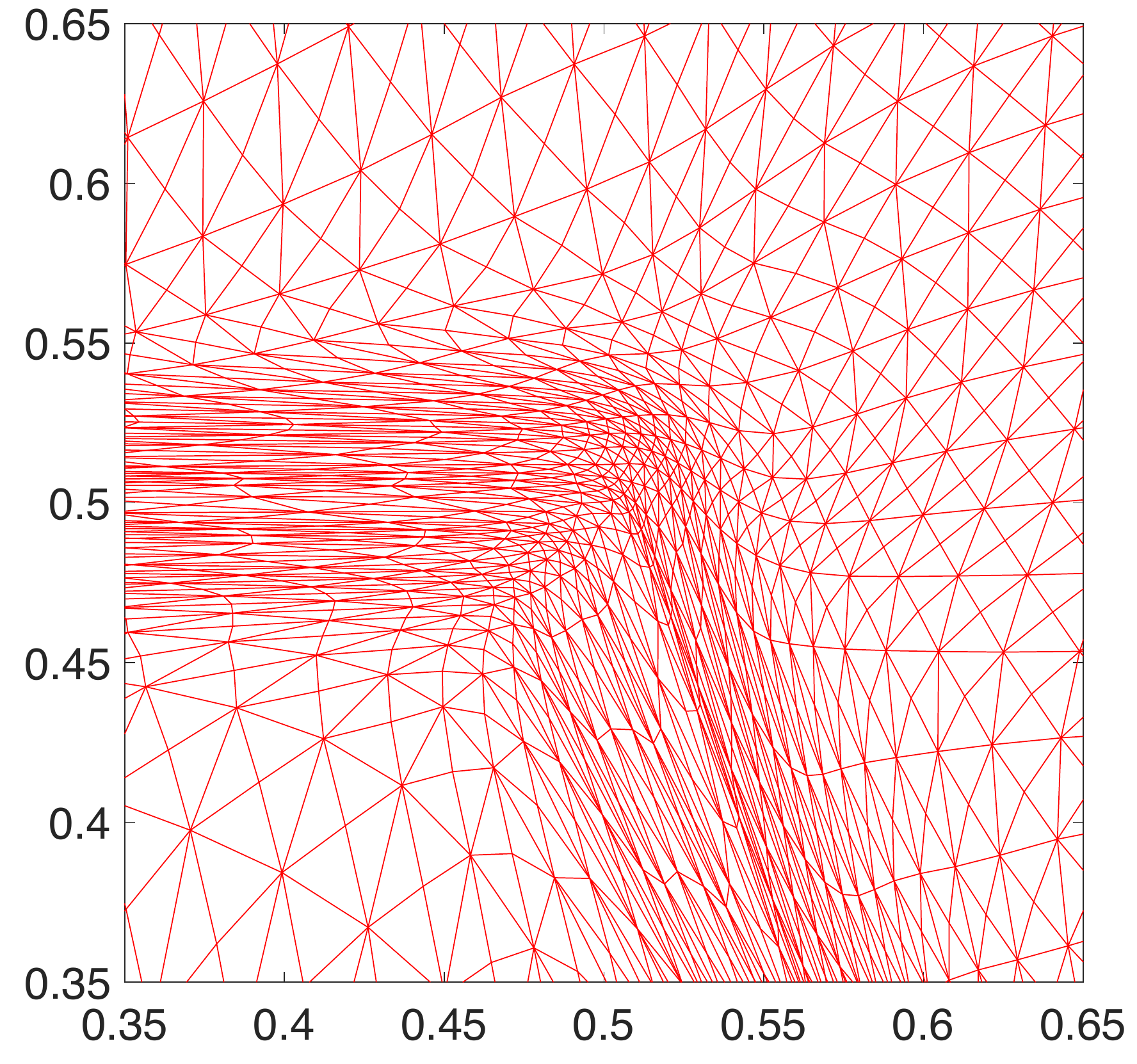}}
\subfigure[close view near the crack tip]{\label{fig:subfig:Scrack_tip}
\includegraphics[width=0.29\linewidth]{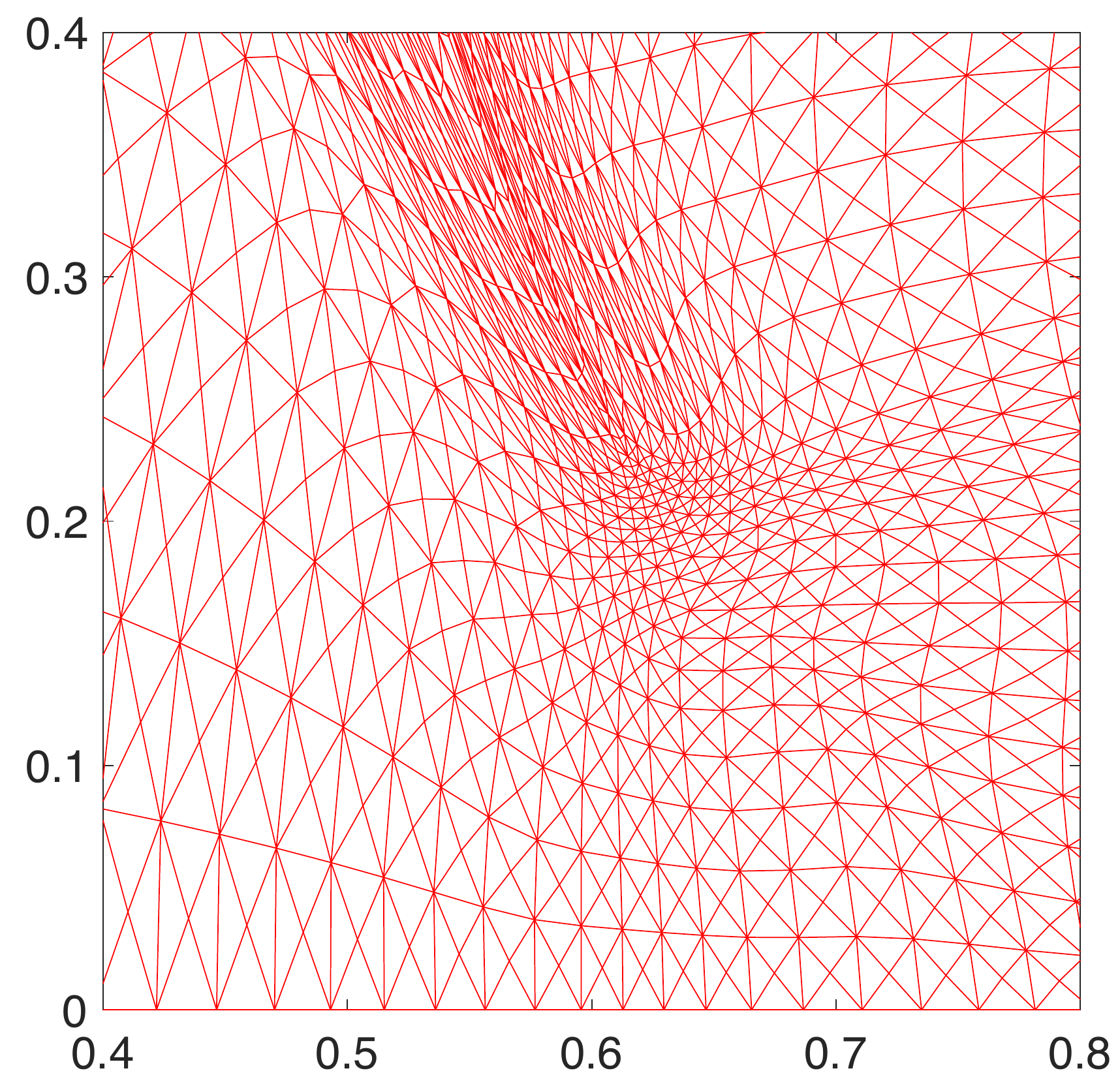}}
\caption{Example 2. The mesh and its close views for the shear crack problem.
($l = 0.00375$ mm and $N=6,400$)}
\label{fig:shear-initial-mesh}
\end{figure}

\begin{figure} [htb]
\centering 
\subfigure[$l = 0.00375$ mm]{\label{fig:subfig:Sld_l75}
\includegraphics[width=0.4\linewidth]{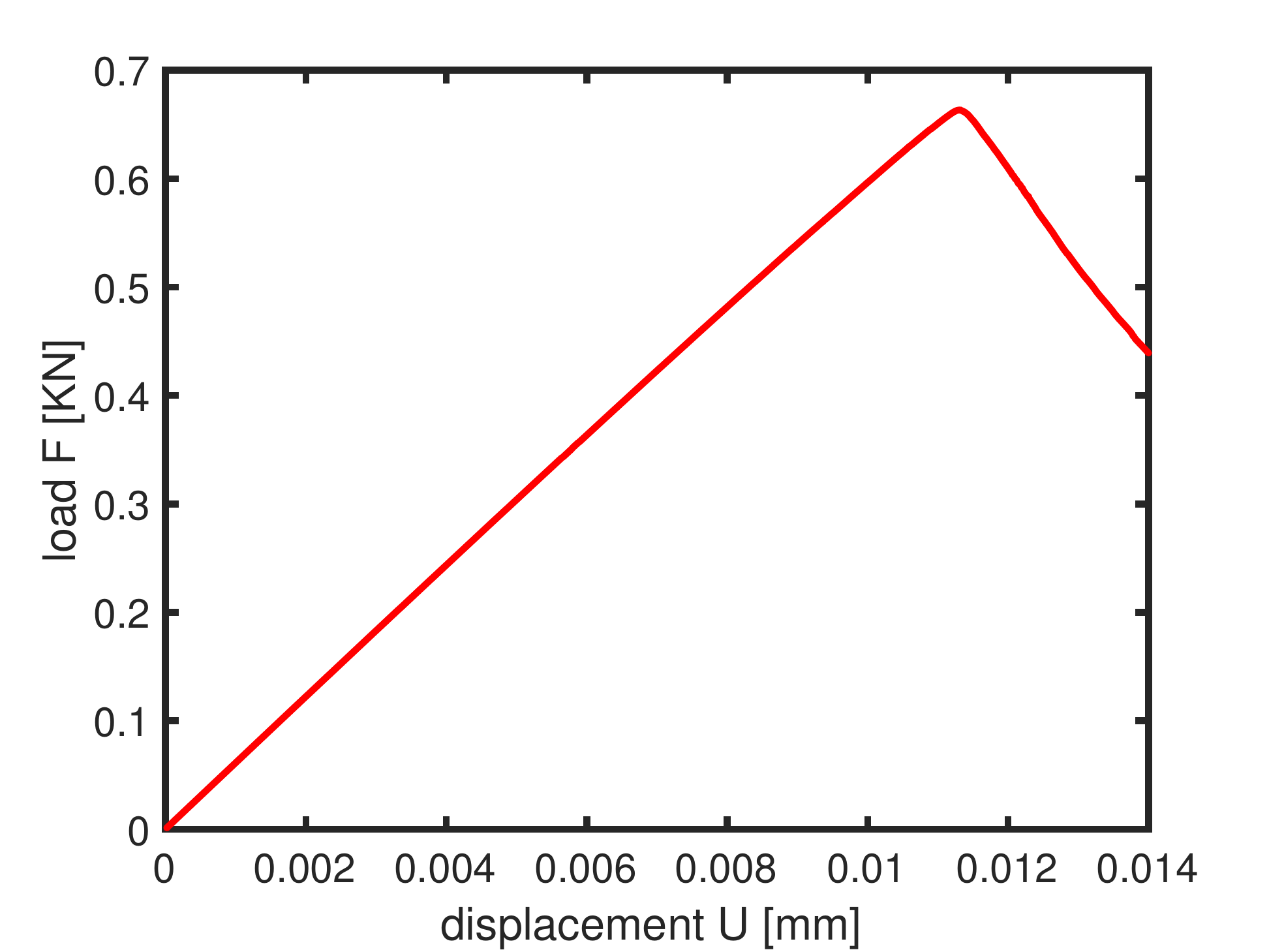}}
\subfigure[$l = 0.0075$ mm]{\label{fig:subfig:Sld_l15}
\includegraphics[width=0.4\linewidth]{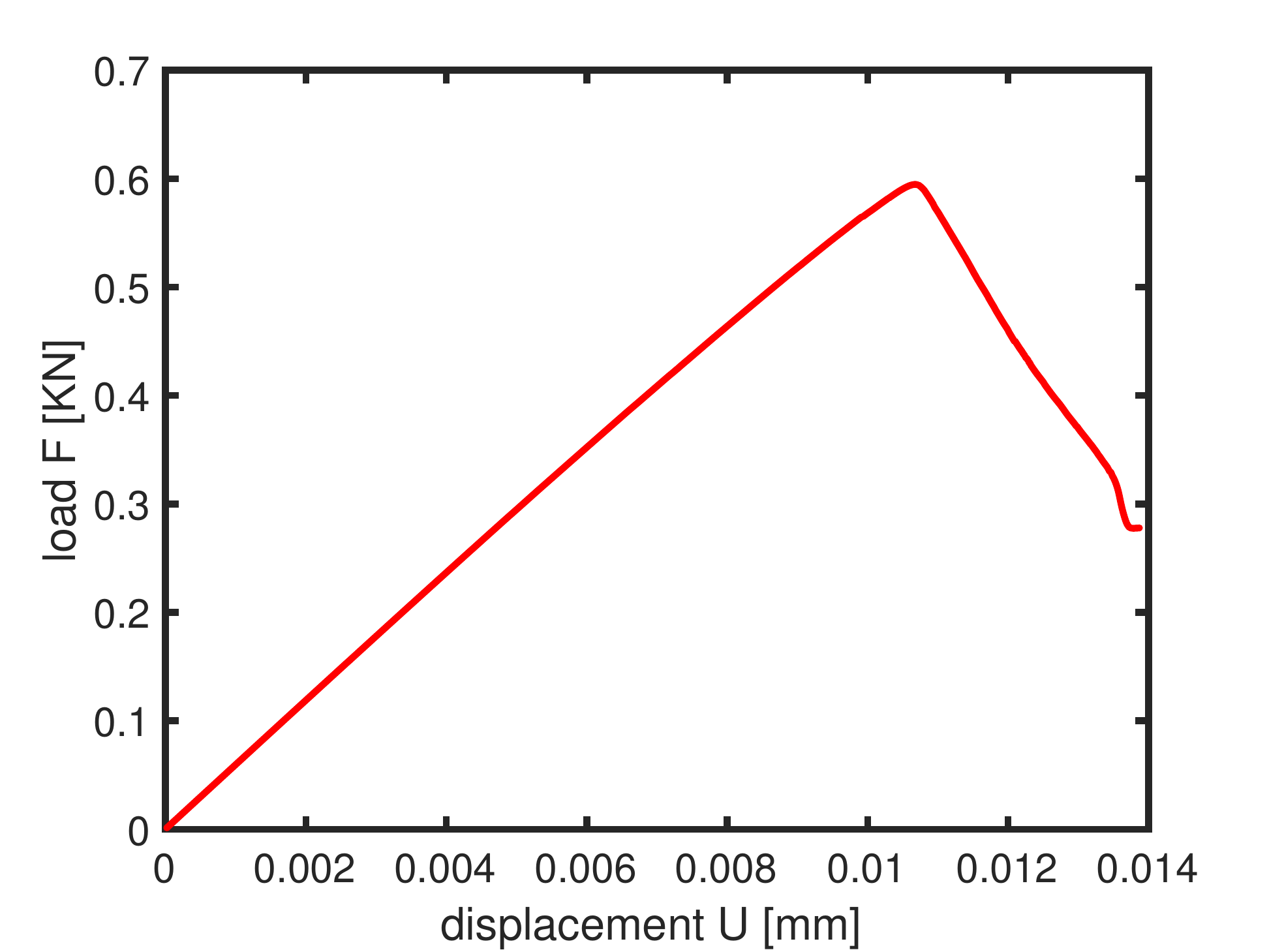}}
\caption{Example 2. The load-deflection curves for the shear crack problem for different values of $l$. ($N=6,400$)}
\label{fig:Sld_Diffl}
\end{figure}

Typical adaptive meshes and contours of the phase-field variable during crack evolution for $l = 0.00375$~mm and
$0.0075$~mm are shown in Fig. \ref{fig:l75_shear} and Fig. \ref{fig:l15_shear}, respectively. The mesh
concentrates dynamically around the crack as it evolves under continuous load. Close views of the meshes around
the crack and crack tip are plotted in Fig. \ref{fig:shear-initial-mesh}. The mesh concentration is adequate
especially during the turning process of the shear crack. This property is important for handling large fracture
networks with complex crack propagation such as joining, branching and nonplanar propagation.
Fig. \ref{fig:Sld_Diffl} shows the load-deflection curves for the shear test with $l = 0.00375$~mm
and $0.0075$~mm using a moving mesh of $N=6,400$. 

\begin{figure} [!htb]
\centering 
\subfigure[No regularization]{\label{fig:subfig:snon_regularization}
\includegraphics[width=0.36\linewidth]{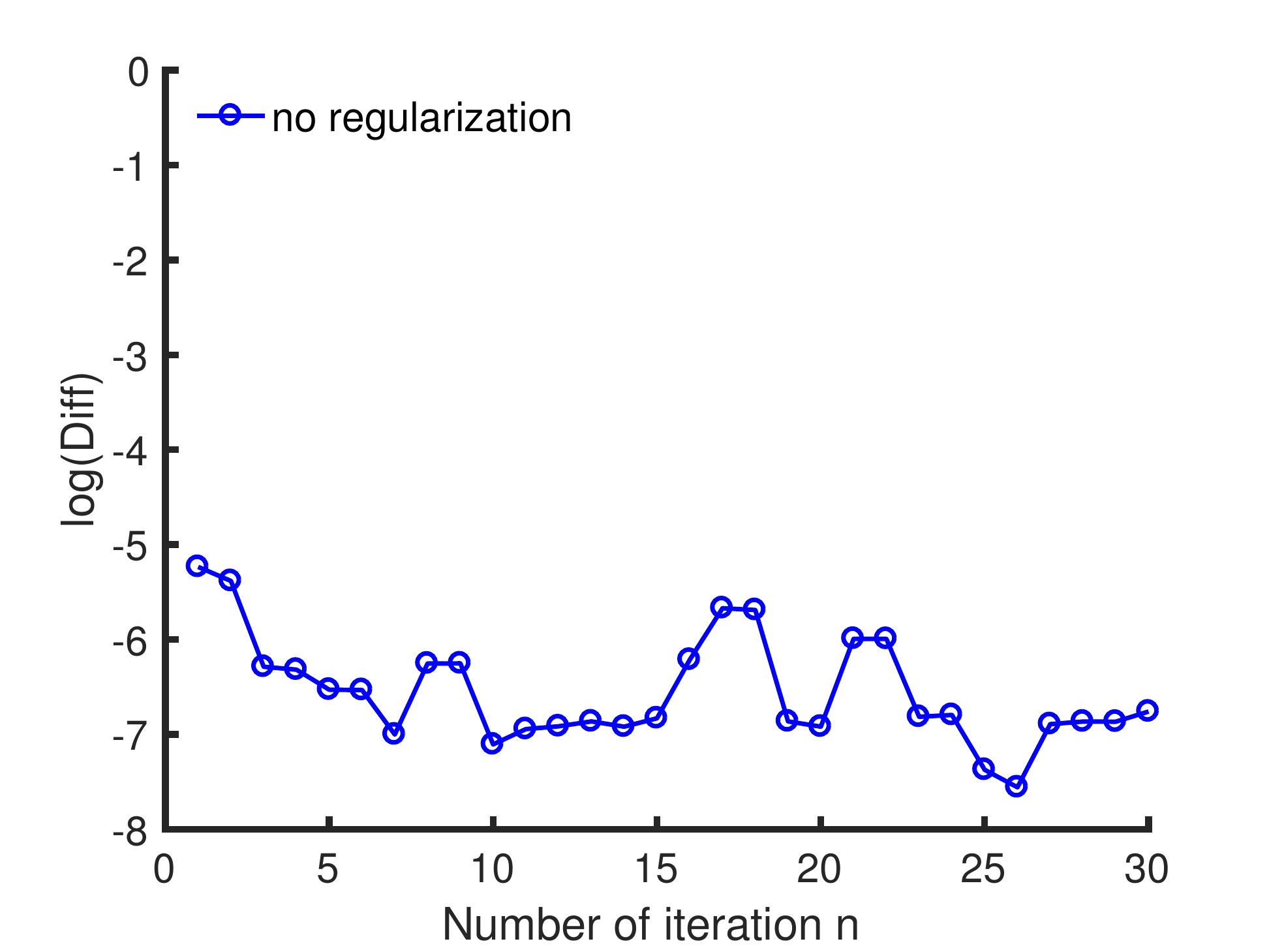}}
\subfigure[Sonic-point]{\label{fig:subfig:ssonic_point}
\includegraphics[width=0.36\linewidth]{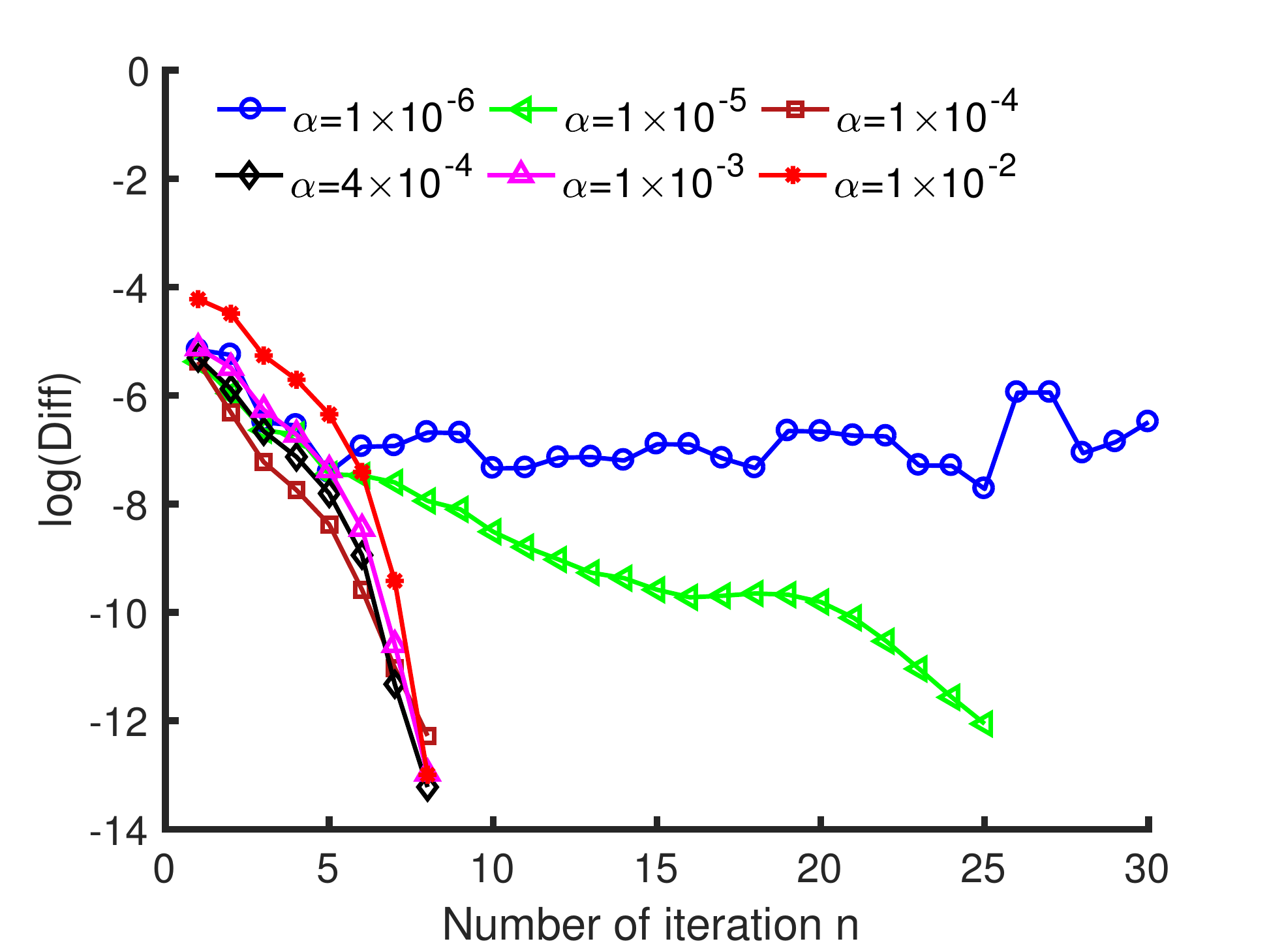}}
\vfill
\subfigure[Exponential convolution]{\label{fig:subfig:sexponential_convolution}
\includegraphics[width=0.36\linewidth]{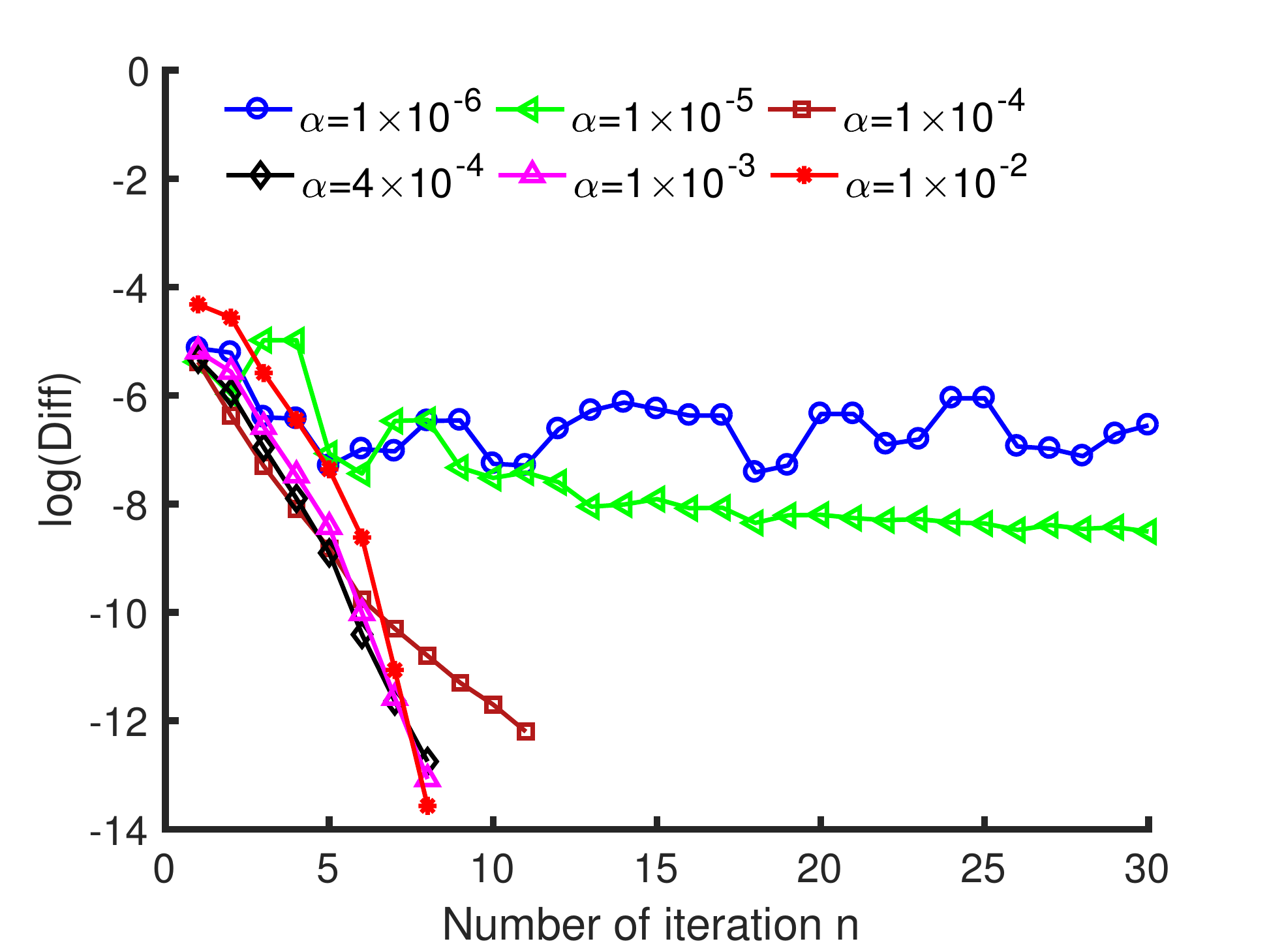}}
\subfigure[Smoothed 2-point convolution]{\label{fig:subfig:ssmoothed_2_point_convolution}
\includegraphics[width=0.36\linewidth]{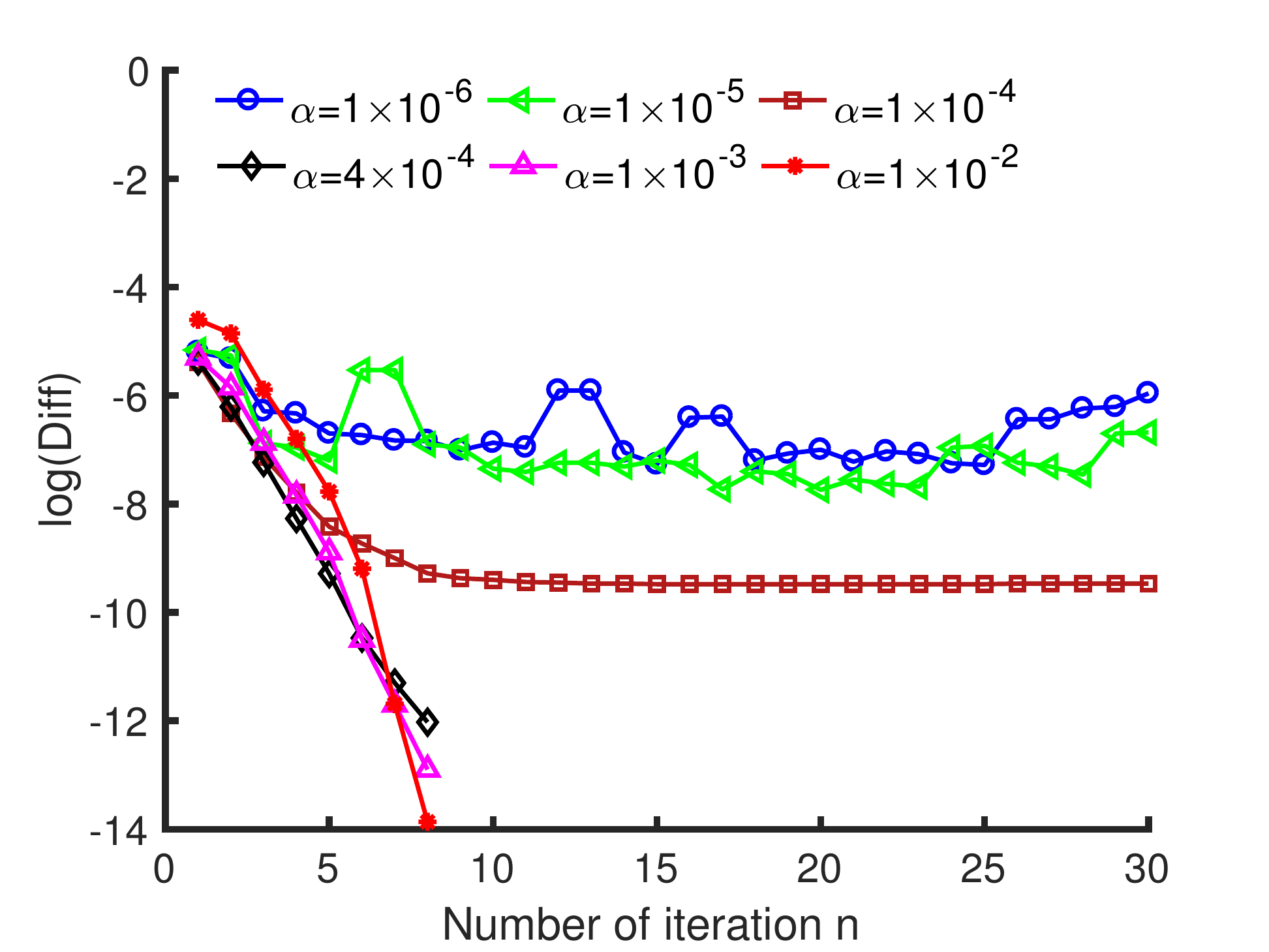}}
\caption{Example 2. Convergence of Newton's iteration using different regularization methods
for the shear test where Diff denotes the $L^2$ norm of the difference between two consecutive approximations.
($l = 0.0075$~mm and $N =6,400$)}
\label{fig:shear_convergence}
\end{figure}

Next, we compare the regularization methods for the shear test problem. Fig. \ref{fig:shear_convergence}
shows the convergence of Newton's iteration using the three regularization methods.
As we can see, all of the methods make Newton's iteration convergent for
$\alpha \ge 4\times10^{-4}$. The convergence improves as $\alpha$ increases.
Moreover, Newton's iteration converges even for $\alpha = 1\times10^{-4}$ when the sonic-point and exponential
convolution methods are used.

\begin{figure} [!htb]
\centering 
\subfigure[$\alpha = 10^{-3}$ for three regularization methods,
without regularization $\alpha = 0$ ($k_l = 10^{-3})$]{\label{fig:subfig:Sld_Diffmeth}
\includegraphics[width=0.45\linewidth]{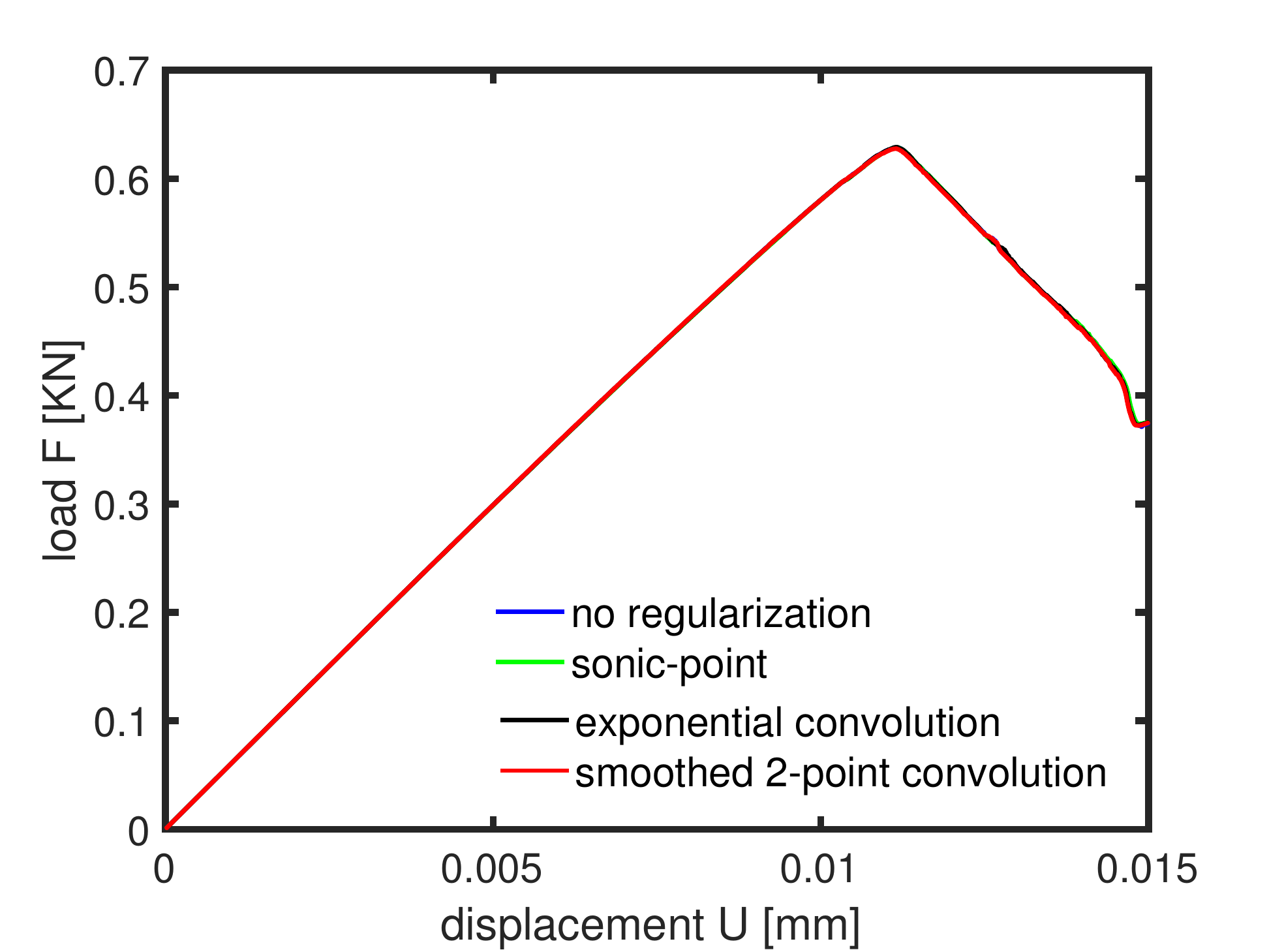}}
\subfigure[sonic-point ($k_l = 0$)]{\label{fig:subfig:Sld_job1}
\includegraphics[width=0.45\linewidth]{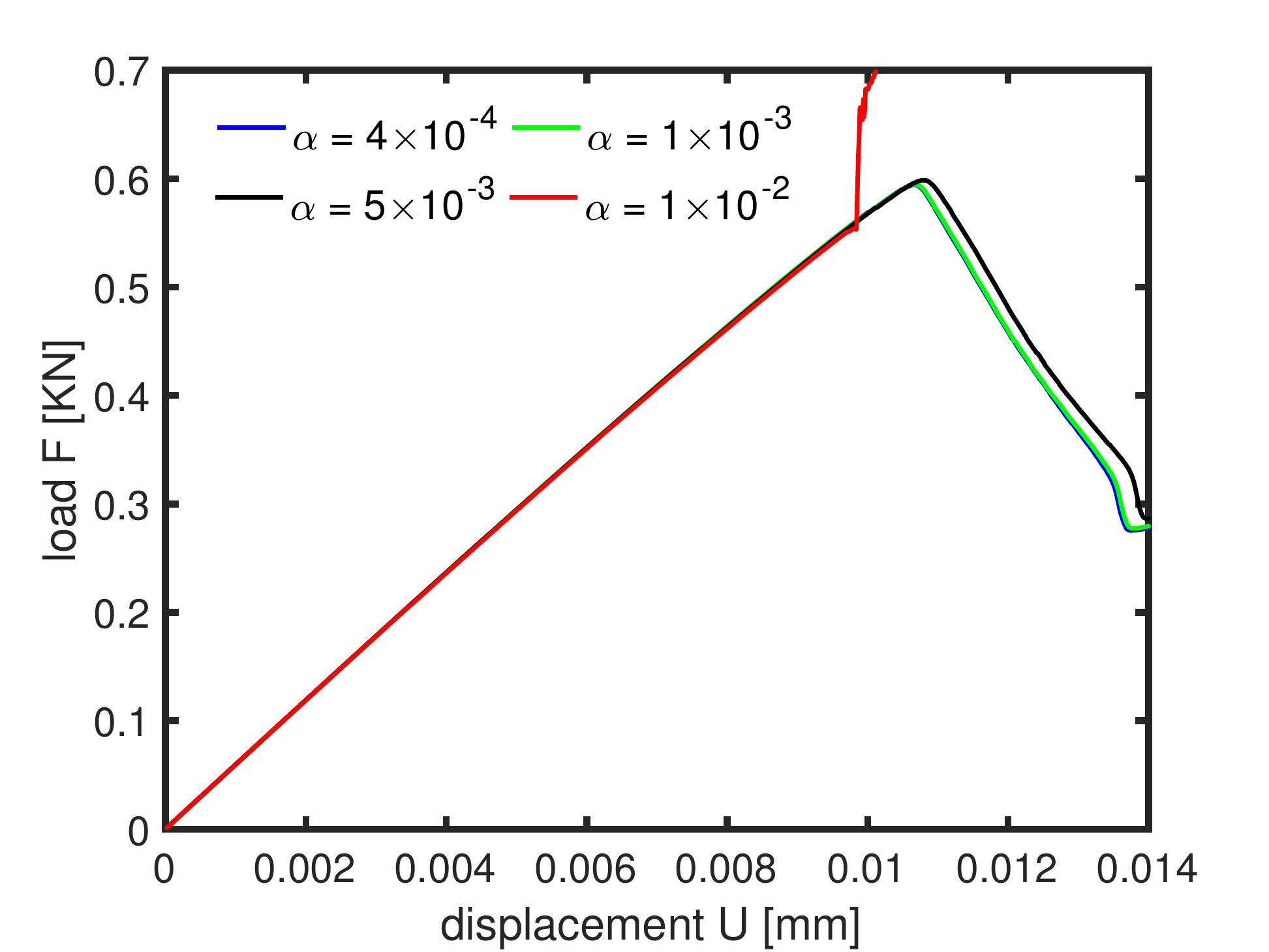}}
\subfigure[exponential convolution ($k_l = 0$)]{\label{fig:subfig:Sld_job2}
\includegraphics[width=0.45\linewidth]{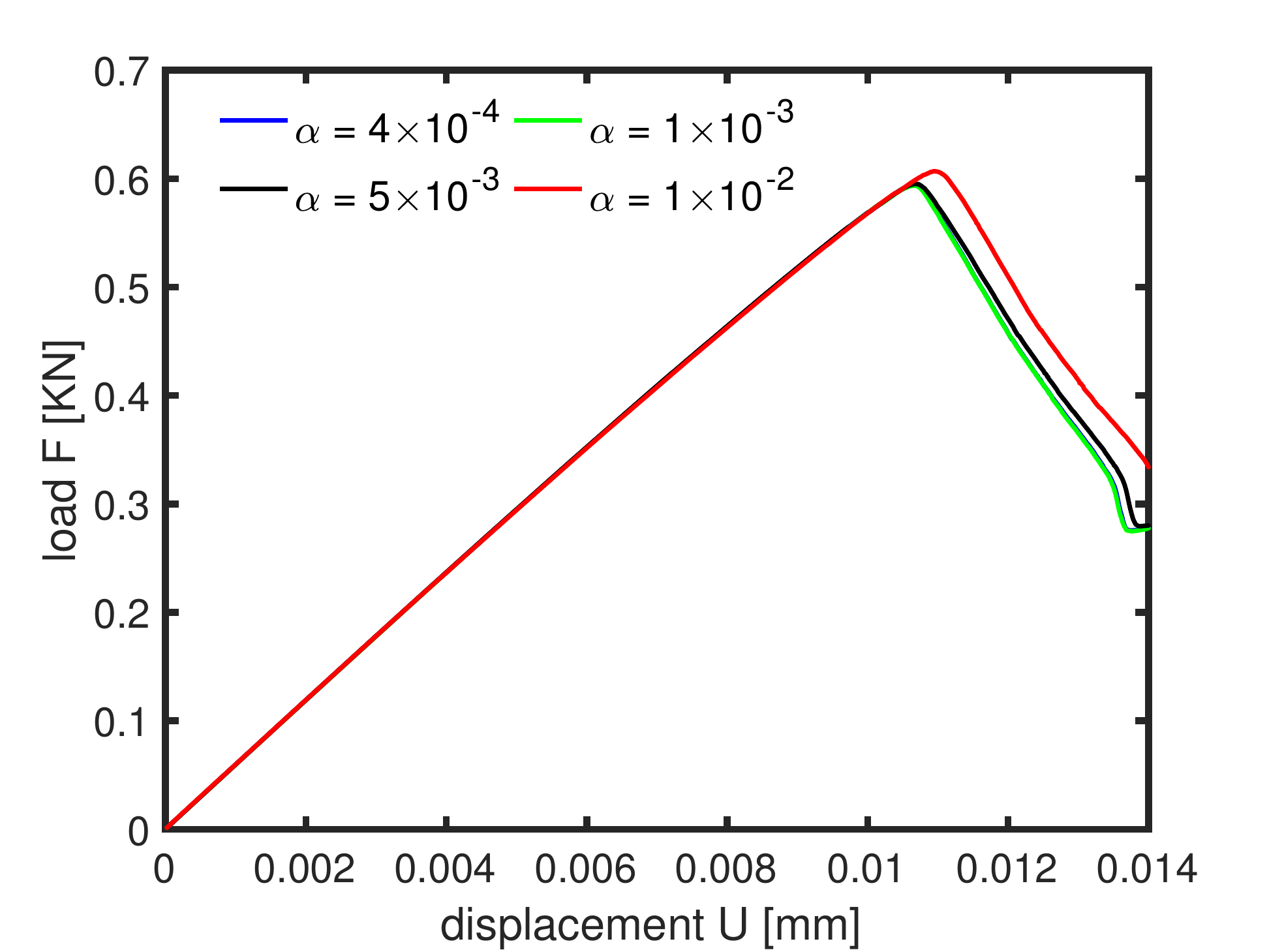}}
\subfigure[smoothed 2-point convolution ($k_l = 0$)]{\label{fig:subfig:Sld_job3}
\includegraphics[width=0.45\linewidth]{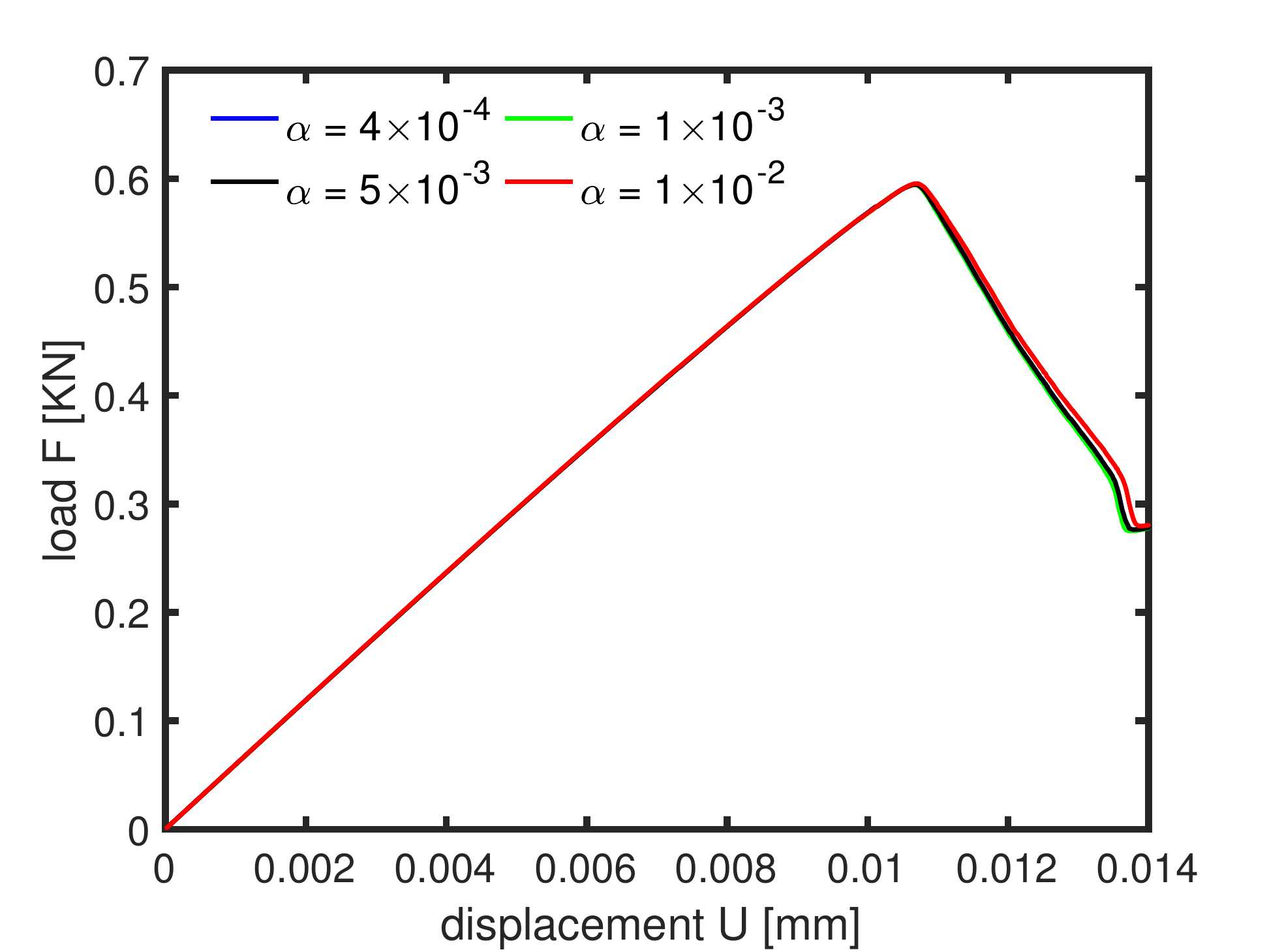}}
\caption{Example 2. The load-deflection curves for the regularization methods
for the shear test with $l=0.0075$~mm and $N = 6,400$. (a) Different regularization methods
with $\alpha = 10^{-3}$; (b) sonic-point method with various $\alpha$; (c)
exponential convolution method with various $\alpha$; (d) smoothed 2-point convolution method
with various $\alpha$.}
\label{fig:shear-effects-of-the-regularization-method}
\end{figure}

Lastly, the load-deflection curves for the shear test with $l=0.0075$~mm and $k_l = 1\times10^{-3}$
are shown in Fig. \ref{fig:subfig:Sld_Diffmeth} for the regularization methods and without regularization.
We can see that the curves for the regularization methods with $\alpha = 1\times10^{-3}$
are almost the same as that without regularization.
The load-deflection curves obtained using the three regularization methods with various values of $\alpha$
are presented in Fig. \ref{fig:subfig:Sld_job1}, \ref{fig:subfig:Sld_job2} and \ref{fig:subfig:Sld_job3}, respectively.
For the sonic-point regularization method, the load-deflection curve becomes unphysical
when $\alpha \ge 1\times10^{-2}$.
It is also interesting to observe that the load for the exponential convolution
and smoothed 2-point convolution regularization is overestimated after crack starts propagating
for large values of $\alpha$ ($\alpha \ge 5\times10^{-3}$).  The effects of $\alpha$ are less significant
with smoothed 2-point convolution than the other two methods. 

\subsection{Example 3. Test with multiple cracks}

\begin{figure} [!htb]
\centering 
\subfigure[junction between two cracks]{\label{fig:subfig:junctionTwo}
\includegraphics[width=0.3\linewidth]{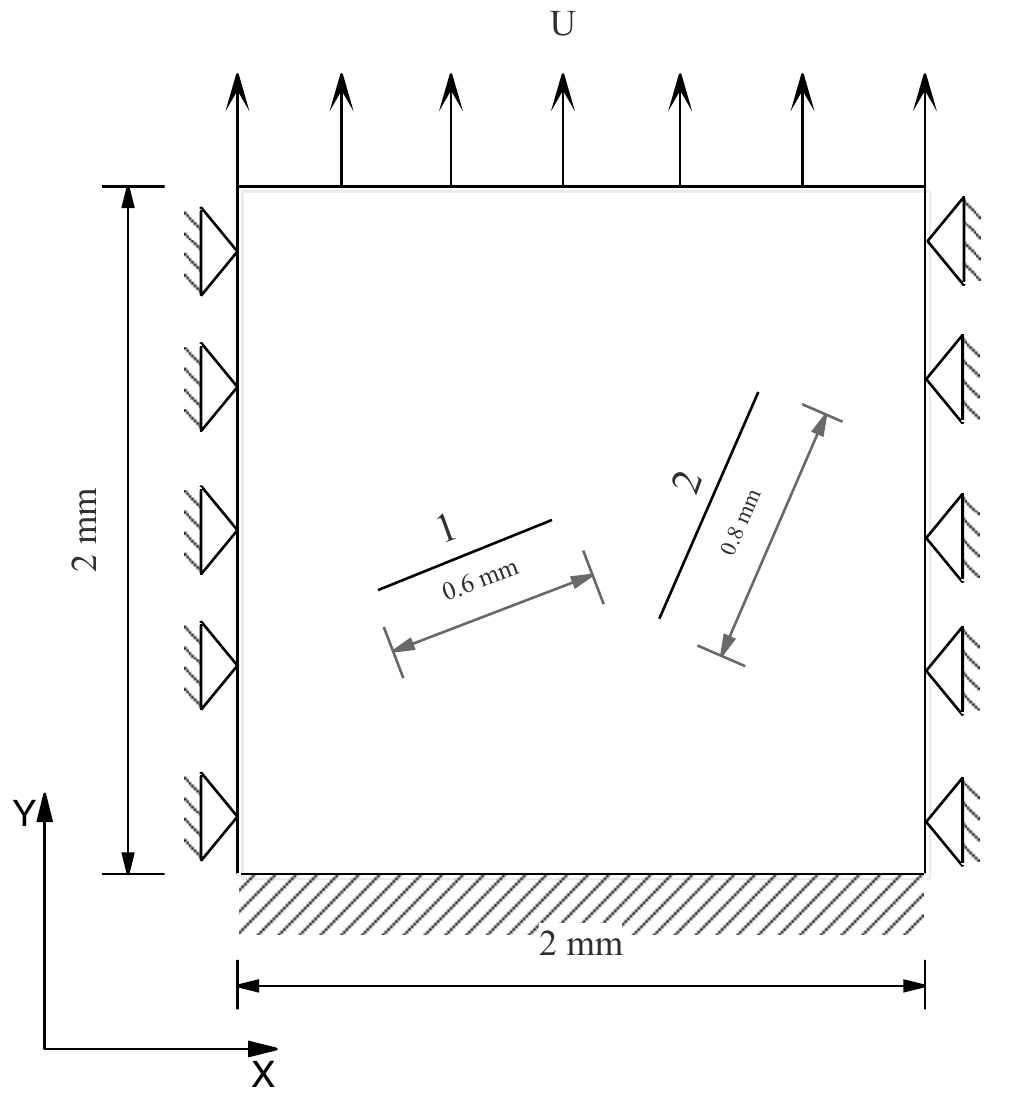}}
\subfigure[random distribution of five cracks]{\label{fig:subfig:junctionMult}
\includegraphics[width=0.3\linewidth]{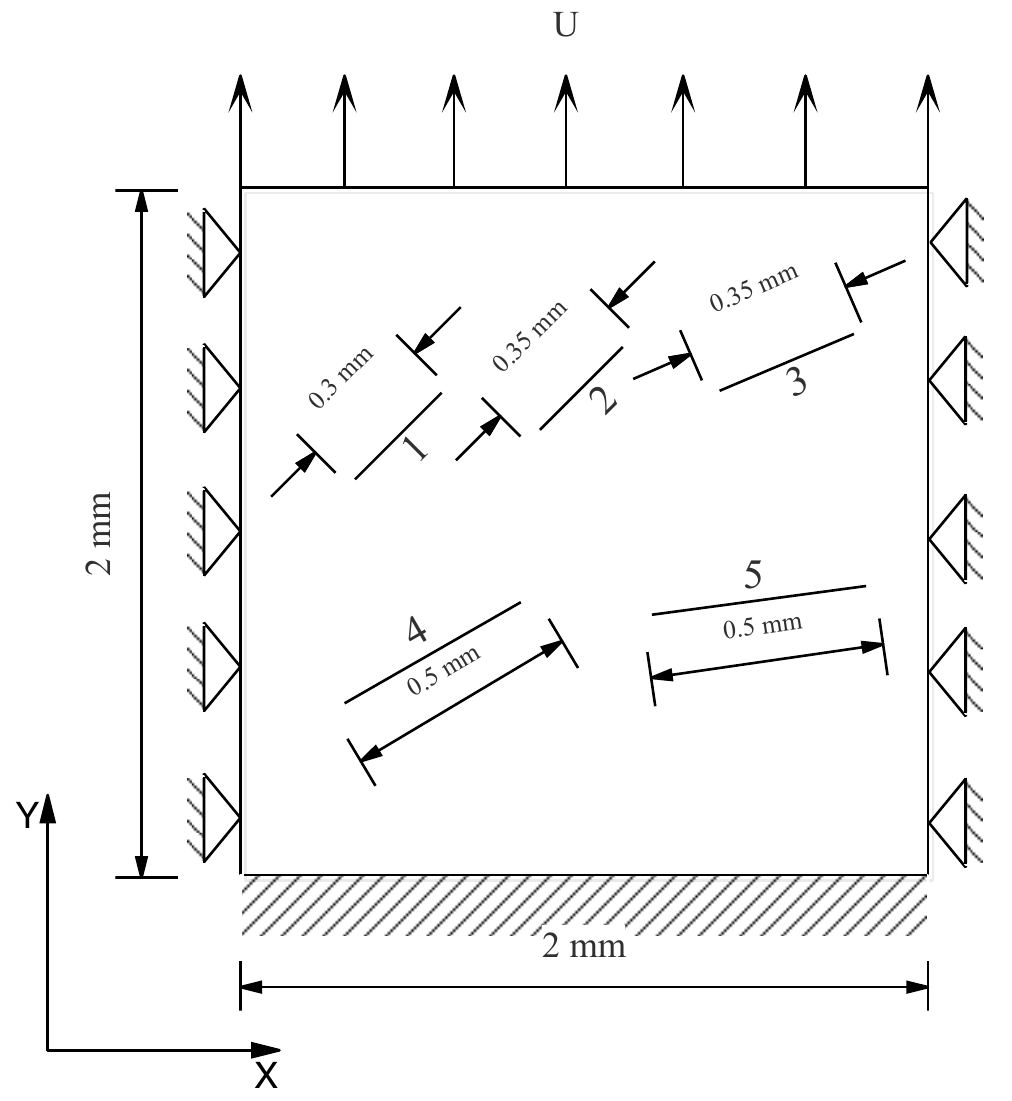}}
\subfigure[random distribution of ten cracks]{\label{fig:subfig:junctionTen}
\includegraphics[width=0.3\linewidth]{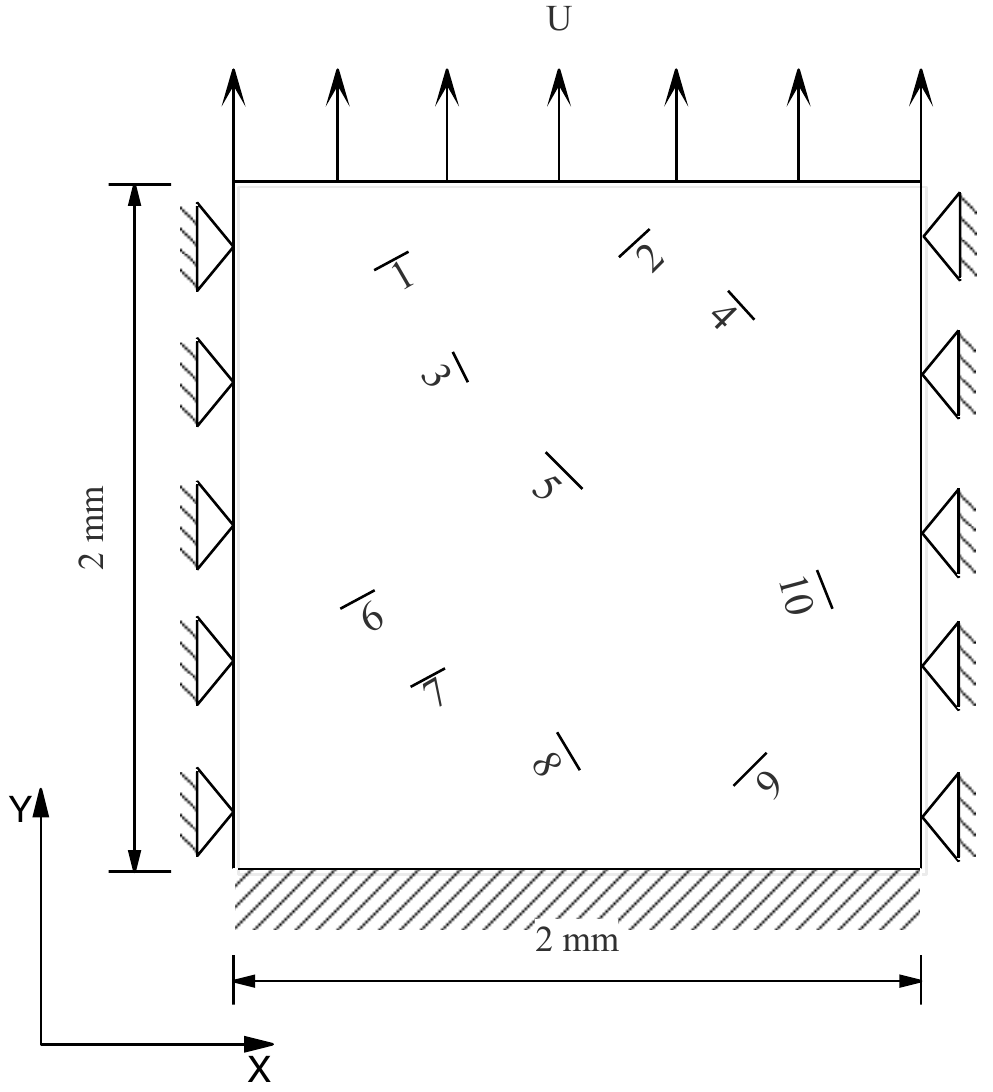}}
\caption{Example 3. Domain and boundary conditions for the test with multiple cracks. (a) two cracks, (b) five cracks,
(c) ten cracks.}
\label{fig:multiple cracks}
\end{figure}

In this example, we first consider two cracks in a square plate of width $2$~mm,
with the domain and boundary conditions shown 
in Fig. \ref{fig:subfig:junctionTwo}, to test the modeling of the junction between two cracks. 
Crack 1 is centered at ($-0.2,0$) with length $0.6$~mm and polar angle (relative to the horizontal direction) $9^{\circ}$,
while Crack 2 is centered at ($0.46,0$) with length $0.8$~mm and  polar angle $65^{\circ}$.
The bottom edge of the domain is fixed and the top edge is fixed along $x$-direction 
while a uniform $y$-displacement $U$ is increased with time to drive the crack propagation. 
The material parameters are the same as previous examples, that is, elastic constants
$\lambda = 121.15$~kN/mm$^{2}$, $\mu = 80.77$~kN/mm$^{2}$, 
and the fracture toughness is $g_c = 2.7 \times 10^{-3} $~kN/mm.
The displacement increment is chosen as $\Delta U = 1\times10^{-4}$~mm for the computation.
The mesh consists of $10,000$ triangular elements and the length scale parameter is chosen as $l = 0.00375$~mm.
A similar configuration with different boundary conditions has been used by Budyn et al. \cite{BZMB2004}.

\begin{figure} 
\centering 
\subfigure[$U = 1.0 \times 10^{-3}$ mm]{\label{fig:subfig:JTM_l75_u1}
\includegraphics[width=0.22\linewidth]{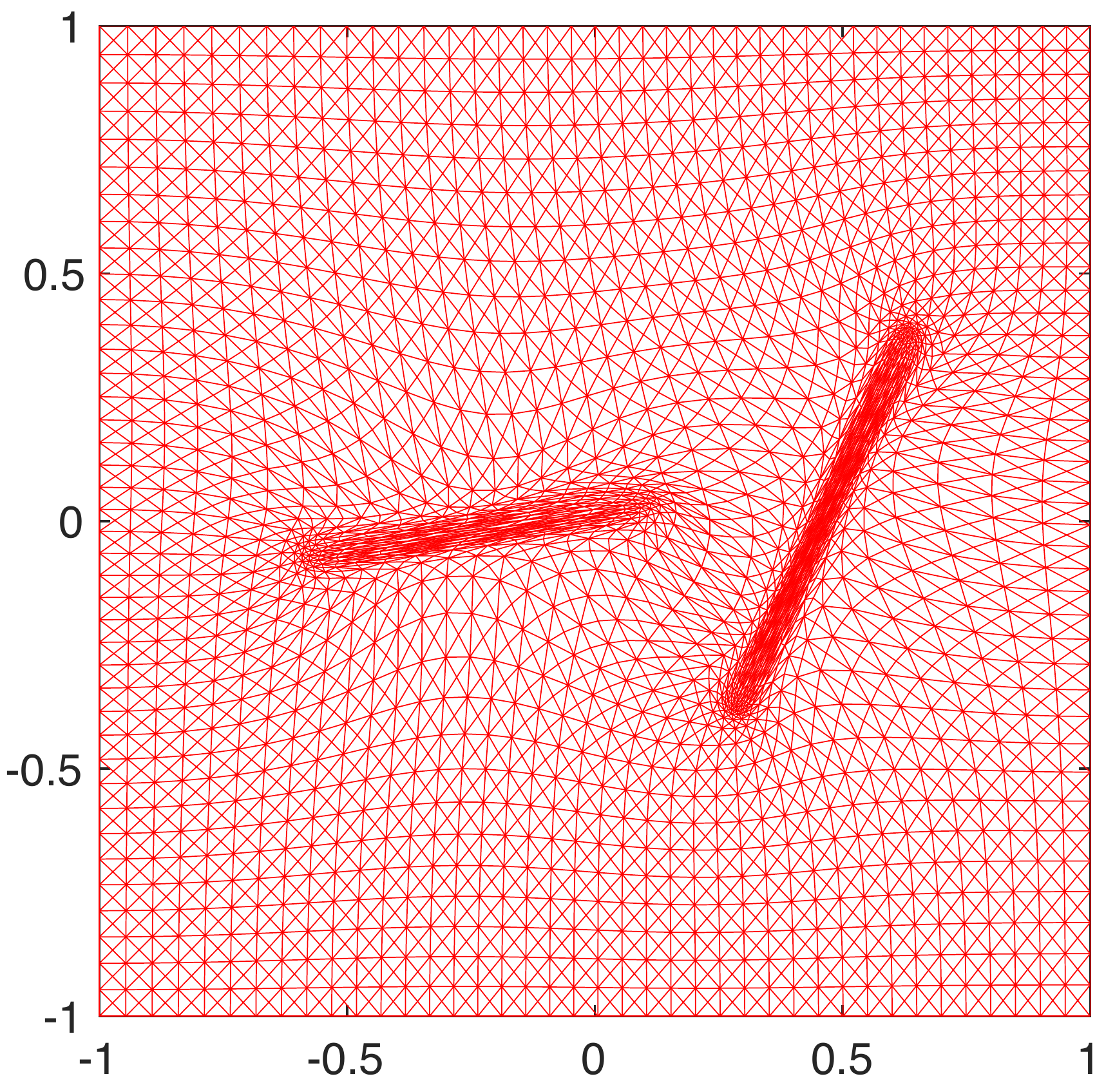}}
\subfigure[$U = 1.0 \times 10^{-2}$ mm]{\label{fig:subfig:JTM_l75_u10}
\includegraphics[width=0.22\linewidth]{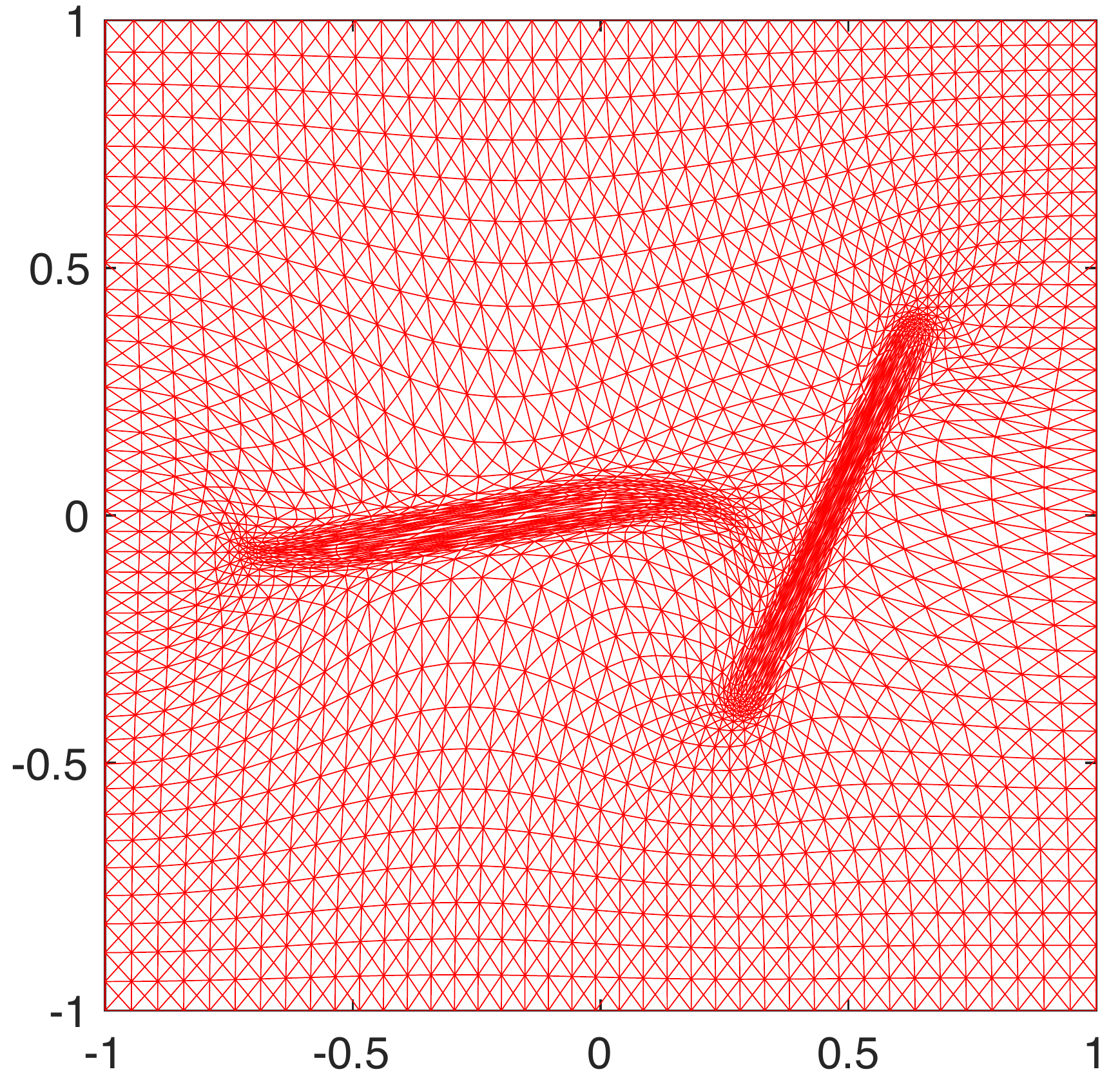}}
\subfigure[$U = 1.1 \times 10^{-2}$ mm]{\label{fig:subfig:JTM_l75_u11}
\includegraphics[width=0.22\linewidth]{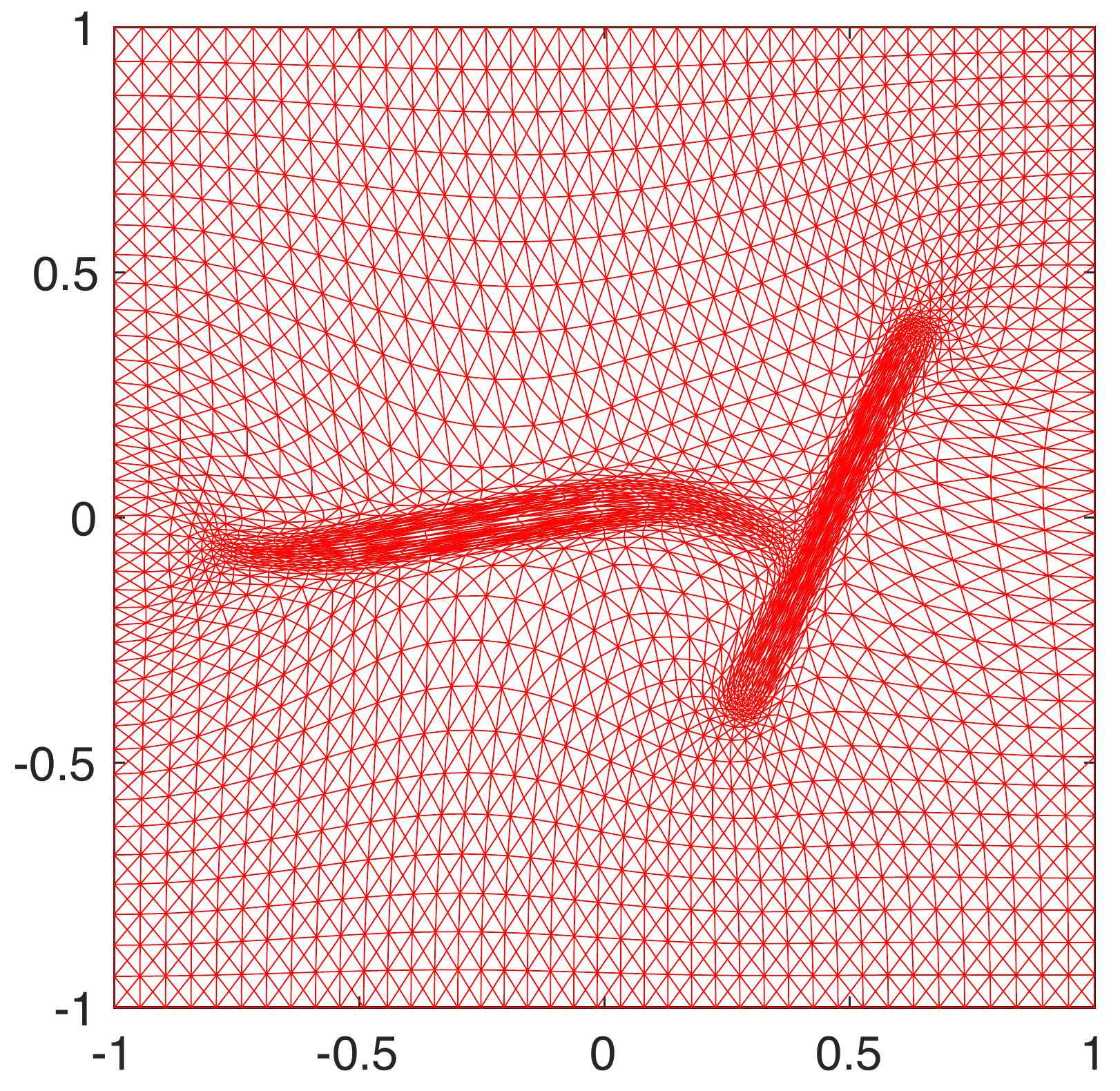}}
\subfigure[$U = 1.6 \times 10^{-2}$ mm]{\label{fig:subfig:JTM_l75_u16}
\includegraphics[width=0.22\linewidth]{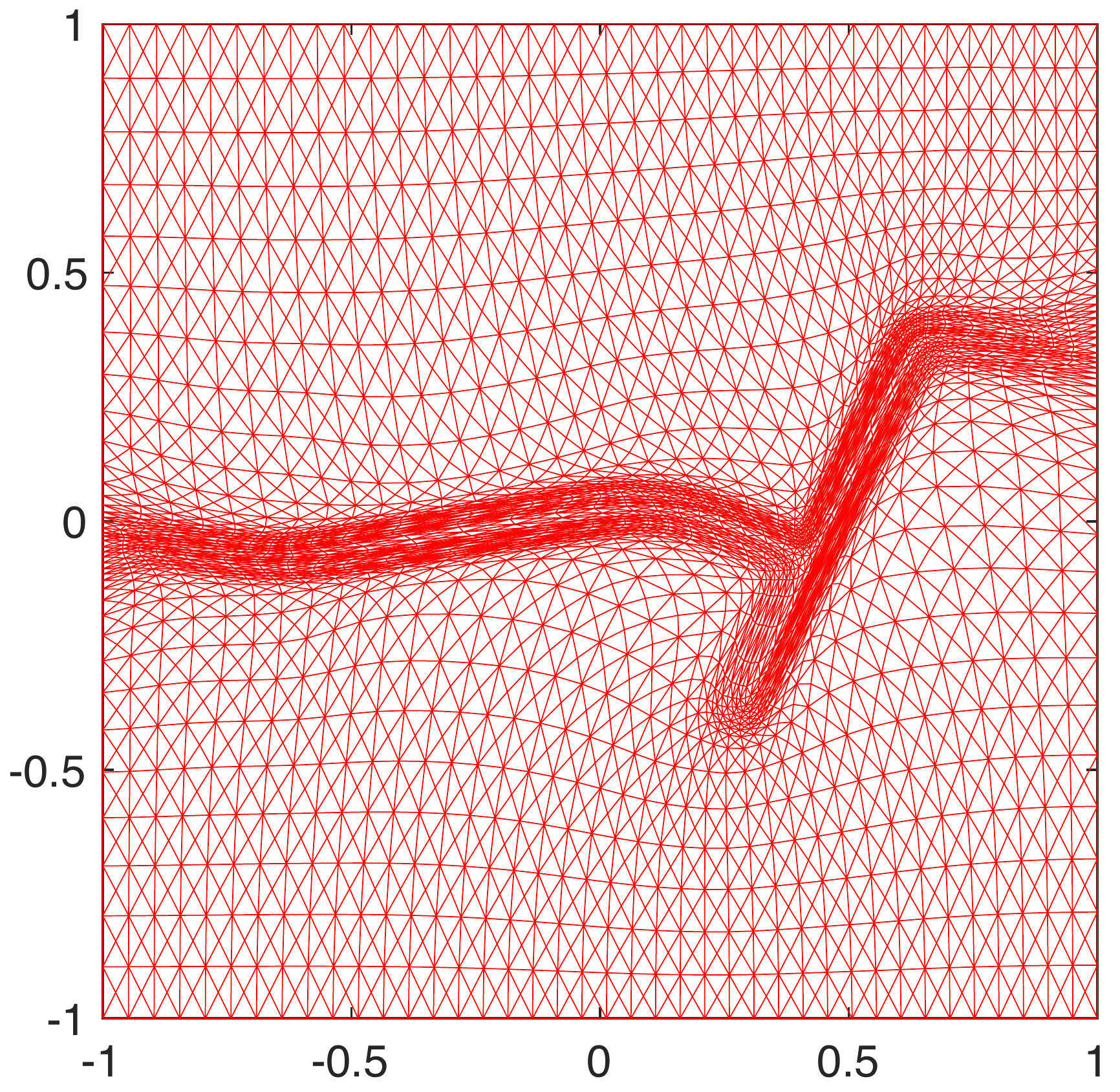}}
\vfill
\subfigure[$U = 1.0 \times 10^{-3}$ mm]{\label{fig:subfig:JTD_l75_u1}
\includegraphics[width=0.22\linewidth]{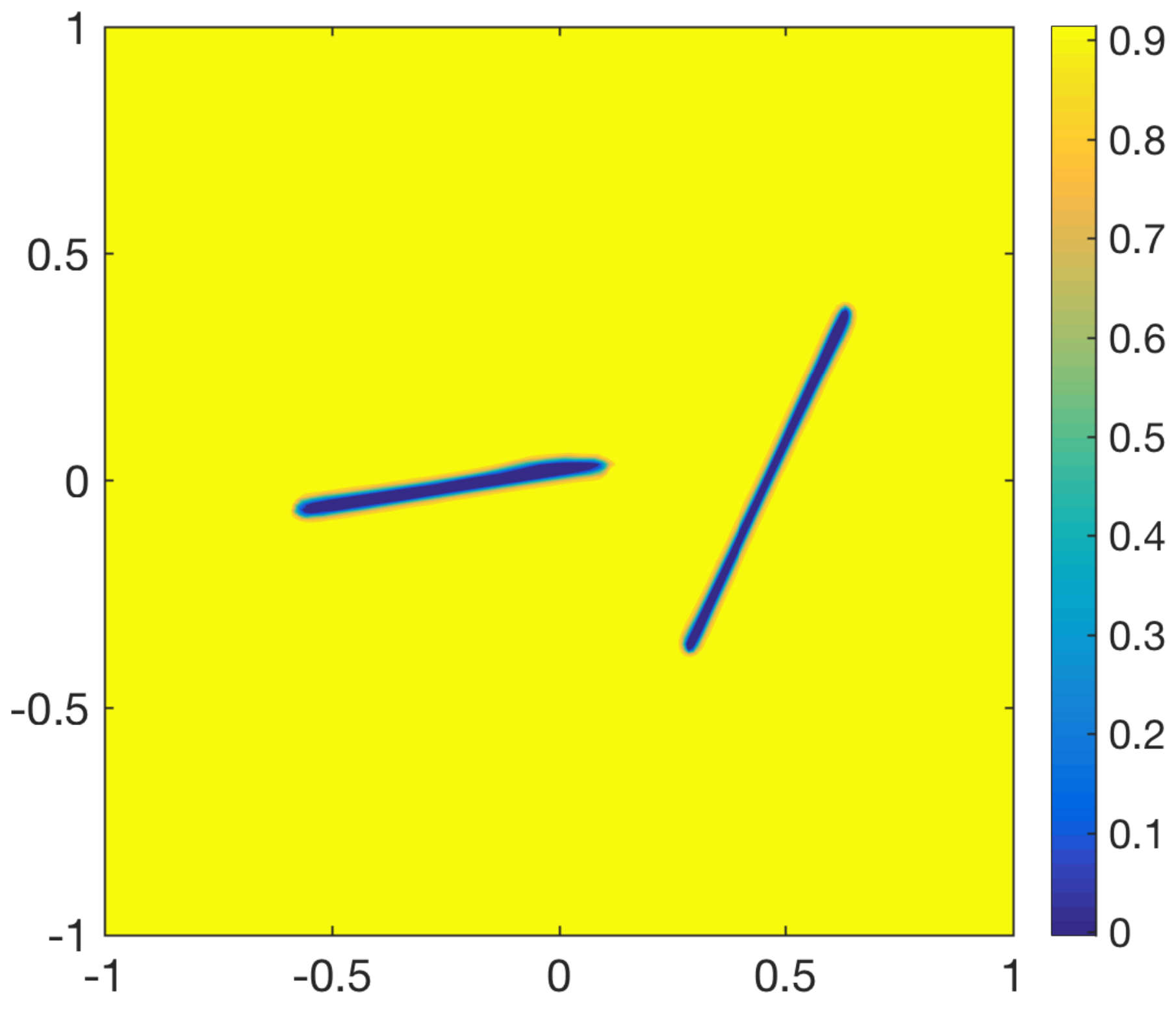}}
\subfigure[$U = 1.0 \times 10^{-2}$ mm]{\label{fig:subfig:JTD_l75_u10}
\includegraphics[width=0.22\linewidth]{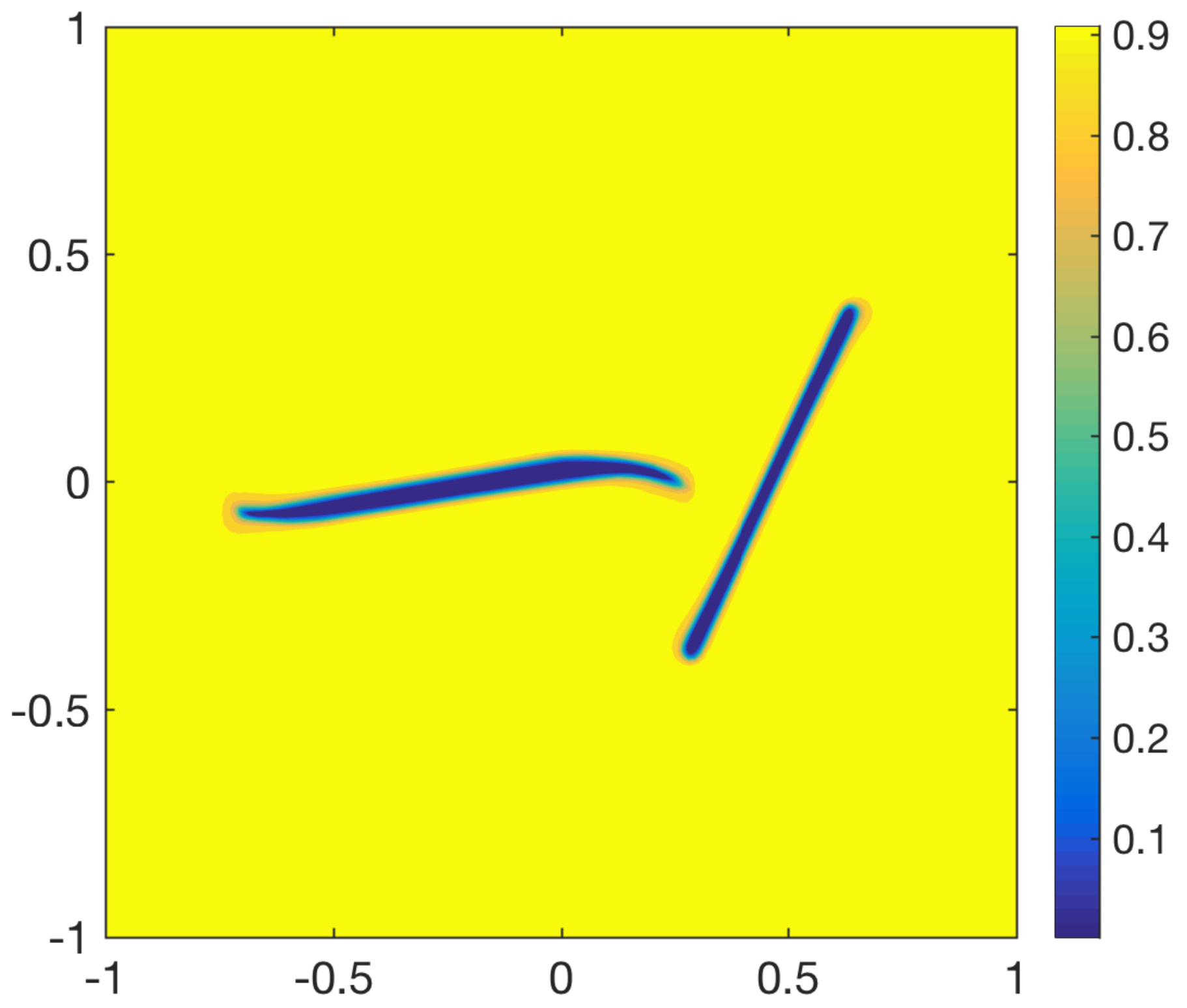}}
\subfigure[$U = 1.1 \times 10^{-2}$ mm]{\label{fig:subfig:JTD_l75_u11}
\includegraphics[width=0.22\linewidth]{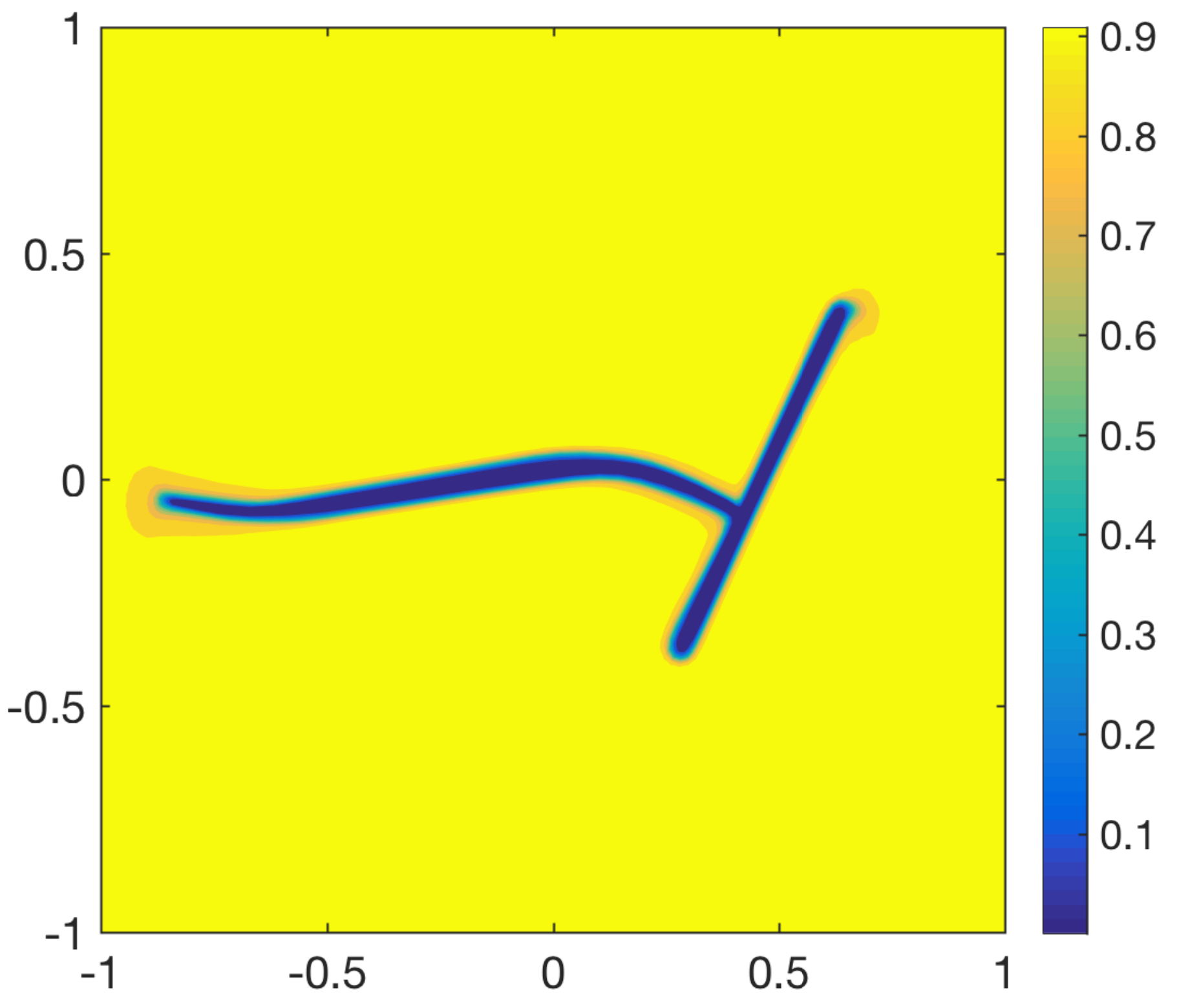}}
\subfigure[$U = 1.6 \times 10^{-2}$ mm]{\label{fig:subfig:JTD_l75_u16}
\includegraphics[width=0.22\linewidth]{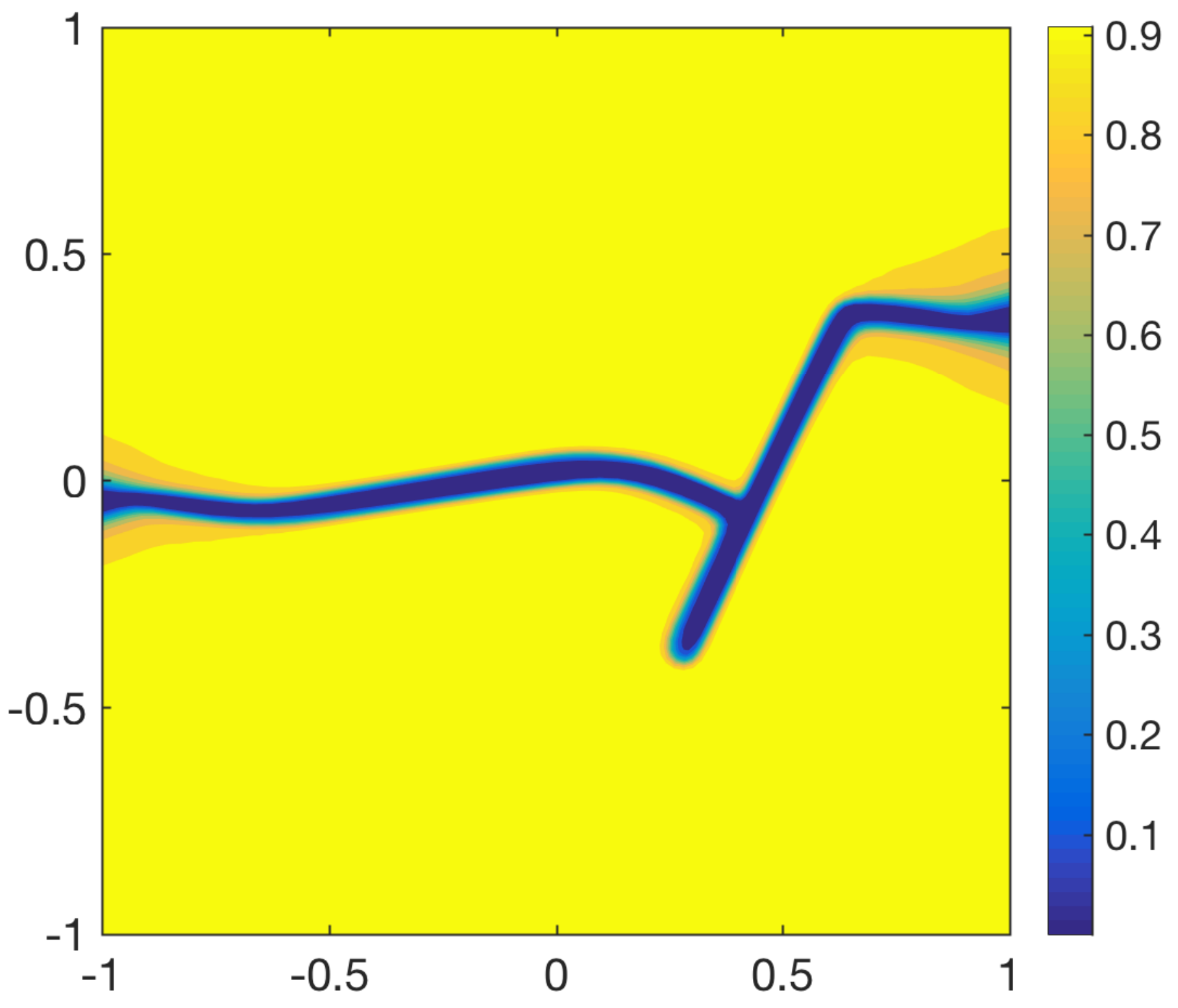}}
\caption{Example 3. The mesh and contours of the phase-field distribution during crack evolution
for the two-crack problem with $l = 0.00375$~mm and $N = 10,000$.}
\label{fig:l75_JunctionTwo}
\end{figure}

Typical adaptive meshes and contours of the phase-field variable during crack evolution are shown in Fig. \ref{fig:l75_JunctionTwo}.
As can be seen, the mesh adapts dynamically to capture the junction process of the two cracks. 
As the load $U$ increases, the evolution can be described as follows.

\begin{itemize}
\item[(a)] $U = 8 \times 10^{-3}$~mm: the tip of Crack 1 activates; 
\item[(b)] $U = 8 \times 10^{-3}$ to $1.0 \times 10^{-2}$~mm: Crack 1 propagates along a curved path
heading towards Crack 2;
\item[(c)] $U = 1.1 \times 10^{-2}$~mm: Crack 1 connects to Crack 2, and one of the tip for Crack 2 activates;
\item[(d)] $U = 1.6 \times 10^{-2}$~mm: both Crack 1 and Crack 2 have propagated to the edge of the plate.
\end{itemize}

Since there is no analytical solution for this problem, we cannot compare
the computed solutions with the exact solution. Nevertheless, they can be justified in physics.
In fact, it is known (e.g., see \cite{Per04})
that the critical stress increases as the polar angle increases and cracks with large polar angles are
harder to initiate. Therefore, Crack 1, which is oriented with a smaller polar angle, is the first to initiate
and propagate. Moreover, the frequency of crack coalescence is strongly related to the crack density
(the crack size and their relative location) within the solid; e.g., see \cite{WAAP96}. 
The higher the crack density is, the more likely two cracks will merge.
Since being close to Crack 2, Crack 1 begins to merge with Crack 2 as it grows.
At the same time, the other side of Crack 1 has propagated to the edge of the plate,
which causes the left of the plate to lose strength.
At later stages, the right side of the tip for Crack 2 activates and propagates
to the edge of the plate.

Next, we consider five initial cracks in the same domain and with the same boundary conditions
as for the two-crack problem, see Fig. \ref{fig:subfig:junctionMult}. The material parameters are the same
as previous examples except $g_c = 2.7 \times 10^{-4} $~kN/mm. 
The lengths of Crack 1, 2, 3, 4, and 5 are $0.3$~mm, $0.35$~mm, $0.35$~mm, $0.5$~mm, and $0.5$~mm, 
centered at ($-0.6,0.3$), ($0,0.5$), ($0.6,0.5$), ($-0.5,-0.4$), and ($0.5,-0.2$),
with polar angle $30^{\circ}$, $45^{\circ}$, $17^{\circ}$, $28.6^{\circ}$, and $9^{\circ}$, respectively.

\begin{figure} 
\centering 
\subfigure[$U = 4.0 \times 10^{-3}$ mm]{\label{fig:subfig:JMM_l75_u4}
\includegraphics[width=0.22\linewidth]{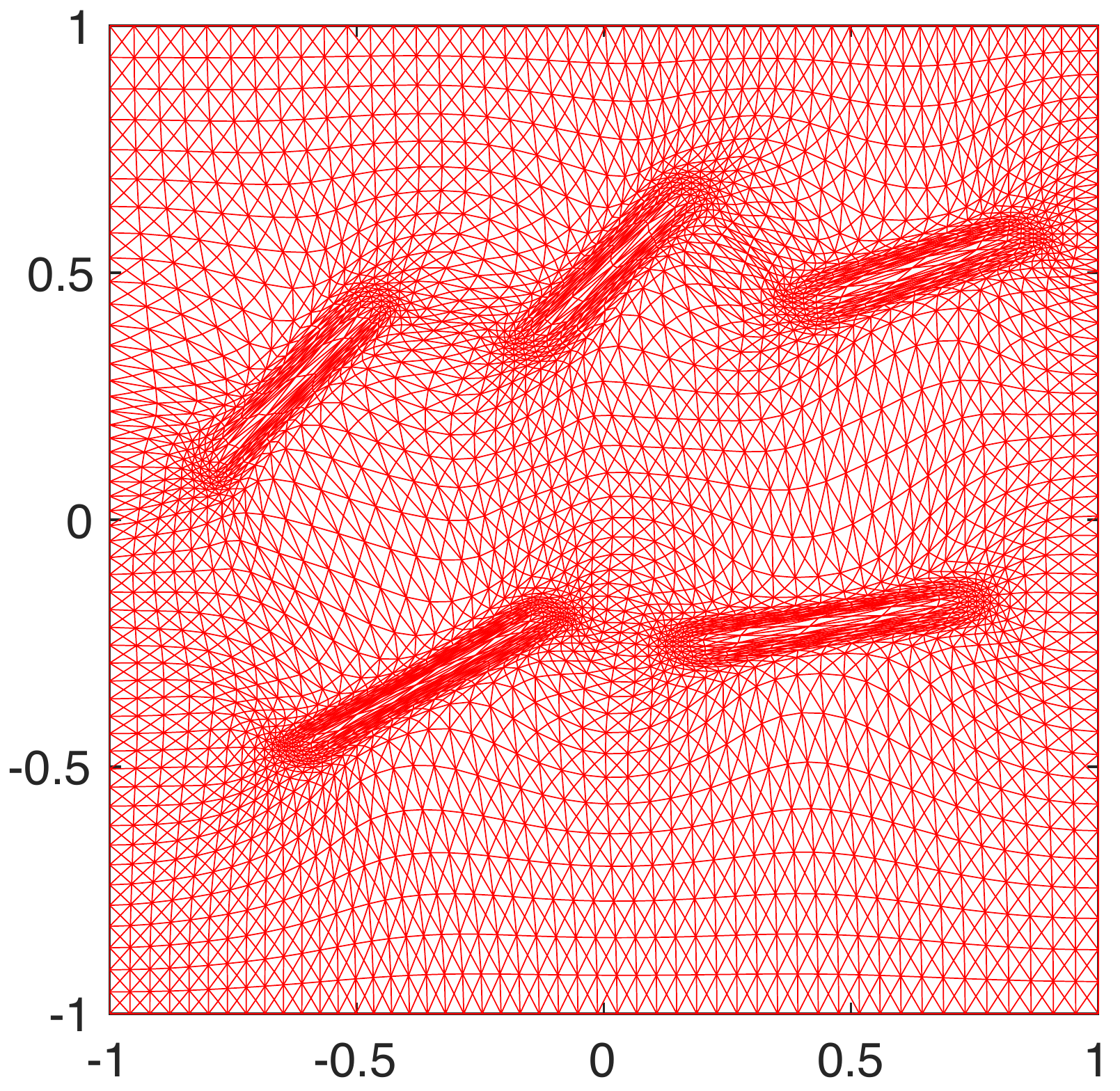}}
\subfigure[$U = 6.0 \times 10^{-3}$ mm]{\label{fig:subfig:JMM_l75_u6}
\includegraphics[width=0.22\linewidth]{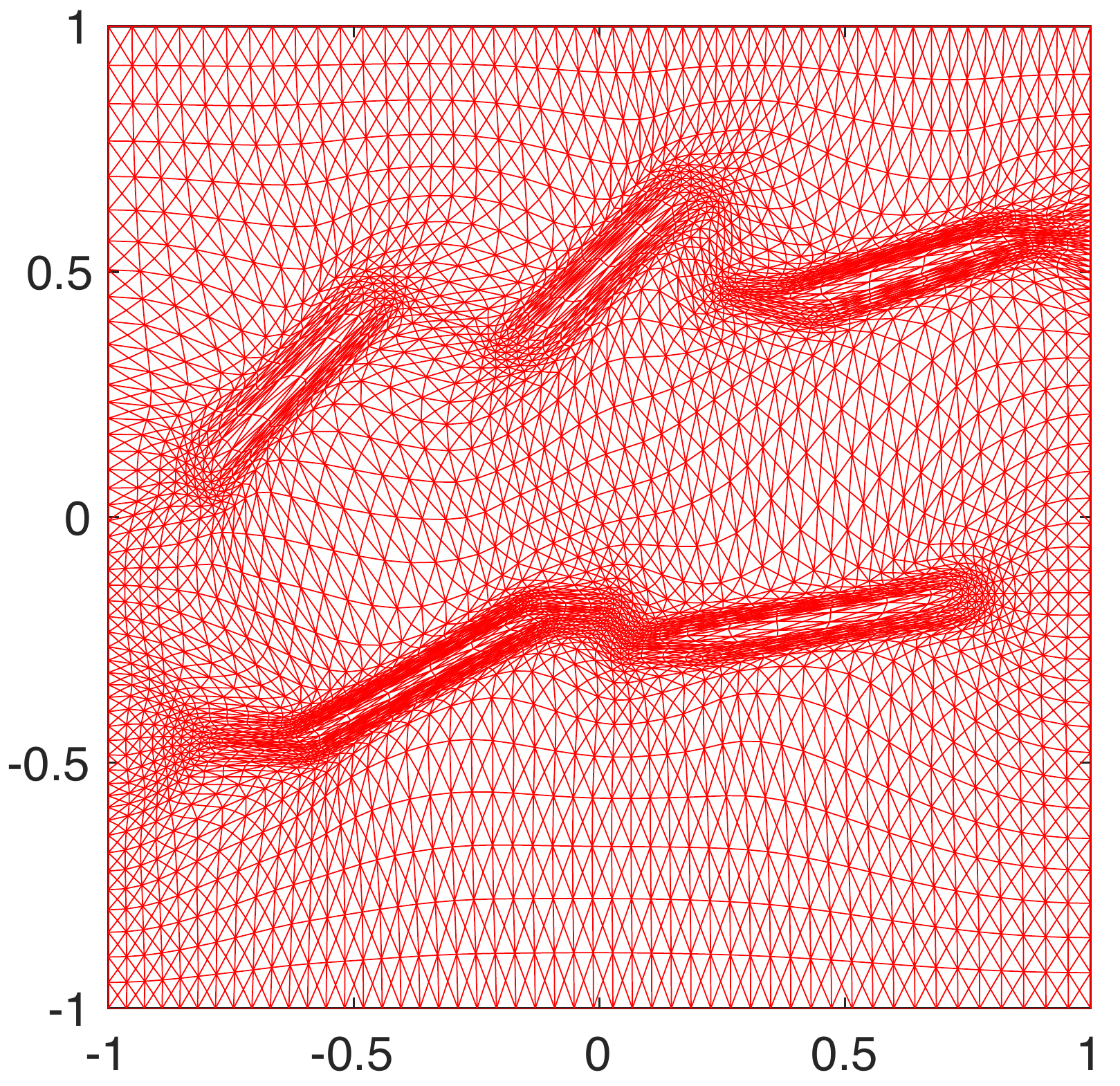}}
\subfigure[$U = 1.0 \times 10^{-2}$ mm]{\label{fig:subfig:JMM_l75_u10}
\includegraphics[width=0.22\linewidth]{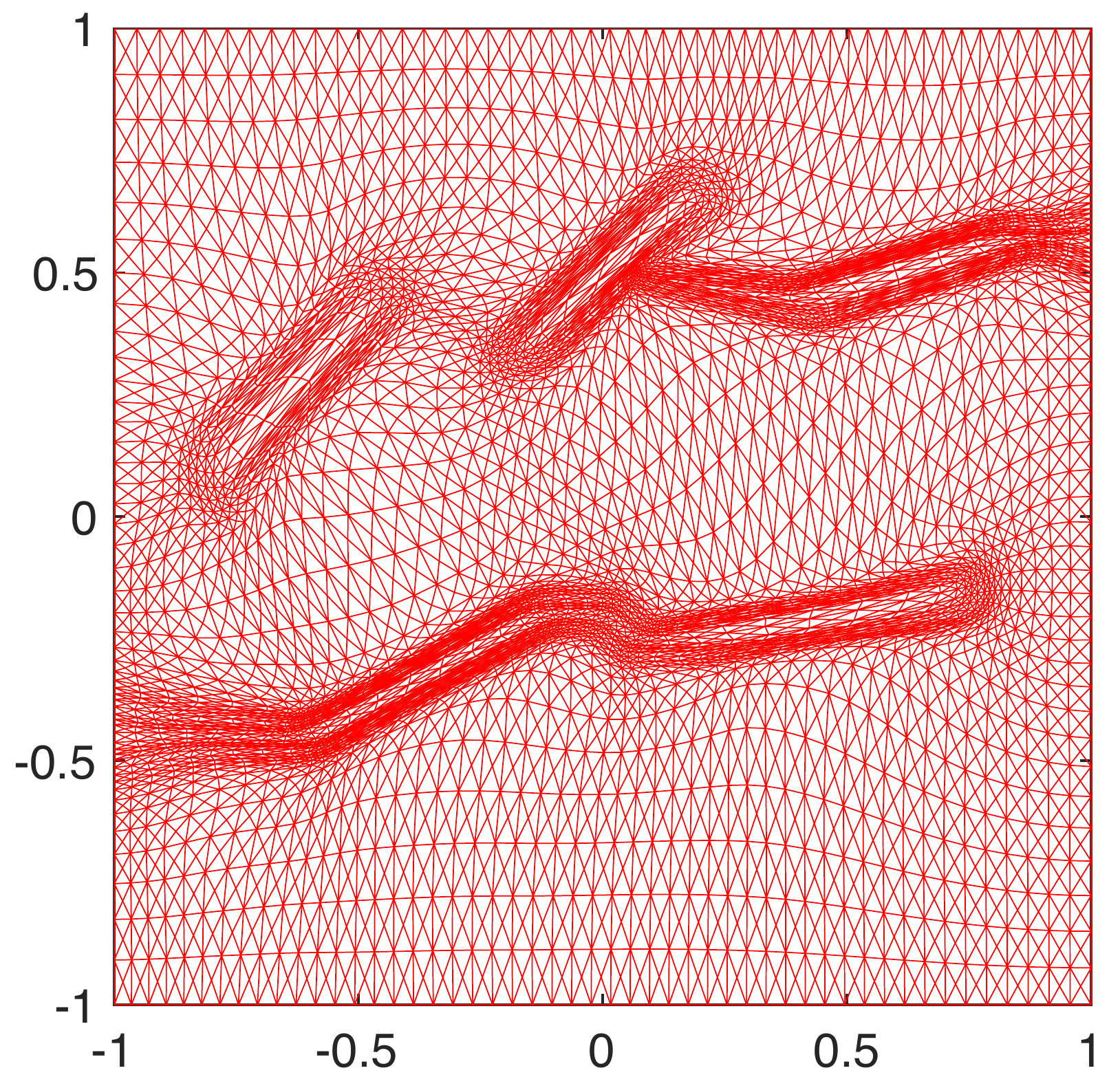}}
\subfigure[$U = 1.7 \times 10^{-2}$ mm]{\label{fig:subfig:JMM_l75_u17}
\includegraphics[width=0.22\linewidth]{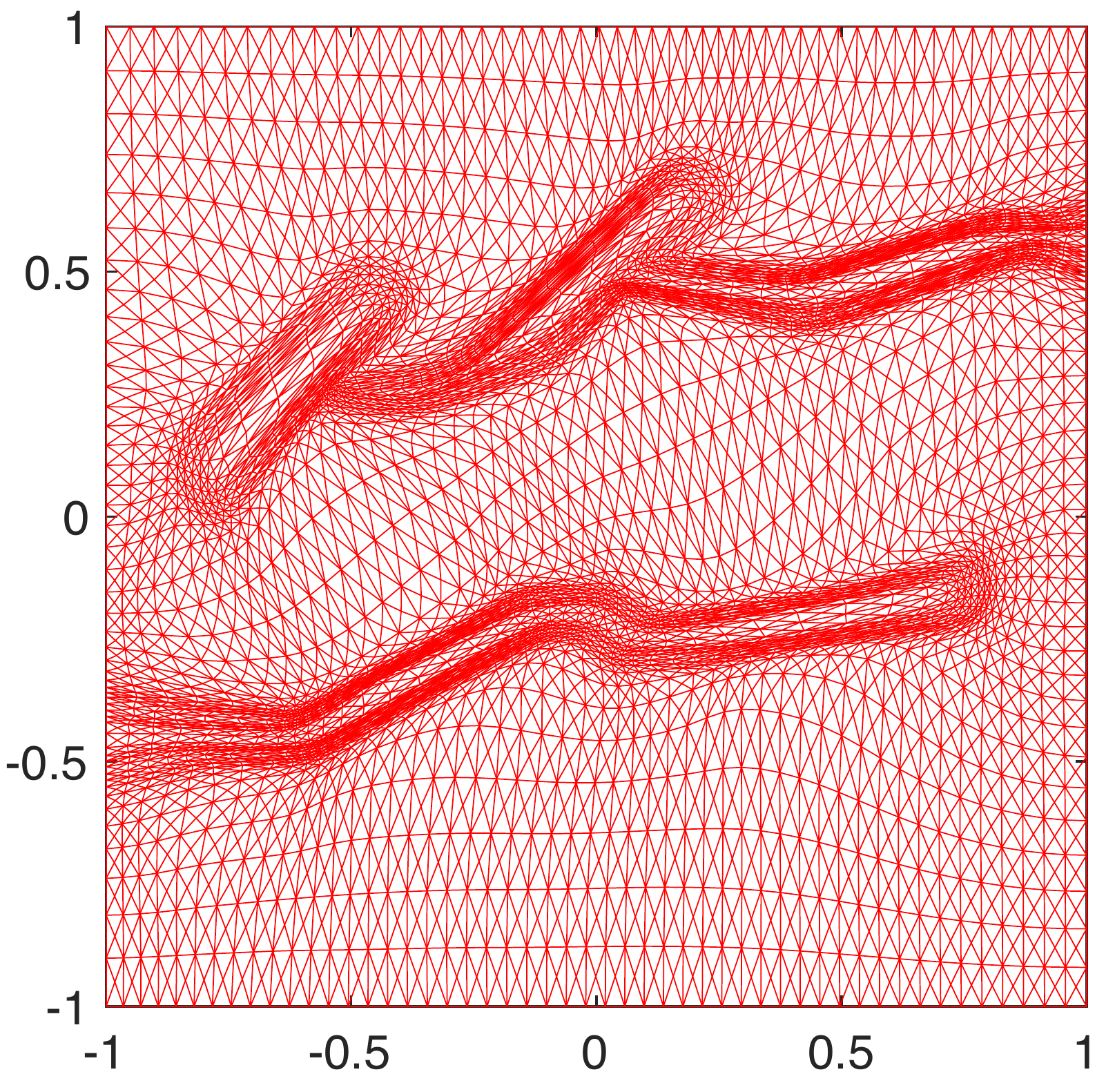}}
\vfill
\subfigure[$U = 4.0 \times 10^{-3}$ mm]{\label{fig:subfig:JMD_l75_u4}
\includegraphics[width=0.22\linewidth]{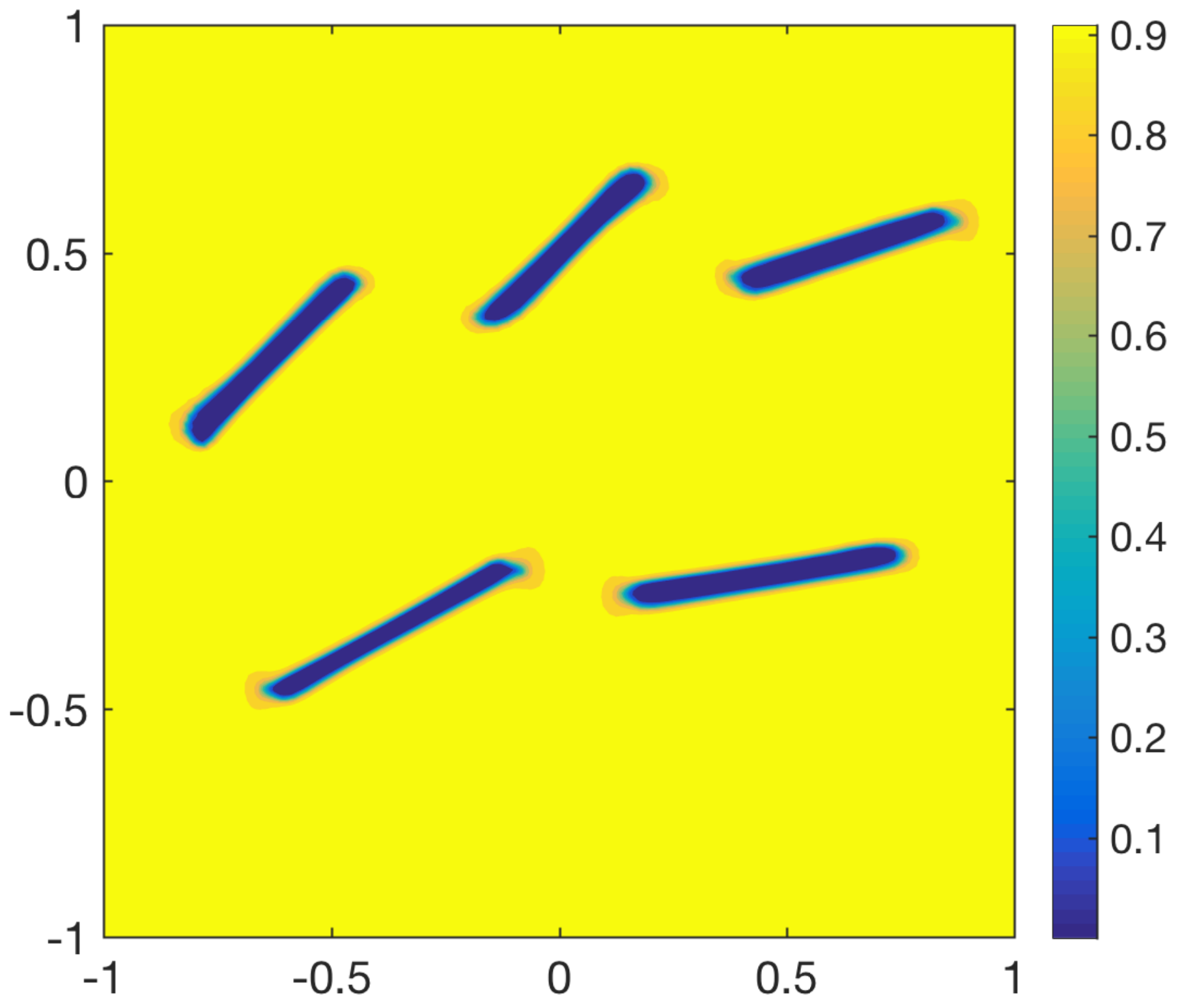}}
\subfigure[$U = 6.0 \times 10^{-3}$ mm]{\label{fig:subfig:JMD_l75_u6}
\includegraphics[width=0.22\linewidth]{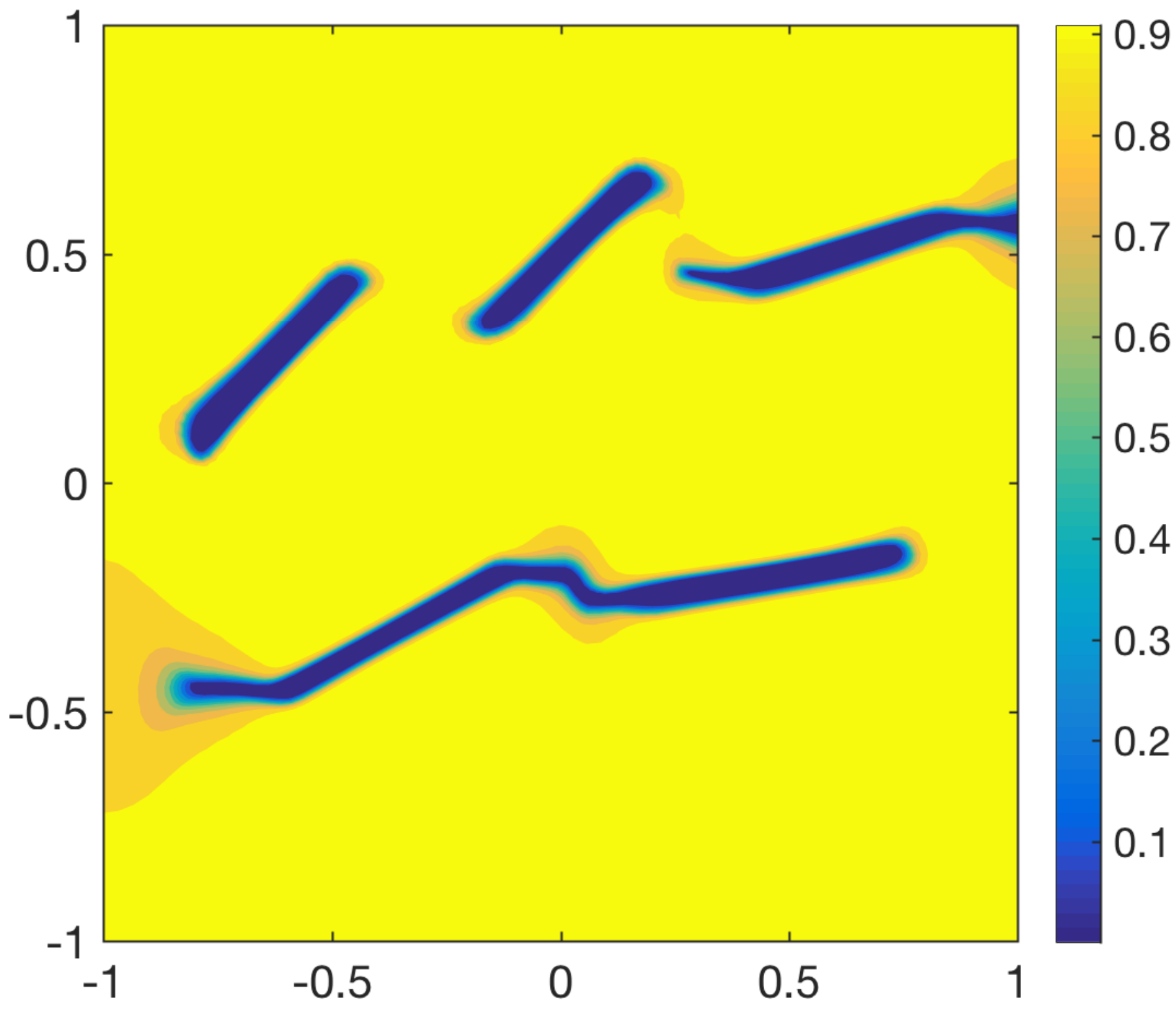}}
\subfigure[$U = 1.0 \times 10^{-2}$ mm]{\label{fig:subfig:JMD_l75_u10}
\includegraphics[width=0.22\linewidth]{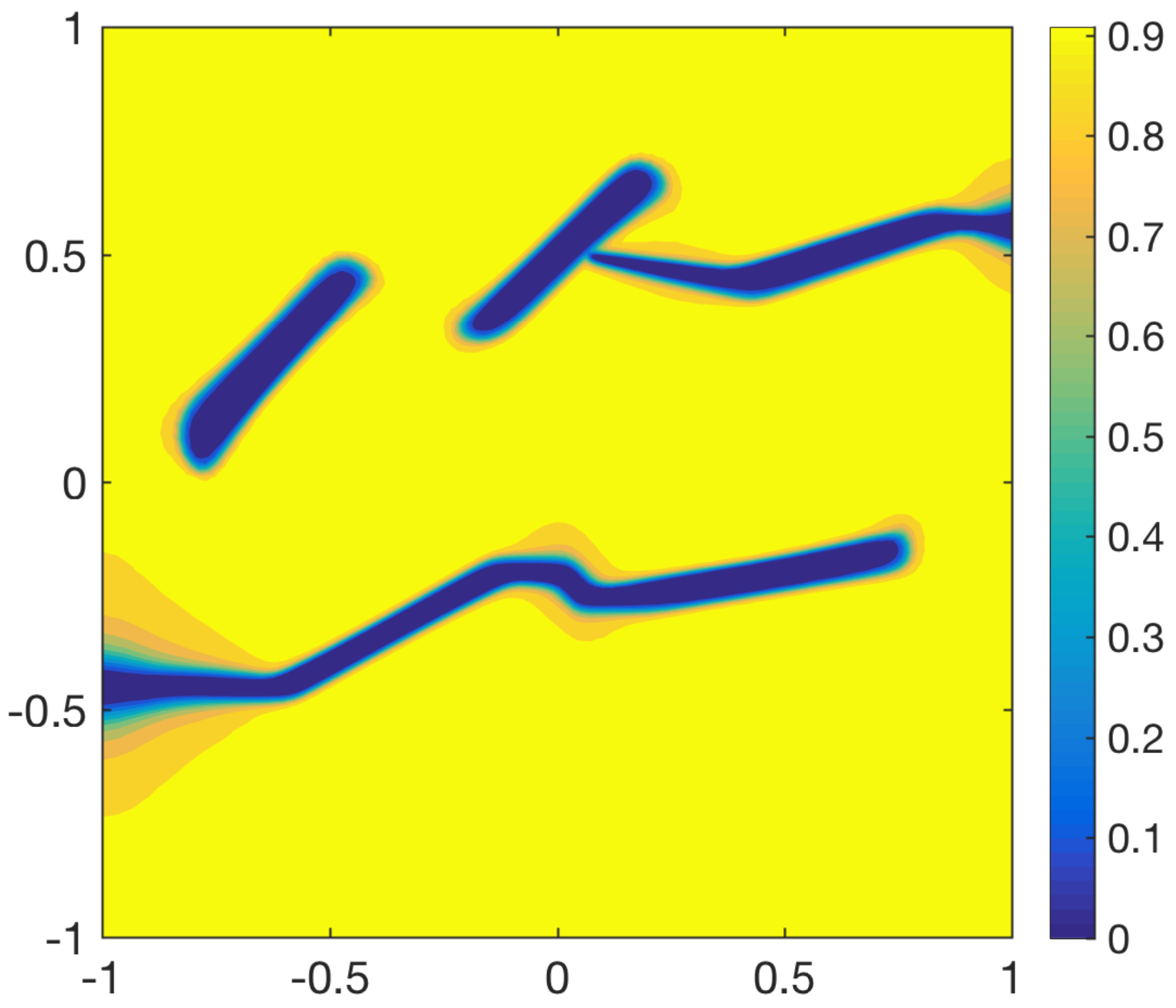}}
\subfigure[$U = 1.7 \times 10^{-2}$ mm]{\label{fig:subfig:JMD_l75_u17}
\includegraphics[width=0.22\linewidth]{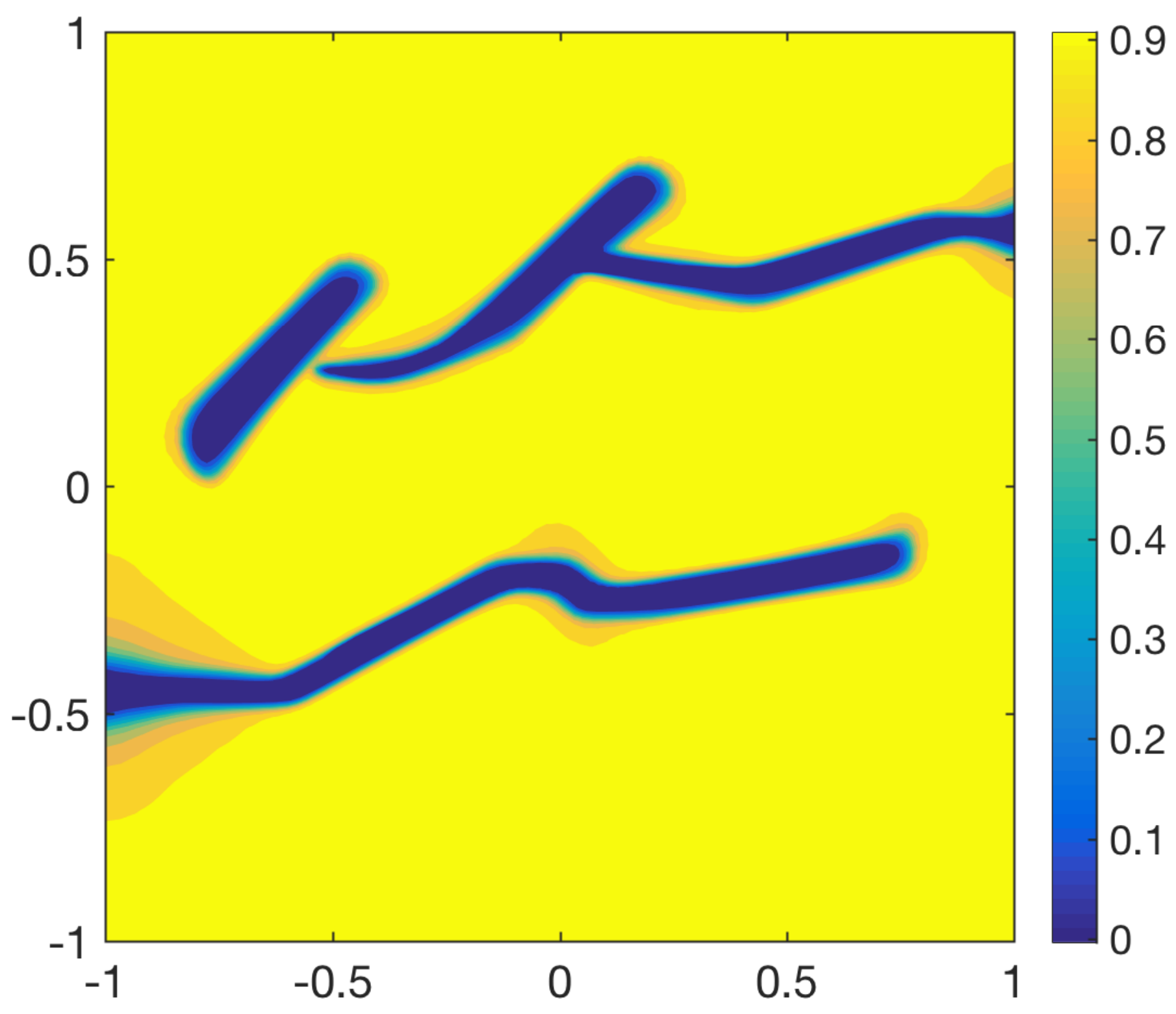}}
\caption{Example 3. The mesh and contours of the phase-field distribution during crack evolution
for the five-crack problem with $l = 0.00375$~mm and $N = 10,000$.}
\label{fig:l75_JunctionMult}
\end{figure}

Typical adaptive meshes and contours of the phase-field variable during crack evolution are shown
in Fig. \ref{fig:l75_JunctionMult}. The evolution can be described as follows.
\begin{itemize}
\item[(a)] $U = 4 \times 10^{-3}$~mm: the tips of Crack 3, 4 and 5 activate;
\item[(b)] $U = 6 \times 10^{-3}$~mm: Crack 4 connects to Crack 5, and Crack 3 has propagated
to the right edge of the plate;
\item[(c)] $U = 1 \times 10^{-2}$~mm: Crack 3 joins Crack 2, and Crack 4 has propagated to the left edge of the plate; 
\item[(d)] $U = 1.7 \times 10^{-2}$~mm: Crack 2 joins Crack 1.
\end{itemize}

Notice that the tips of Crack 3, 4, and 5 activate earlier due to their smaller polar angles.
Compared with other cracks, Crack 4 and 5 have longer length and a closer distance, 
which results in the early crack merging.
As the load increases, Crack 4 connects to Crack 5 and Crack 3 from both ends.
At later stages, Crack 2 begins to propagate after Crack 3 joining Crack 2.
Throughout the entire process, Crack 1 is limited to propagation due to its higher polar angle and
the interaction of other cracks.

\begin{figure} 
\centering 
\subfigure[$U = 1.0 \times 10^{-3}$ mm]{\label{fig:subfig:JTeM_l75_u10}
\includegraphics[width=0.22\linewidth]{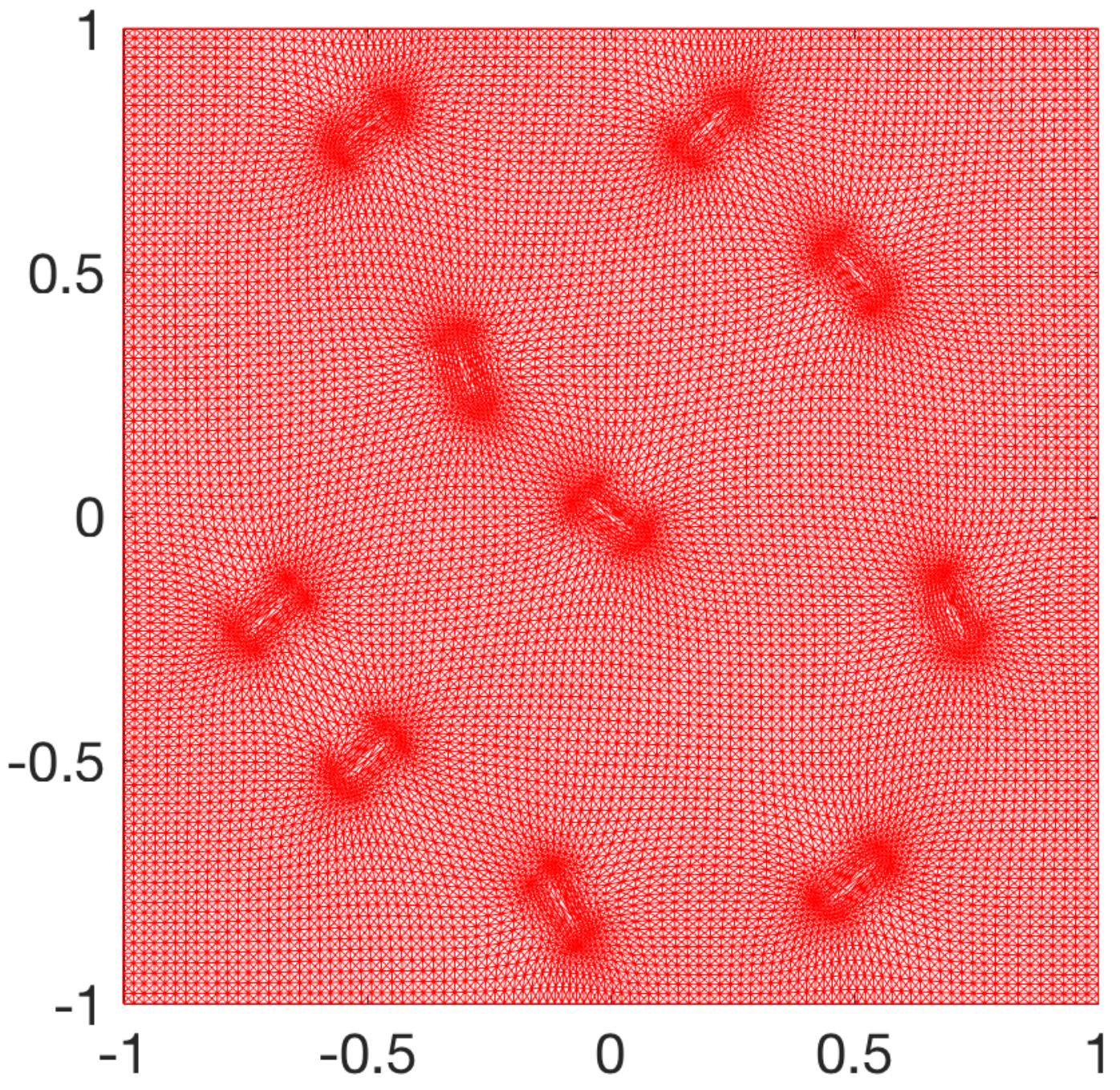}}
\subfigure[$U = 7.0 \times 10^{-3}$ mm]{\label{fig:subfig:JTeM_l75_u70}
\includegraphics[width=0.22\linewidth]{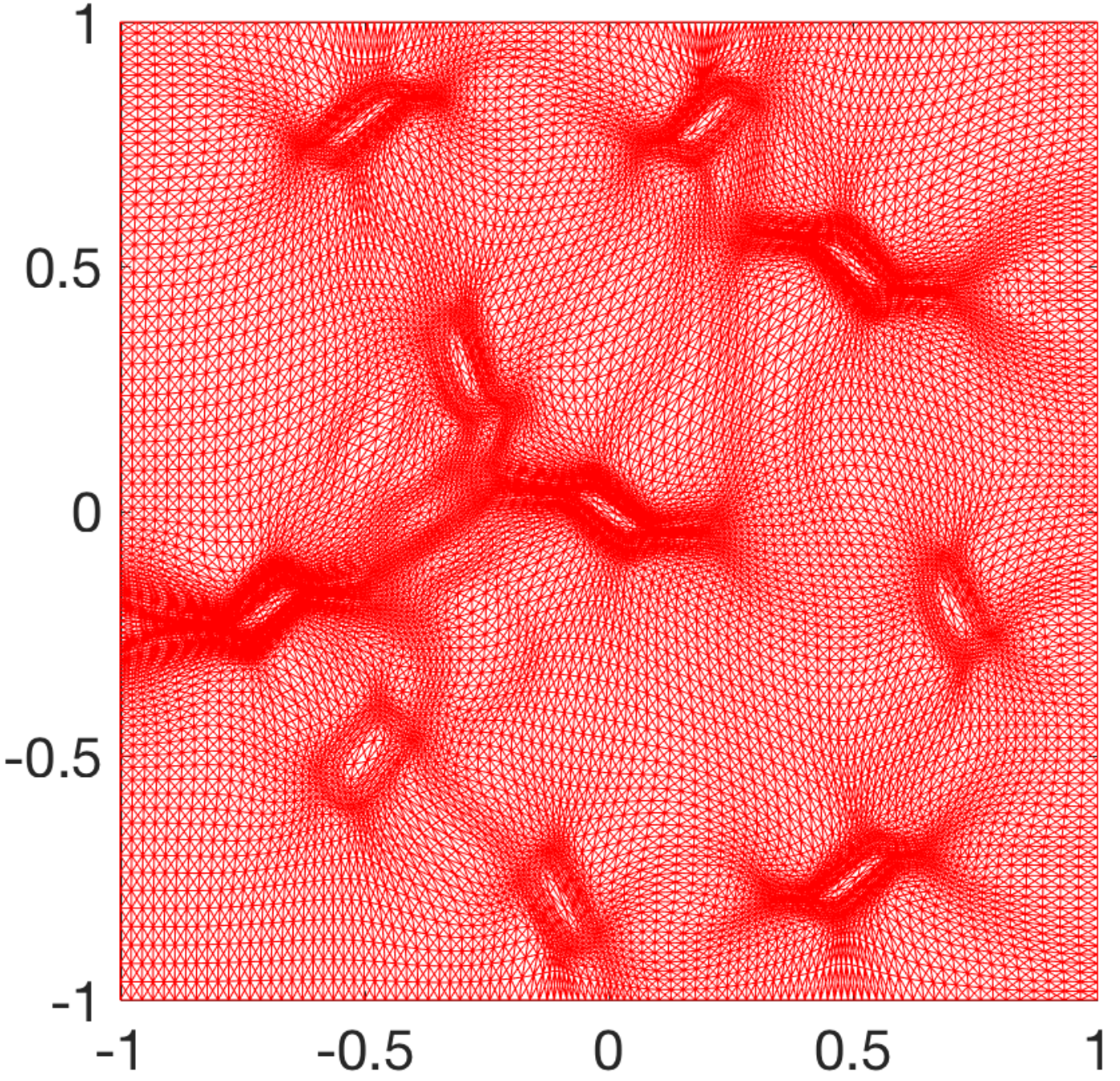}}
\subfigure[$U = 8.0 \times 10^{-3}$ mm]{\label{fig:subfig:JTeM_l75_u80}
\includegraphics[width=0.22\linewidth]{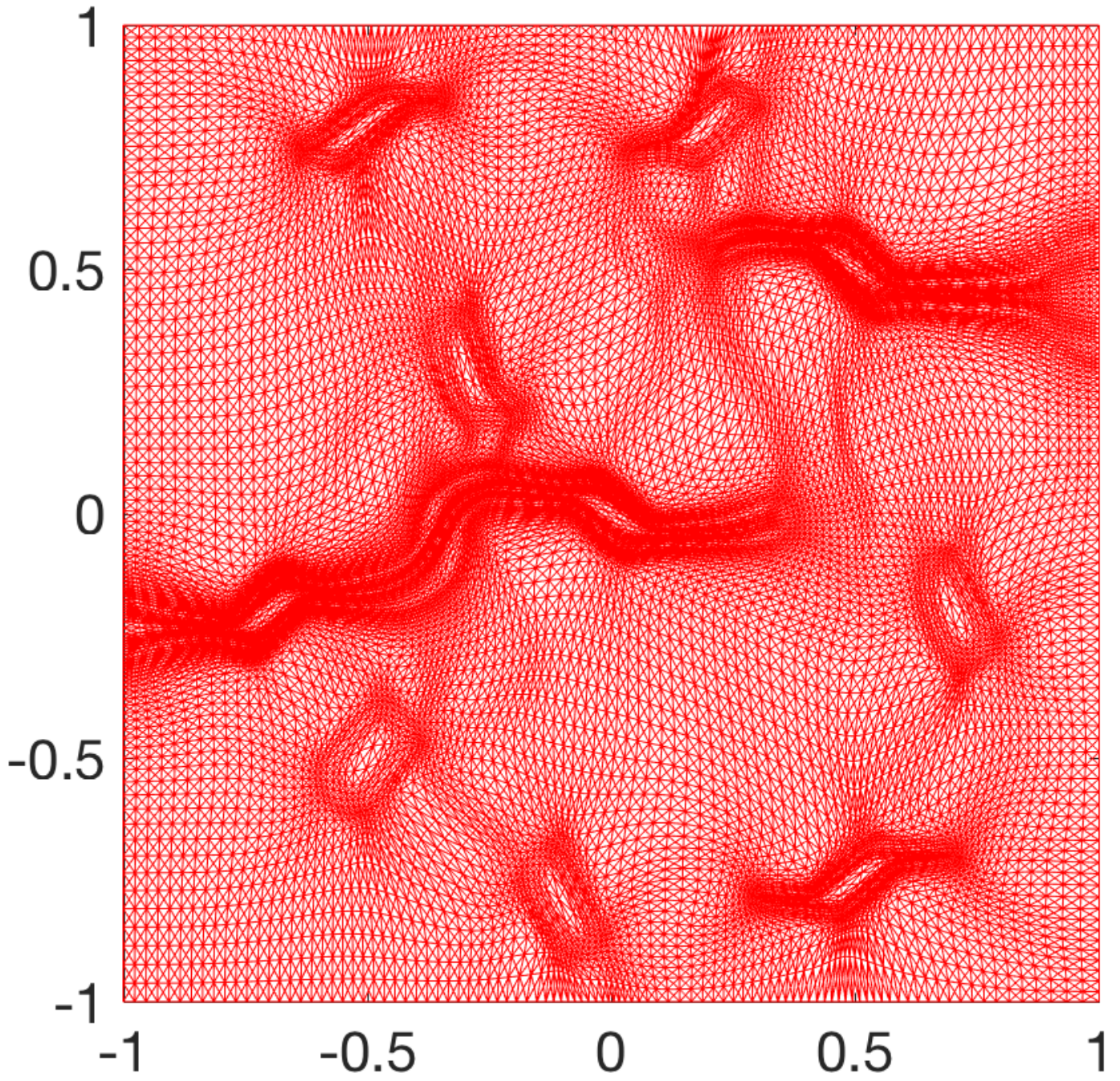}}
\subfigure[$U = 2.4 \times 10^{-2}$ mm]{\label{fig:subfig:JTeM_l75_u240}
\includegraphics[width=0.22\linewidth]{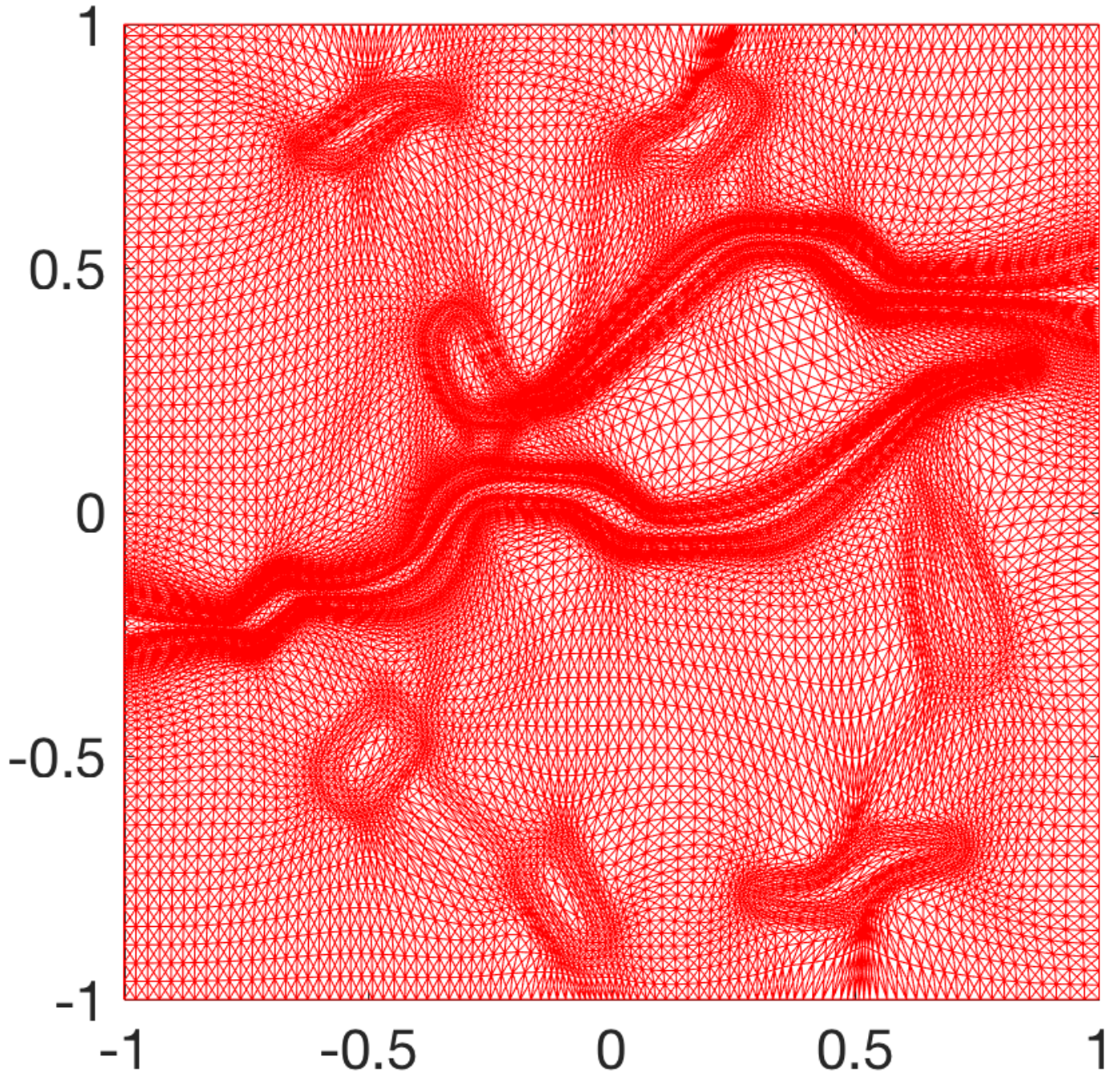}}
\vfill
\subfigure[$U = 1.0 \times 10^{-3}$ mm]{\label{fig:subfig:JTeD_l75_u10}
\includegraphics[width=0.22\linewidth]{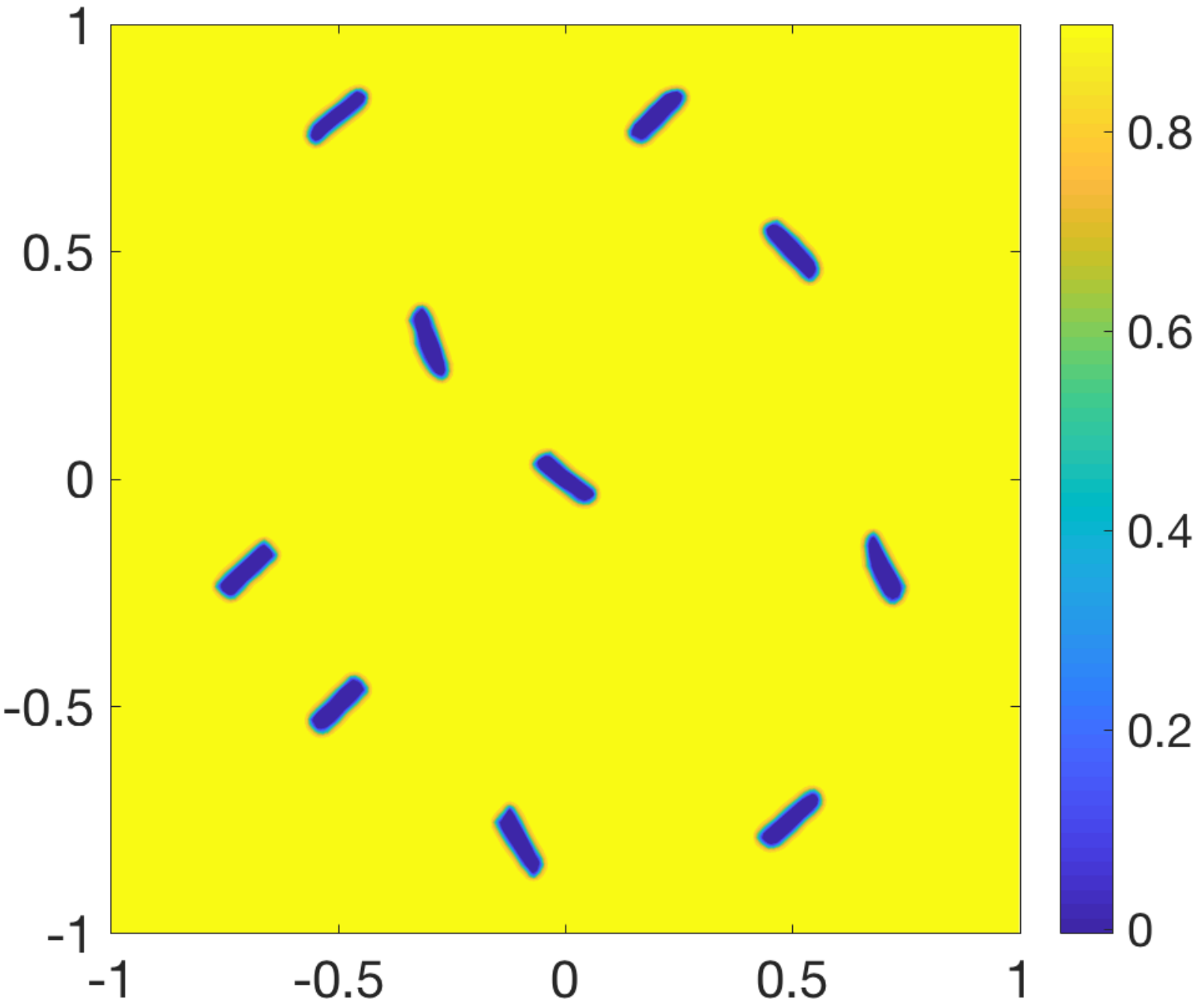}}
\subfigure[$U = 7.0 \times 10^{-3}$ mm]{\label{fig:subfig:JTeD_l75_u70}
\includegraphics[width=0.22\linewidth]{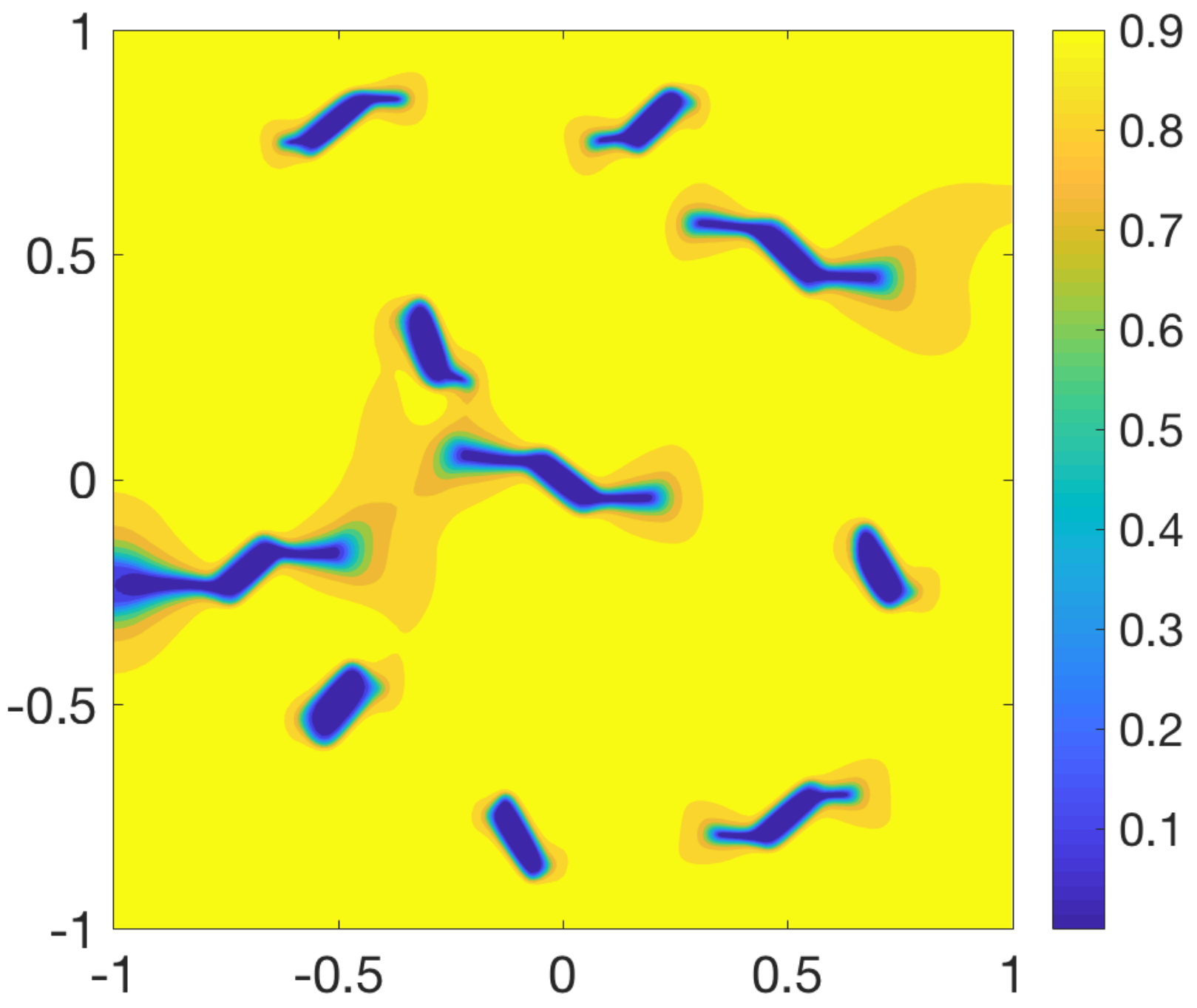}}
\subfigure[$U = 8.0 \times 10^{-3}$ mm]{\label{fig:subfig:JTeD_l75_u80}
\includegraphics[width=0.22\linewidth]{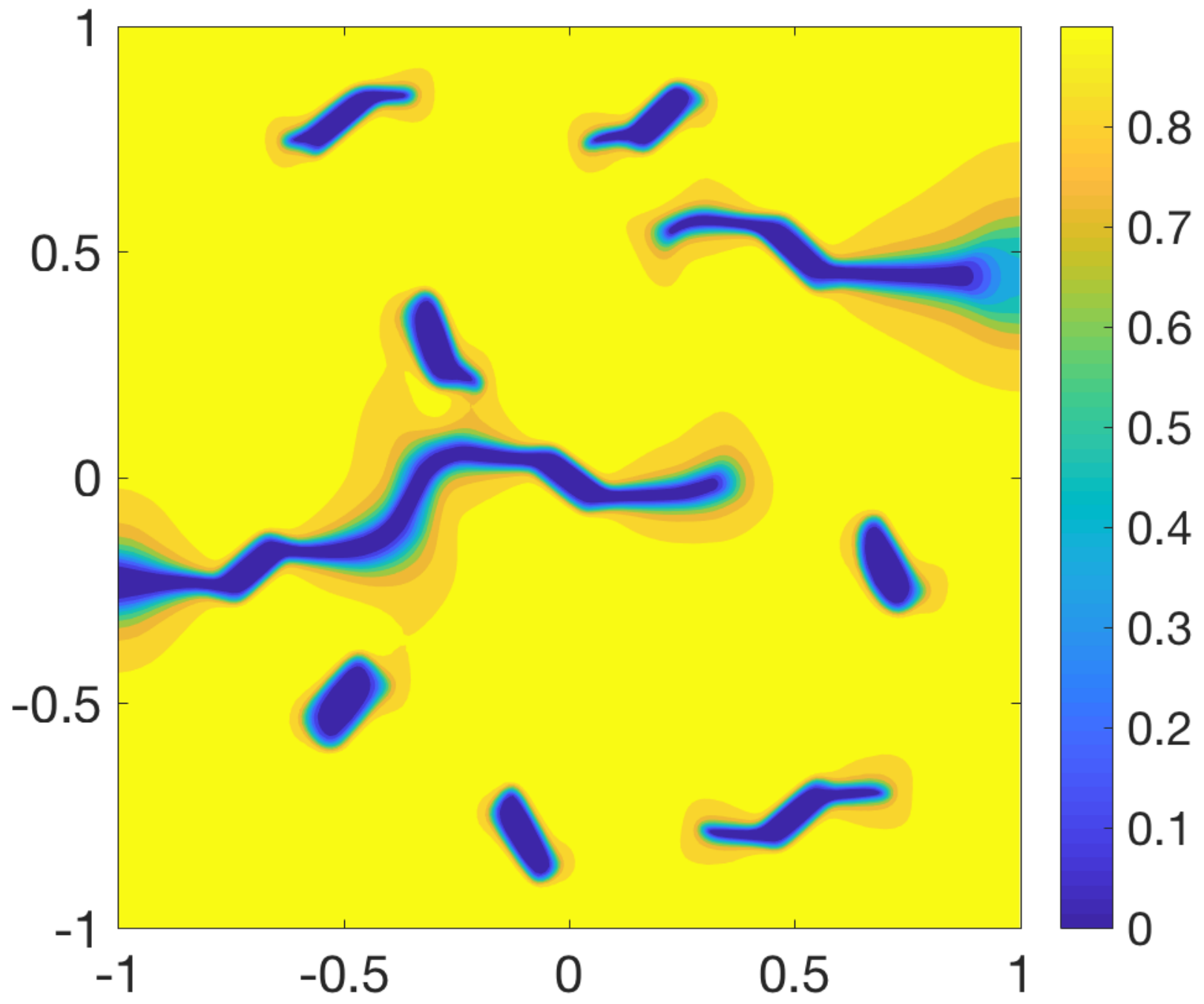}}
\subfigure[$U = 2.4 \times 10^{-2}$ mm]{\label{fig:subfig:JTeD_l75_u240}
\includegraphics[width=0.22\linewidth]{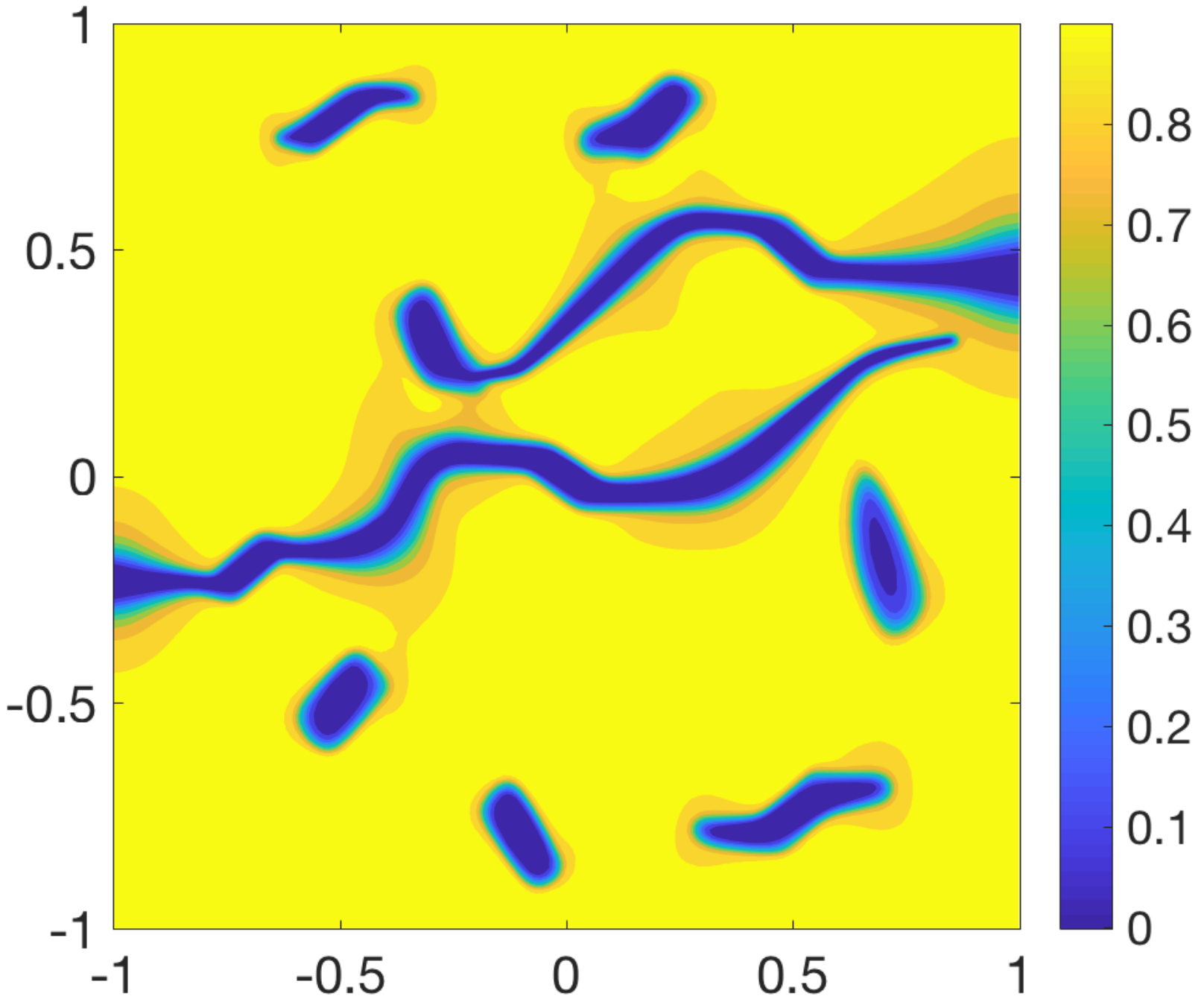}}
\caption{Example 3. The mesh and contours of the phase-field distribution during crack evolution
for the ten-crack problem with $l = 0.00375$~mm and $N = 40,000$.}
\label{fig:l75_JunctionTen}
\end{figure}

Finally, we consider a more complex situation with ten initial cracks
and with the same domain and boundary conditions as for the two-crack problem;
see Fig. \ref{fig:subfig:junctionTen}. The material parameters are the same
as in the previous examples except $g_c = 2.7 \times 10^{-4} $~kN/mm.
All of Crack 1 through Crack 10 have the same length $0.1$~mm.
They are centered at ($-0.5,0.8$), ($0.2,0.8$), ($-0.3,0.3$), ($0.5,0.5$), ($0,0$),
($-0.7,-0.2$), ($-0.5,-0.5$), ($-0.1,-0.8$), ($0.5,-0.75$), and ($0.7,-0.2$),
with polar angle $40^{\circ}$, $45^{\circ}$, $109^{\circ}$, $132^{\circ}$, $143^{\circ}$, $40^{\circ}$,
$45^{\circ}$, $120^{\circ}$, $40^{\circ}$, and $115^{\circ}$, respectively.
Typical adaptive meshes and contours of the phase-field variable during crack evolution are shown
in Fig. \ref{fig:l75_JunctionTen}. The evolution can be described as follows.
\begin{itemize}
\item[(a)] $U = 6 \times 10^{-3}$~mm: the tips of Crack 4, 5, 6, and 9 activate;
\item[(b)] $U = 8 \times 10^{-3}$~mm: Crack 5 connects to Crack 6, and Crack 5 has propagated
to the left edge of the plate;
\item[(c)] $U = 2.4 \times 10^{-2}$~mm: Crack 4 joins Crack 3, and Crack 4 has propagated to the right edge of the plate.
\end{itemize}

The results for both the two-, five- and ten-crack problems also demonstrate that
the MMPDE moving mesh method captures the crack propagation successfully for
multiple and complex crack systems.

\section{Conclusions}
\label{SEC:conclusion}

In the previous sections we have studied the moving mesh finite element solution of phase-field models
for brittle fracture and investigated regularization methods to improve the convergence of Newton's iteration
for solving the nonlinear system resulting from finite element discretization. 
The MMPDE moving mesh method has been used to track the crack propagation and
improve the computational efficiency. Numerical examples have been presented to demonstrate
the effectiveness of the moving mesh method to dynamically concentrate mesh elements around
propagating cracks and its ability to handle multiple and complex crack systems.

A distinguished feature in the numerical simulation of brittle fracture using phase-field modeling
is the decomposition of the strain tensor in the elastic energy which is necessary to account for
the reduction in the stiffness of an elastic solid due to the presence of cracks.
The decomposition makes the energy functional non-smooth and increases the nonlinearity of
the governing equation. Computationally, this causes the non-existence of the Jacobian matrix
of the nonlinear discrete equation \eqref{discrete-u} at some places and the failure of Newton's
iteration to converge (cf. Figs. \ref{fig:subfig:non_regularization} and \ref{fig:subfig:snon_regularization}).
Three methods, the sonic-point, exponential convolution, and smoothed 2-point convolution,
have been proposed in Section~\ref{SEC:convergence} to regularize the decomposition of the strain tensor.
Numerical examples have demonstrated that all of these methods can effectively improve the convergence
of Newton's iteration for relatively small values of the regularization parameter $\alpha$
but without compromising the accuracy of numerical simulation. 
Generally speaking, the larger $\alpha$ is, the faster Newton's iteration converges.
The sonic-point method works with smaller $\alpha$ than the other two methods
while the smoothed 2-point convolution method has slightly smaller effects on the load-deflection curves
for the same $\alpha$ than the other two methods.

Another benefit of the regularization is that the finite-difference approximation of the Jacobian matrix,
which is known to work for smooth problems, also works well in the current situation with regularization.


\vspace{20pt}

{\bf Acknowledgment.}
The authors would like to thank Dr. Junbo Cheng for bringing their attention to the sonic-point regularization
method in computational fluid dynamics. F.Z. was supported by China Scholarship Council (CSC)
and China University of Petroleum - Beijing (CUPB)  for his research visit to the University of Kansas
from September of 2015 to September of 2017. F.Z. is thankful to Department of Mathematics
of the University of Kansas for the hospitality during his visit.
The authors are grateful to the anonymous referees for their valuable comments in improving
the quality of the paper.


\end{document}